\documentclass{article}

\usepackage[english]{babel}
\usepackage[utf8]{inputenc}
\usepackage[T1]{fontenc}
\usepackage[a4paper,left=2cm,right=2cm,bottom=3cm]{geometry}
\usepackage{indentfirst}
\usepackage{amsmath,amssymb,amsfonts,amsthm,mathtools}
\numberwithin{equation}{section}
\usepackage{enumitem}
\usepackage{booktabs}
\usepackage{array}
\usepackage{color}
\definecolor{viridian}{rgb}{0.25,0.51,0.43}
\usepackage[colorlinks=true,linkcolor=cyan,citecolor=viridian]{hyperref}
\hypersetup{pdfcreator={},pdfproducer={LaTeX},hypertexnames=false,linktocpage=true}
\usepackage{prettyref}
\usepackage[nameinlink]{cleveref}
\usepackage{natbib}
\setcitestyle{numbers,square}
\allowdisplaybreaks[1]

\title{Lyapunov spectrum, Harder--Narasimhan filtration and thin monodromy}
\author{Fei Yu\thanks{School of Mathematical Sciences, Zhejiang University, Hangzhou 310058, China, \href{mailto:yufei@zju.edu.cn}{yufei@zju.edu.cn}}  \mbox{ and} Junming Zhang\thanks{Chern Institute of Mathematics and LPMC, Nankai University, Tianjin 300071, China, \href{mailto:junmingzhang@mail.nankai.edu.cn}{junmingzhang@mail.nankai.edu.cn}}}
\date{\vspace{-1em}}

\newtheorem{theorem}{Theorem}[section]
\newtheorem{corollary}[theorem]{Corollary}
\newtheorem{proposition}[theorem]{Proposition}
\newtheorem{lemma}[theorem]{Lemma}

\newtheorem{remark}[theorem]{Remark}
\theoremstyle{definition}
\newtheorem{example}[theorem]{Example}
\theoremstyle{plain}

\newtheorem{assumption}[theorem]{Assumption}

\newtheorem{convention}[theorem]{Convention}

\newrefformat{lem}{Lemma~\ref{#1}}
\newrefformat{coro}{Corollary~\ref{#1}}
\newrefformat{thm}{Theorem~\ref{#1}}
\newrefformat{prop}{Proposition~\ref{#1}}
\newrefformat{rem}{Remark~\ref{#1}}
\newrefformat{example}{Example~\ref{#1}}
\newrefformat{defn}{Definition~\ref{#1}}
\newrefformat{section}{Section~\ref{#1}}
\newrefformat{conj}{Conjecture~\ref{#1}}
\newrefformat{apdx}{Appendix~\ref{#1}}
\newrefformat{fact}{Fact~\ref{#1}}

\DeclareRobustCommand{\romitem}[1]{\textup{(\romannumeral#1\relax)}}

\DeclareRobustCommand{\upRoman}[1]{%
 \ifmmode
 \mathrm{\uppercase\expandafter{\romannumeral#1}}%
 \else
 \textup{\uppercase\expandafter{\romannumeral#1}}%
 \fi
}

\newcommand{\iu}{{\sqrt{-1}}}

\newcommand{\EE}{{\mathcal{E}}}

\newcommand{\KK}{{\mathcal{K}}}
\newcommand{\LL}{{\mathcal{L}}}

\begin{document}

\pagenumbering{gobble}
\maketitle
\pagenumbering{arabic}
\setcounter{section}{0}
\setcounter{page}{1}
\vspace{-1em}

\begin{abstract}
	We study the equality case of the inequalities between the Lyapunov spectrum and the Harder--Narasimhan slopes of a semi-simple local system with non-expanding cusp monodromy over a finite-area hyperbolic surface. This gives a nonlattice-transversality alternative of the monodromy group. On the other hand, we generalize Filip's Assumptions A and B to give classes of polarized variation of Hodge structures that achieve the equality. As an application, we recover the thinness of C-32 in $\mathrm{Sp}_6\mathbb{R}$. Besides, we provide an upper bound, called average Cartan spectrum, that dominates the Lyapunov spectrum and study its rigidity.
\end{abstract}

\tableofcontents

\section{Introduction}

\subsection{From Kontsevich--Zorich to thin monodromy}

Lyapunov exponents measure the asymptotic growth of flat sections, while Harder--Narasimhan slopes measure the instability of a holomorphic bundle. The starting point of comparing them is the Kontsevich--Zorich cocycle on the Hodge bundle over Teichm\"uller dynamics. Kontsevich's Hodge-theoretic formula and Zorich's deviation theory connected its Lyapunov spectrum with the geometry of the Hodge bundle \cite{Kontsevich97,Zorich97}. Forni proved the curvature formula for the sum of the nonnegative Kontsevich--Zorich exponents \cite{Forni02}. On a Teichm\"uller curve $X=\overline X\setminus D$, let $V^{1,0}_\bullet$ be the rank-$n$ Hodge bundle with its induced parabolic structure, and write $\KK_{\overline X}(D)=\Omega^1_{\overline X}(\log D)$ for the logarithmic cotangent bundle. Let $\lambda_1\geqslant\cdots\geqslant\lambda_n\geqslant0$ denote the nonnegative Lyapunov exponents of the cocycle.
Parabolic Chern--Weil theory turns the curvature formula into the normalized degree identity
\[\lambda_1+\cdots+\lambda_n=\frac{2\operatorname{pardeg}V^{1,0}_\bullet}{\deg\KK_{\overline X}(D)},\]
with the hyperbolic normalization fixed below. This identity connects Lyapunov data with algebro-geometric degrees of Hodge bundles.

Yu and Zuo subsequently made the relation with the finer Weierstrass and Harder--Narasimhan geometry explicit for Teichm\"uller curves \cite{YuZuo13}. Yu then formulated the polygon-comparison conjecture \cite{Yu18}: the Lyapunov polygon should lie on or above the normalized Harder--Narasimhan polygon. For a subbundle $\mathcal F_{i,\bullet}$ in the parabolic Harder--Narasimhan filtration of $V^{1,0}_\bullet$, with $r_i=\operatorname{rank}\mathcal F_i$, this predicts
\[\lambda_1+\cdots+\lambda_{r_i}\geqslant\frac{2\operatorname{pardeg}\mathcal F_{i,\bullet}}{\deg\KK_{\overline X}(D)}.\]

The connection between equality and thin monodromy was observed by Kontsevich in a 2013 lecture on K\"ahler random walks. For the fourteen hypergeometric Calabi--Yau variations of weight three, he noted numerical equality in the seven thin cases and strict inequality in the seven arithmetic cases \cite{Kontsevich13}. See Sarnak's Notes \cite{sarnak2014notes} for more discussions on thin matrix groups. Also see \cite{SinghVenkataramana14,BravThomas14} for more discussions on arithmeticity or thinness of hypergeometric monodromy groups in $\mathrm{Sp}_4\mathbb{R}$.

Eskin, Kontsevich, M\"oller, and Zorich proved a stronger inequality for arbitrary holomorphic subbundles of flat bundles over hyperbolic curves with non-expanding cusp monodromy, thereby proving Yu's conjecture \cite{EKMZ18}. They formulated Kontsevich's observation as a general rigidity problem: outside special geometric cases, should sharpness force monodromy to be thin in its ambient arithmetic group? Their paper also conjectured degree formulas for sums of Lyapunov exponents of hypergeometric Calabi--Yau variations of weight three.
Filip \cite{Filip21} proved the conjectured equality for the sum of Lyapunov exponents in the seven thin cases, using a Morse function method. 

On the other hand, Deroin--Dujardin \cite{deroin2017complex} used Garnett's theory of foliated harmonic measures and currents \cite{Garnett83} to study the relationship between the Lyapunov exponent and the degree of complex projective structures, which can be viewed as in the case of rank $2$ bundle over a curve. And soon Daniel--Deroin \cite{DanielDeroin19} developed this to a Brownian motion formalism on K\"ahler manifolds, generalized the comparison above to compact K\"ahler manifolds and gave a criterion of equality. Also see \cite{deroin2005hypersurfaces}.

In this paper, we study what happens when the equality holds for punctured Riemann surfaces. Also, we will give some classes of variation of Hodge structure that make the equality holds. They are generalized from Filip's Assumption A and Assumption B \cite{Filip21}. 

\subsection{An analytical degree interpretation of the gap}

Let $\overline{X}$ be a closed Riemann surface with a finite subset $D\subset\overline{X}$. Set $X=\overline{X}\setminus D$. We assume that $2g(\overline{X})-2+|D|>0$, where $g(\overline{X})$ denotes the genus of $\overline{X}$. Let $\KK_X$ be the canonical line bundle over $X$, and then $\KK_{\overline{X}}(D)$ is its logarithmic extension to $\overline{X}$. When the latter is regarded as a filtered line bundle with trivial filtration, we write $\KK_{\overline{X}}(D)_\bullet$. 

The Riemann surface $X$ then admits a unique complete conformal hyperbolic metric $g_{\mathrm{hyp}}$ with finite area. Throughout the paper we make the metric normalized with curvature $-4$, as in \cite{EKMZ18}. 
We work with a complex local system $\mathbb{V}$ over $X$ and the associated semisimple flat bundles $(V,\nabla)$, whose cusp monodromy has all eigenvalues on the unit circle. Then it is said to have \emph{non-expanding cusp monodromy}. $(V,\nabla)$ then admits the natural Deligne extension $\mathcal{V}_\bullet$ equipped with holomorphic structure $\nabla^{0,1}$. On the other hand, for such flat bundles, \cite{EKMZ18} proved the existence of Lyapunov spectrum $\boldsymbol{\lambda}(\mathbb{V})=(\lambda_j(\mathbb{V}))_{j=1}^{\operatorname{rank}(V)}$ of the corresponding cocycle lifting the unit-speed geodesic flow over the unit tangent bundle $T^1X$. 

\cite{EKMZ18} proved the following Lyapunov--HN inequality, i.e. for a rank-$k$ filtered holomorphic subbundle $\mathcal{F}_\bullet\subset\mathcal{V}_\bullet$, 
\begin{equation}\label{eq:LyaHN}
	\frac{2\operatorname{pardeg}\mathcal F_\bullet}{\deg\KK_{\overline{X}}(D)}\leqslant\lambda_1(\mathbb{V})+\cdots+\lambda_k(\mathbb{V}).
\end{equation}

For parabolic flat bundles, a Harder--Narasimhan filtration exists and is unique by \cite[Proposition~5.4]{Yokogawa93}. Applied the inequality above to the subbundles in the Harder--Narasimhan filtration, this gives dominance of the normalized Harder--Narasimhan slopes $2\boldsymbol{\mu}_{\mathrm{HN}}/\deg\KK_{\overline{X}}(D)$ by the Lyapunov spectrum $\boldsymbol{\lambda}$.

To study the equality, we use the Brownian motion formalism developed in \cite{DanielDeroin19}. The estimates of tame harmonic metric on $(V,\nabla)$ imply that Brownian motion with generator $\frac12\Delta_{\mathrm{hyp}}$ and the unit-speed geodesic flow have the same Lyapunov spectrum (see \prettyref{prop:lyap-coincidence}). And the maximal Lyapunov exponent can be illustrated as the Lyapunov exponent $\lambda(\nu)$ (see \eqref{eq:lambda-nu-def} for its definition) of a harmonic probability measure $\nu$ on the projectivization of the flat bundle. 

With this approach, we test the inequality for a saturated holomorphic line subbundle of the dual flat bundle, equipped with its induced filtration.

\begin{theorem}[Theorem~\ref{thm:main-A} \& Theorem~\ref{thm:strict-zero-divisor}]\label{thm:intro-A}
	Let $(V,\nabla)$ be a flat bundle as above and let $\mathcal{A}_\bullet\subset\mathcal{V}_\bullet^\vee$ be a filtered holomorphic line subbundle of its dual Deligne extension, with the induced filtration. Suppose that $(V,\nabla)$ is irreducible. Then for every harmonic probability measure $\nu$ on the flat projective bundle $p\colon\mathbb{P}(V)\to X$,
	\[\lambda(\nu)-\frac{2\operatorname{pardeg}\mathcal{A}_\bullet}{\deg\KK_{\overline{X}}(D)}=\delta_{\mathcal{A}}(\nu)\geqslant0,\]
	where $\delta_{\mathcal{A}}(\nu)$ is the mass of a regularized Poincar\'e--Lelong term, defined in \eqref{eq:PL-positive-term}. 
	
	Moreover, if $\delta_{\mathcal{A}}(\nu)$ vanishes, then \[\operatorname{supp}\nu\cap\Sigma_{\mathcal{A}}=\varnothing,\]
	where $\Sigma_{\mathcal{A}}$ denotes the zero divisor of the evaluation section, induced by $\mathcal{A}\to\mathcal{V}^\vee$, of $\mathcal{O}_{\mathbb{P}(V)}(1)\otimes p^*\mathcal{A}^\vee$ on $\mathbb{P}(V)$.
\end{theorem}
The term $\delta_{\mathcal{A}}(\nu)$ was already considered as a candidate on non-compact base in \cite[Section 3.5]{DanielDeroin19}, where it is called the \emph{analytic degree}. The analytic degree is defined using a tame harmonic metric $h$ on the local system $\mathbb{V}$. This gives decompositions
\[\nabla=D_h+\Phi+\Phi^{*_h},\]
where $D_h$ is $h$-unitary, and $\Phi$, called \emph{Higgs field}, is holomorphic $1$-form taking value in $\operatorname{End}(V)$ with respect to $D_h^{0,1}$.

As a corollary, we reprove \eqref{eq:LyaHN}, see Corollary~\ref{cor:HN-Lyap}. On a compact base and for a general connected semisimple Lie group $G_{\mathbb{R}}$, Proposition~\ref{prop:principal-HN-Lyap} expresses these inequalities intrinsically using the Harder--Narasimhan type of a flat $G_{\mathbb{R}}$-principal bundle.

One could also compare this mass with the dynamical degree in Daniel--Deroin \cite{DanielDeroin19} and the error term given by the bad locus in Costantini \cite{Costantini20}, in particular when the Riemann surface $X$ is compact.

\subsection{Nonlattice-transversality alternative}

The equality case of \eqref{eq:LyaHN} would imply a nonlattice-transversality alternative for the monodromy group. Here we assume that $G_{\mathbb R}$ is a connected semisimple real Lie group with finite center and without compact factors, and let $G_{\mathbb C}$ be its universal complexification. Write $\rho\colon\pi_1(X)\to G_{\mathbb R}$ for the monodromy and $\Gamma:=\rho(\pi_1(X))$ for its image. 

For an irreducible complex representation $W$ of $G_{\mathbb C}$, let $\mathcal{W}_\bullet$ be the Deligne extension of the associated local system $\mathbb{W}$, whose cusp monodromy is assumed to be non-expanding, as above.
A filtered holomorphic line subbundle $\mathcal A_\bullet\subset\mathcal W_\bullet^\vee$, with the induced filtration, defines by flat trivialization a holomorphic map
\[f_{\mathcal A}\colon\widetilde X\to\mathbb P(W^\vee).\]
For $S\subset\mathbb P(W)$, we write
\[S^*=\{[\alpha]\in\mathbb P(W^\vee):\mathbb P(\ker\alpha)\cap S=\varnothing\}\]
for its \emph{dual complement}. Thus $[\alpha]\in S^*$ precisely when $\alpha(v)\neq0$ for every $[v]\in S$. We call this condition transversality to $S$. We have the following nonlattice-transversality alternative without assuming Zariski density.

\begin{theorem}[Corollary~\ref{cor:lattice-transversality}]\label{thm:intro-C}
	Let $\mathcal S\subset\mathbb P(W)$ be a compact $G_{\mathbb R}$-orbit. If
	\[\lambda_1(\mathbb{W})=\frac{2\operatorname{pardeg}\mathcal A_\bullet}{\deg\KK_{\overline X}(D)},\]
	then either $\Gamma$ is not a lattice in $G_{\mathbb R}$, or
	\[f_{\mathcal A}(\widetilde X)\subset\mathcal S^*.\]
\end{theorem}

For Hermitian groups, this alternative restricts the image of a flag-valued map to the compact dual, as follows. Let $G_{\mathbb R}$ be a simple Hermitian real group with finite center and without compact factors, and let $\mathcal D\subset\check{\mathcal D}$ be its Borel-embedded symmetric domain. Choose a positive cominuscule $\omega$ and a positive integer $m$ such that there is an irreducible representation $G_{\mathbb{C}}\to V_{m\omega}$ with highest weight $m\omega$. 

\begin{theorem}[Theorem~\ref{thm:main-B}]\label{thm:intro-D}
	Put $W=V_{m\omega}^\vee$. If there is a filtered holomorphic line subbundle $\mathcal{A}_\bullet\subset\mathcal W_\bullet^\vee$ whose projective map is flag-valued, i.e.
	\[f_{\mathcal{A}}\colon\widetilde X\longrightarrow \check{\mathcal D},\]
	and
	\[\lambda_1(\mathbb{W})=\frac{2\operatorname{pardeg}\mathcal{A}_\bullet}{\deg\KK_{\overline{X}}(D)},\]
	then either $\Gamma$ is not a lattice in $G_{\mathbb R}$, or $f_\mathcal{A}$ takes values in a single Hermitian component: $\mathcal D$ when $G_\mathbb{R}$ is of non-tube type, and either $\mathcal D$ or the opposite domain $\mathcal D^{\mathrm{opp}}$ when $G_\mathbb{R}$ is of tube type.
\end{theorem}

For classical Hermitian groups, we discuss the realization of $f_\mathcal{A}$ when it takes value in a single Hermitian component in Section~\ref{subsec:classical-hermitian-models}. Furthermore, the transversality alternative could be excluded in lots of cases. See Theorem~\ref{thm:symplectic-contact} and Theorem~\ref{thm:unitary-contact} for detailed discussions.

\subsection{Generalization of Assumption A and Assumption B}

In \cite{Filip21}, Assumption A and Assumption B are introduced for polarized real variation of Hodge structure to obtain examples making \eqref{eq:LyaHN} an equality for some $k$. Here we generalize his assumptions and results following his strategy developed in Zhang \cite{Zhang25} and Li--Zhang \cite{LiZhang26}.

\begin{theorem}[Theorem~\ref{thm:main-C}]\label{thm:intro-E}
	Let $\mathbb V$ be a polarized complex weight-three variation with monodromy $\Gamma$ in $\mathrm{SU}_{m,n}$ and Hodge numbers $(m-1,1,1,n-1)$, where $m\geqslant n\geqslant2$.
	Put $F^2=V^{3,0}\oplus V^{2,1}$, of rank $m$. Suppose that $\Phi_{2,1}$ extends to an isomorphism of filtered line bundles on $\overline X$,
	\[\Phi_{2,1}\colon V^{2,1}_\bullet\xrightarrow{\sim}V^{1,2}_\bullet\otimes\KK_{\overline X}(D),\]
	where $V^{p,q}_\bullet$ denotes the metric extension of the Hodge bundle. Then
	\[\lambda_1(\mathbb V)+\cdots+\lambda_m(\mathbb V)=\frac{2\operatorname{pardeg}F^2}{\deg\KK_{\overline{X}}(D)}.\]
	There is also a proper closed set $\Lambda$, which is  $\Gamma$-invariant, of maximal isotropic complex $n$-planes for the flat Hermitian polarization such that $F^2(z)\cap L=0$ for every $L\in\Lambda$ and $z\in\widetilde X$.
	
	Furthermore, the monodromy $\Gamma$ is not a lattice in $\mathrm{SU}_{m,n}$. If $m=n$ and the variation is the complexification of a real polarized variation, then $\Gamma<\mathrm{Sp}_{2n}\mathbb R$ is not a lattice in $\mathrm{Sp}_{2n}\mathbb R$ either. In either case, if $\Gamma$ is contained in an arithmetic lattice $\Gamma_{\mathrm{ar}}$ of the ambient group, then $[\Gamma_{\mathrm{ar}}:\Gamma]=\infty$. If $\Gamma$ is also Zariski dense in that group, it is thin.
\end{theorem}

The other ordering of $m,n$ is covered by duality. When $m=n=2$ and the variation arises from a polarized real variation of Hodge structure, the assumption on $\Phi_{2,1}$ is exactly Filip's Assumption A \cite[Definition~2.1.9]{Filip21}. This theorem generalizes it to complex variation and higher rank.

\begin{remark}
	The Anosov property of such monodromy is studied by Bronstein--Davalo--Li \cite{BDL26}.
\end{remark}

The proof also yields a domain of discontinuity in this flag manifold, and the harmonic map into $\mathrm{SU}_{m,n}/\mathrm S(\mathrm U_m\times\mathrm U_n)$ is a bi-Lipschitz embedding, see Corollary~\ref{cor:weight-three-geometry}. When $X$ is compact, the domain has compact quotient and the word-metric orbit map is quasi-isometric. Similar domain of discontinuity and bi-Lipschitz embedding for $4$-cyclic Higgs bundles of rank $(1,1,1,n-1)$ were built in \cite{LiZhang26}.

An orthogonal analogue generalizes the rank-five setting of Assumption~B in \cite[Definition~5.4.1]{Filip21}.

\begin{theorem}[Corollary~\ref{cor:orthogonal-B-thinness}, Proposition~\ref{prop:orthogonal-B-equality} \& Proposition~\ref{prop:orthogonal-B-strictness}]\label{thm:intro-orthogonal-B}
	Let $\mathbb V$ be a polarized real weight-four variation with monodromy $\Gamma$ and Hodge numbers $(n,1,1,1,n)$, where $n\geqslant1$. The signature of its polarization $Q$ is $(2n+1,2)$. Suppose that the monodromy $\Gamma$ lies in $\mathrm{SO}(Q)$. Put $F^3=V^{4,0}\oplus V^{3,1}$, of rank $n+1$. Suppose that the logarithmic extension
	\[\Phi_{2,2}\colon V^{2,2}_\bullet\xrightarrow{\sim}V^{1,3}_\bullet\otimes\KK_{\overline X}(D)\]
	is an isomorphism of filtered line bundles on $\overline X$, where $V^{p,q}_\bullet$ denotes the metric extension of the Hodge bundle. Then $\Gamma$ is not a lattice in $\mathrm{SO}(Q)$, and
	\[\sum_{j=1}^{n+1}\lambda_j(\mathbb V)\geqslant\frac{2\operatorname{pardeg}F^3}{\deg\KK_{\overline X}(D)}.\]
	This is an equality when $n=1$, and it is strict when $n>1$ and $\Gamma$ is Zariski dense in $\mathrm{SO}(Q)$.
	If $\Gamma$ is contained in an arithmetic lattice $\Gamma_{\mathrm{ar}}$ of $\mathrm{SO}(Q)$, then $[\Gamma_{\mathrm{ar}}:\Gamma]=\infty$; if it is also Zariski dense, it is thin.
\end{theorem}

Note that the equality between the sum of Lyapunov exponents and the normalized degree when $n=1$ has been claimed in \cite[Section 5.4.4]{Filip21}. Here we will give a rigorous proof on that. When $n>1$, the strict inequality implies a possibility of negation to the general equivalence between thin monodromy and the Lyapunov--HN equality.  

We give two applications of above theorems on hypergeometric variations.

\paragraph{Hypergeometric Calabi--Yau variations}
Recall that Filip \cite{Filip21} proved the conjectured Lyapunov--HN equality \cite{Kontsevich13,EKMZ18} for seven thin cases of hypergeometric Calabi--Yau variations of weight three. For seven arithmetic cases, we give a proof for the numerical strict inequality \cite{Kontsevich13} using Theorem~\ref{thm:symplectic-contact}. See Corollary~\ref{cor:fourteen-CY} for detailed discussions. 

\paragraph{C-32 in $\mathrm{Sp}_6\mathbb{R}$}
The degree-six symplectic hypergeometric groups were studied by Bajpai--Dona--Singh--Singh \cite{BDSS21}, motivated by a question of Katz about whether the maximally unipotent degree-six family follows a similar pattern as degree-four symplectic hypergeometric groups. Bajpai--Dona--Nitsche \cite{BDN25} later settled many further degree-six cases, including many thin cases, by computer-assisted ping-pong. In \cite[Table 3]{BDN25}, they list the three remaining degree-six symplectic cases, C-32, C-47, and C-55, for which neither arithmeticity nor thinness was known. In May 2026, Riestenberg--Taha--Trettel--Wienhard, in an earlier version of \cite{RTTW26}, proved the arithmeticity of C-47 and C-55. They also conjectured the thinness of C-32. Then in October 2026, they proved this in the newest version of \cite{RTTW26} by computer-assisted ping-pong. In this paper, we use Theorem~\ref{thm:intro-E} to recover the thinness of C-32. Also, we obtain a Lyapunov--HN equality of the corresponding local system. See Section~\ref{subsec:C32-example} for detailed discussions.

\subsection{Average Cartan spectrum}

Besides the lower bound of Lyapunov spectrum given by the Harder--Narasimhan slopes, we also give an upper bound $\mathbf{C}$, called \emph{average Cartan spectrum}, that dominates the Lyapunov spectrum $\boldsymbol{\lambda}$. 

Let $G_{\mathbb R}$ be a connected semisimple real Lie group with finite center and without compact factors. We use the adapted tame harmonic metric $h$ on a $G_{\mathbb{R}}$-local system $\mathbb{V}$ with non-expanding cusp monodromy. Then we have the monodromy equivariant harmonic map $u\colon\widetilde{X}\to G_{\mathbb{R}}/K_{\mathbb{R}}$, where $K_{\mathbb{R}}$ denotes a fixed maximal compact subgroup of $G_{\mathbb{R}}$. And $\Phi+\Phi^{*_h}$, where $\Phi$ denotes the corresponding Higgs field, can be interpreted as the differential of $u$. 

Let $\kappa$ be the Cartan projection of $G_\mathbb{R}$ mapping to the closed positive Weyl chamber. For every $v\in T^1X$, define
\[H(v):=\kappa\!\left(\exp(-\mathrm{d}u(v))\right).\]
Here $\mathrm{d}u(v)$ is written in any local frame of the $K_{\mathbb R}$-reduction $h$. The definition is independent of this choice because $\kappa$ is invariant under $K_{\mathbb R}$-conjugation. Equivalently, $H(v)$ is the unique dominant element in the $K_{\mathbb R}$-orbit of $-\mathrm{d}u(v)$. The \emph{average Cartan spectrum} is
\[\mathbf{C}=\int_{T^1X}H(v)\,\mathrm{d}\mu_L(v),\] 
where $\mu_L$ is the normalized Liouville measure.

\begin{theorem}[Theorem~\ref{thm:main-D}]\label{thm:intro-F}
	Let $\boldsymbol{\lambda}$ be the Lyapunov spectrum of the $G_{\mathbb R}$-cocycle (see Section~\ref{subsec:Lya}). Then $\boldsymbol{\lambda}\preceq\mathbf{C}$.
	Additionally, if the Higgs field $\Phi$ is nilpotent, as it is for a variation of Hodge structure, then
	\[\boldsymbol{\lambda}=\mathbf{C}\quad\Longleftrightarrow\quad u\text{ is totally geodesic}.\]
	Let $\eta\colon G_{\mathbb R}\to G_{\mathbb C}$ be the universal complexification map, and put $G_{\mathbb R}^{\mathrm{lin}}:=\eta(G_{\mathbb R})$ and $\rho^{\mathrm{lin}}:=\eta\circ\rho$.
	If these equivalent conditions hold and $u$ is nonconstant, then
	\[\rho^{\mathrm{lin}}(\gamma)=\iota(j(\gamma))\rho_c(\gamma),\]
	where $j\colon\pi_1(X)\to\mathrm{SL}_2\mathbb R$ lifts the uniformizing representation, $\iota\colon\mathrm{SL}_2\mathbb R\to G_{\mathbb R}^{\mathrm{lin}}$ is a Lie-group homomorphism defining the totally geodesic disk, and $\rho_c$ takes values in its compact pointwise stabilizer, which centralizes $\iota(\mathrm{SL}_2\mathbb R)$.
\end{theorem}

We also classify the tame harmonic bundle underlies a polarized variation of Hodge structure with nonconstant period map and $\boldsymbol{\lambda}=\mathbf{C}$ in Section~\ref{subsec:hodge-reconstruction}. For maximal Hodge bundles over closed Riemann surfaces, which are studied by Biquard--Collier--Garc\'ia-Prada--Toledo \cite{BCGT23}, their rigidity theorem supplies the horizontal uniformizing factor and compact factor. See Section~\ref{subsec:maximal-hodge-comparison} for more discussions.

\begin{remark}
	Another candidate of upper bound is $\boldsymbol{\varepsilon}$, called \emph{eigenvalues of curvature}, introduced by Forni \cite{Forni02} to study Lyapunov exponents of the Hodge bundle. For the weight-one Hodge bundle over a Teichm\"uller curve, it is proven by Yu \cite{Yu18} that $2\boldsymbol{\mu}_{\mathrm{HN}}/\deg(\KK_{\overline{X}}(D))\preceq\boldsymbol{\varepsilon}$ and $\boldsymbol{\lambda}$ is bounded by $2\boldsymbol{\varepsilon}$ from above. One could verify that $\boldsymbol{\varepsilon}\preceq\mathbf{C}$, so it is still unknown that whether $\boldsymbol{\lambda}\preceq\boldsymbol{\varepsilon}$ or not.
\end{remark}

\paragraph{Organization}
In Section~\ref{sec:setup}, we recall the tame harmonic metric on the local system and the definition of Lyapunov spectra for the geodesic flow and the Brownian motion. We also proved the existence of Lyapunov-realizing harmonic measure. Then, we give the analytical degree interpretation of the gap between Lyapunov spectrum and Harder--Narasimhan slopes in Section~\ref{sec:PL-identity}. Section~\ref{sec:hermitian-transversality} establishes the nonlattice-transversality alternative in general and its specialization for Hermitian groups. We then generalize Filip's Assumptions A and B in Section~\ref{sec:CY-applications}, derive the geometric consequences and symplectic hypergeometric applications. Finally, we study the average Cartan spectrum in Section~\ref{sec:Lyapunov-Cartan}.

\paragraph{Declaration of AI use}
The idea of Theorem~\ref{thm:main-A}, Theorem~\ref{thm:strict-zero-divisor} and Section~\ref{sec:Lyapunov-Cartan} appeared to the authors from an iterative, guided conversation with ChatGPT 5.5 Plus in July 2026. Also, a baby version of Theorem~\ref{thm:contact-limit-set}, with the arithmeticity assumption, arose from that conversation. These contents were soon verified, corrected, revised, and improved by the authors. ChatGPT 6 Astra is used to fulfill the proof of Theorem~\ref{thm:universal-hermitian-avoidance} and to improve the language of this manuscript. 

The authors note that Theorem~\ref{thm:main-A} and Theorem~\ref{thm:strict-zero-divisor} were almostly proved in \cite{DanielDeroin19} and the proof in this manuscript follows their idea. After a retrospective literature review aimed to identify mathematical precedents for the observations suggested by the large language model, the authors found that the comparison between the average Cartan spectrum and Lyapunov spectrum is already well-studied in a more classical set-up, see e.g. Ruelle \cite[(1.7)]{ruelle1982large}. 

All mathematical statements and cited references in the manuscript have been independently checked by the authors, who take full responsibility for its content.

\paragraph{Acknowledgement}

F. Yu would like to thank Dawei Chen, Anton Zorich and Kang Zuo for helpful discussions. J. Zhang is grateful to Peter Smillie for introducing \cite{DanielDeroin19} to him in December 2025 and many enlightening discussions, which inspire him to give a proof of Proposition~\ref{prop:lyap-coincidence}\romitem{2}. J. Zhang would also like to thank Minhao Fan for organizing a reading seminar in March 2026 on \cite{le2016brownian} that helps him to learn the Brownian motion on the Euclidean space and Feiyang Yu for helpful discussions about the Brownian motion on Riemannian manifolds.

F. Yu is supported by the Fundamental Research Funds for the Central Universities 226-2024-00136. J. Zhang is supported by NSF of China grant No. 125B2007 and partially supported by the National Key R\&D Program of China No. 2022YFA1006600, the Fundamental Research Funds for the Central Universities, and Nankai Zhide foundation. 

\section{Finite-area local systems, tame harmonic metrics, and Lyapunov exponents}
\label{sec:setup}

Let $\overline{X}$ be a closed Riemann surface with a finite subset $D\subset\overline{X}$. Set $X=\overline{X}\setminus D$. We assume that $2g(\overline{X})-2+|D|>0$, where $g(\overline{X})$ denotes the genus of $\overline{X}$. Let $\KK_X$ be the canonical line bundle over $X$, and then $\KK_{\overline{X}}(D)$ is its logarithmic extension to $\overline{X}$.
When the latter is regarded as a filtered line bundle with trivial filtration, we write $\KK_{\overline{X}}(D)_\bullet$. The Riemann surface $X$ then admits a unique complete conformal hyperbolic metric with finite area. Throughout the paper we make the metric normalized with curvature $-4$, as in \cite{EKMZ18}. If $z=x+\iu y$ is a coordinate on the upper half-plane, then $g_{\mathrm{hyp}}=(\mathrm{d}x^2+\mathrm{d}y^2)/(4y^2)$ and $\operatorname{Area}(X)=\pi\deg\KK_{\overline{X}}(D)/2$. We write $\omega_{\mathrm{hyp}}=\mathrm{d}A_{\mathrm{hyp}}$ for its K\"ahler, equivalently area, form. 
We denote the normalized hyperbolic area measure by $\mu_X$. Write $\mu_L$ for the normalized Liouville probability measure on its unit tangent bundle $T^1X$. More precisely, if $\sigma_x$ is normalized angular measure on $T_x^1X$, then $\mu_L$ is characterized by $\mathrm{d}\mu_L(v)=\mathrm{d}\mu_X(x)\,\mathrm{d}\sigma_x(v)$ for $v\in T_x^1X$. It is invariant under the geodesic flow.

\subsection{Local system and variation of Hodge structure}

Let $\mathbb{V}$ be a complex local system of rank $r$ over $X$ with associated flat bundle $(V,\nabla)$. The $(0,1)$-part of its flat connection defines the de Rham holomorphic bundle $\mathcal{V}=(V,\nabla^{0,1})$. We assume that $\mathbb{V}$ is \emph{semisimple}, i.e. its monodromy representation $\pi_1(X)\to\mathrm{GL}(V)$ is completely reducible, and that its cusp monodromy is \emph{non-expanding}: every eigenvalue of each local monodromy (the image of a simple closed curve around a cusp) $T_p$ has absolute value one. These hypotheses are used in the harmonic-metric and Lyapunov theory below, see \cite{Simpson90,Mochizuki07,EKMZ18}.

For the non-expanding local system fixed above, the associated flat bundle admits a Deligne extension across $D$ on which the connection has logarithmic poles \cite[\upRoman{2}, Proposition~5.4]{Deligne70} as follows. Fix a cusp and let $T$ denote analytic continuation of flat sections around a positively oriented loop. Write $T=T_s\mathrm e^N$, with $N$ nilpotent and $[T_s,N]=0$. We use the Deligne extension whose residue eigenvalues lie in $[0,1)$. Accordingly, write $V=\bigoplus_{\alpha\in[0,1)}V_\alpha$, where $T_s|_{V_\alpha}=\mathrm e^{-2\pi i\alpha}\operatorname{id}$. If $q=\mathrm e^{2\pi iz}$ and $v(z)$ is a multivalued flat section in $V_\alpha$, then $\mathrm e^{2\pi i\alpha z}\mathrm e^{-zN}v(z)$ is single-valued. These sections generate the Deligne extension with corresponding parabolic weight $\alpha$. We denote the resulting filtered Deligne extension by $\mathcal{V}_\bullet$.

For a finite-dimensional vector space $E$, a \emph{finite decreasing filtration indexed by $\mathbb Z$} is a family of subspaces $F^pE$ such that
\[F^{p+1}E\subset F^pE,\qquad F^pE=E\ \text{for }p\ll0,\qquad F^pE=0\ \text{for }p\gg0.\]
Thus increasing the index makes the subspace smaller, and repetitions are allowed. Its \emph{$p$-th graded piece} is the quotient
\[\operatorname{Gr}_F^pE:=F^pE/F^{p+1}E.\]
For a filtered space $E$ and a subspace $S$, we use
\[F^pS=F^pE\cap S,\qquad F^p(E/S)=(F^pE+S)/S,\qquad\operatorname{Gr}_F^pE=F^pE/F^{p+1}E.\]

We use the following convention for parabolic degree \cite{Simpson90,EKMZ18}. Let $\mathcal E_\bullet$ be a filtered holomorphic bundle on $\overline{X}$, with underlying holomorphic bundle $\mathcal E$. At each $p\in D$, let $\mathcal E_p=F_{p,1}\supsetneq\cdots\supsetneq F_{p,m_p}\supsetneq F_{p,m_p+1}=0$ be the parabolic flag, and attach weights $0\leqslant\alpha_{p,1}<\cdots<\alpha_{p,m_p}<1$ to the graded pieces $F_{p,j}/F_{p,j+1}$. Then
\[\operatorname{pardeg}(\mathcal E_\bullet)
=\deg(\mathcal E)+\sum_{p\in D}\sum_{j=1}^{m_p}\alpha_{p,j}\dim(F_{p,j}/F_{p,j+1}),
\qquad
\mu_{\mathrm{par}}(\mathcal E_\bullet)=\frac{\operatorname{pardeg}(\mathcal E_\bullet)}{\operatorname{rank}\mathcal E}.\]
For a filtered line $\LL_\bullet$ of weight $\beta_p$ at $p$, this reduces to $\operatorname{pardeg}(\LL_\bullet)=\deg\LL+\sum_{p\in D}\beta_p$. 

\paragraph{Harder--Narasimhan filtrations and steps.}
A parabolic bundle $\mathcal B_\bullet$ is \emph{semistable} if every nonzero proper holomorphic subbundle, with its induced parabolic structure, has parabolic slope at most $\mu_{\mathrm{par}}(\mathcal B_\bullet)$.
Its \emph{Harder--Narasimhan filtration} is the unique filtration
\[
0=\mathcal B^{(0)}_\bullet
\subsetneq\mathcal B^{(1)}_\bullet
\subsetneq\cdots\subsetneq
\mathcal B^{(\ell)}_\bullet=\mathcal B_\bullet
\]
by saturated holomorphic subbundles with the induced parabolic structures, whose successive quotients $\mathcal Q^{(j)}_\bullet:=\mathcal B^{(j)}_\bullet/\mathcal B^{(j-1)}_\bullet$ are semistable and have strictly decreasing parabolic slopes. One can see \cite[Proposition~5.4]{Yokogawa93} for references.
We call the subbundles $\mathcal B^{(j)}_\bullet$ the \emph{HN steps}.
Writing $r=\operatorname{rank}\mathcal B$, the vector of \emph{HN slopes}
\[
\boldsymbol{\mu}_{\mathrm{HN}}(\mathcal B_\bullet):=\left(\mu_1(\mathcal B_\bullet),\ldots, \mu_r(\mathcal B_\bullet)\right)
\]
is obtained by repeating each $\mu_{\mathrm{par}}(\mathcal Q^{(j)}_\bullet)$ according to $\operatorname{rank}\mathcal Q^{(j)}_\bullet$, in decreasing order. Thus, for $r_j=\operatorname{rank}\mathcal B^{(j)}$,
\[
\operatorname{pardeg}\mathcal B^{(j)}_\bullet
=\sum_{i=1}^{r_j}\mu_i(\mathcal B_\bullet).
\]
When $D=\varnothing$, these are the ordinary HN filtration, steps, and slopes.

\paragraph{Polarized Hodge filtrations and period maps.}
We fix the following conventions. Let $U_{\mathbb R}$ be a finite-dimensional real vector space, let $w\in\mathbb Z$, let $\mathbf h=(h^{p,q})_{p+q=w}$ with $h^{p,q}=h^{q,p}$, and let $Q$ be a nondegenerate real bilinear form satisfying that $Q(v,u)=(-1)^wQ(u,v)$. A \emph{polarized real Hodge structure of weight $w$ and Hodge numbers $\mathbf h$} is a decreasing filtration $F^\bullet$ of $U_{\mathbb C}$, which is the complexification of $U_{\mathbb{R}}$, such that
\begin{enumerate}[label=\textup{(\roman*)}]
	\item $\dim \operatorname{Gr}_F^p U_{\mathbb C}=h^{p,w-p}$;

	\item $U_{\mathbb C}=F^p\oplus\overline{F^{w-p+1}}$ for every $p$;

	\item $Q(F^p,F^{w-p+1})=0$ for every $p$;

	\item if
		\[U_F^{p,q}:=F^p\cap\overline{F^q},\qquad
		C_F|_{U_F^{p,q}}=\sqrt{-1}^{\,p-q}\operatorname{id},\]
		 then the Hermitian form
		\[h_F(u,v):=(-1)^{w(w-1)/2}Q(C_Fu,\overline v)\]
is positive definite.
\end{enumerate}
Condition \romitem{2} gives the Hodge decomposition $U_{\mathbb C}=\bigoplus_{p+q=w}U_F^{p,q}$, and conditions \romitem{3}--\romitem{4} are the Hodge--Riemann bilinear relations.
We denote the corresponding classifying space by $\mathcal D_{\mathbb R}(U_{\mathbb R},Q,\mathbf h)$.

For a complex local system we use the parallel complex convention. If $U$ is complex and $H$ is a nondegenerate Hermitian form, a \emph{polarized complex Hodge structure} of weight $w$ and type $\mathbf h$ is a decomposition
\[U=\bigoplus_{p+q=w}U^{p,q},\qquad \dim U^{p,q}=h^{p,q},\]
which is $H$-orthogonal and for which
$(-1)^qH$ is positive definite on $U^{p,q}$. Its Hodge filtration is
$F^p=\bigoplus_{r\geqslant p}U^{r,w-r}$. We write
$\mathcal D_{\mathbb C}(U,H,\mathbf h)$ for the corresponding
classifying space.

Let $Y$ be a complex manifold. A polarized real variation of Hodge structure consists of a flat real local system $\mathbb U_{\mathbb R}$, a flat polarization $Q$, and holomorphic subbundles $F^p\subset\mathcal U:=\mathbb U_{\mathbb R}\otimes_{\mathbb R}\mathcal O_Y$ such that every fiber $F_y^\bullet$ is a polarized real Hodge filtration of fixed weight and Hodge numbers and $\nabla F^p\subset F^{p-1}\otimes\Omega_Y^1$ for every $p$. The last condition is \emph{Griffiths transversality}. A polarized complex variation is defined in the same way from a flat complex local system with flat Hermitian polarization and a fiberwise polarized complex Hodge decomposition whose associated filtration is holomorphic and satisfies the Griffiths transversality.

Equivalently, after flat trivialization over the universal cover $\widetilde Y$, the Hodge filtration is encoded by a period map $\mathcal{P}\colon\widetilde Y\to\mathcal D$, where $\mathcal D$ is the appropriate real or complex classifying space. The variation axioms are precisely that $\mathcal{P}$ is holomorphic, is equivariant for the monodromy action, and is \emph{horizontal}, meaning that its differential satisfies the Griffiths transversality. Conversely, every monodromy-equivariant holomorphic horizontal map to $\mathcal D$ pulls back the tautological Hodge filtration and descends to a polarized variation on $Y$. See \cite{GriffithsSchmid69,Schmid73}.

The monodromy theorem \cite[Lemma 4.5]{Schmid73}, attributed to Borel, yields that polarized complex variation of Hodge structure over the punctured Riemann surface $X$ automatically has non-expanding cusp monodromy. 

\subsection{The tame harmonic metric and cusp estimates}

For an Hermitian metric $h$ on the flat bundle $(V,\nabla)$, write
\[
\nabla=D_h+\Phi+\Phi^{*_h},
\]
where $D_h$ is $h$-unitary, $\Phi$ is of type $(1,0)$, and $\Phi^{*_h}$ is its $h$-adjoint. 
The metric is \emph{harmonic} if $D_h^{0,1}\Phi=0$. Then $F_{D_h}+[\Phi,\Phi^{*_h}]=0$, where $F_{D_h}$ denotes the curvature form of $D_h$. Thus $\mathcal E=(V,D_h^{0,1})$ is a holomorphic bundle and $\Phi$ is holomorphic, which is called \emph{Higgs field}. The metric is called \emph{tame} if, near each puncture with coordinate $w$ centered at the puncture, the eigenvalues of $\Phi(w\cdot\partial/\partial w)$ remain bounded as $w\to0$. The de Rham bundle $\mathcal{V}$ has holomorphic structure $\nabla^{0,1}=D_h^{0,1}+\Phi^{*_h}$.

For a variation of Hodge structure, write
\[\Phi_{p,q}:=\Phi|_{V^{p,q}}:
V^{p,q}\longrightarrow V^{p-1,q+1}\otimes\KK_X.\]
An Hermitian metric $h$ on a holomorphic bundle is \emph{acceptable} \cite[\S3]{Simpson90} if, near each puncture with coordinate $w$ centered at the puncture, its Chern curvature $F_h$ satisfies
\[
|F_h|\leqslant f+\frac{C}{|w|^2(\log|w|)^2},
\qquad f\in L^p,\quad p>1,
\]
where the norms and integrability are taken with respect to the
Euclidean metric in $w$.

Following \cite[\S2.3]{EKMZ18}, we denote the metric extension by $\Xi_h$. Let $\mathcal{B}$ be a holomorphic bundle on a punctured disk $U^*=U\setminus\{p\}$, with coordinate $q$ centered at $p$ and Hermitian metric $h$. For $a\in\mathbb R$, set
\begin{equation}
\Xi_h(\mathcal{B})^{\geqslant a}(U):=
\left\{s\in H^0(U^*,\mathcal{B}):
\|s\|_h=O\left(|q|^{a-\varepsilon}\right)
\text{ for every }\varepsilon>0\right\}.
\label{eq:metric-extension}
\end{equation}
The resulting filtered sheaf $\Xi_h(\mathcal{B})$ is the \emph{metric extension} of $(\mathcal{B},h)$. For bundles arising from a tame harmonic bundle, these sheaves are locally free, decrease with $a$, and satisfy
\[
\Xi_h(\mathcal{B})^{\geqslant a+1}
=\Xi_h(\mathcal{B})^{\geqslant a}(-p).
\]
When the metric extension is fixed, we also write $\mathcal{B}_\bullet=\Xi_h(\mathcal{B})$ and $\mathcal{B}_a=\Xi_h(\mathcal{B})^{\geqslant a}$.
Its underlying holomorphic bundle is $\mathcal{B}_0$, and the images of $\mathcal{B}_a$ in its fiber at $p$, for $0\leqslant a<1$, give the parabolic flag and weights. In particular, a local frame whose norm is $|w|^a$ up to powers of $|\log|w||$ has weight $a$ modulo integers.
We say that $h$ is \emph{adapted} to a prescribed filtered bundle $\mathcal{B}_\bullet$ if $\Xi_h(\mathcal{B})^{\geqslant a}=\mathcal{B}_a$ for every $a$. 

Simpson's tame nonabelian Hodge correspondence \cite[Section~6]{Simpson90} then supplies a harmonic metric $h$ on the ordinary local system $\mathbb{V}$ adapted to the filtered Deligne extension $\mathcal{V}_\bullet$. We fix this metric, so $\Xi_h(\mathcal{V})=\mathcal{V}_\bullet$, and write $\mathcal{E}_\bullet:=\Xi_h(\mathcal{E})$ for its metric extension on the holomorphic structure $D_h^{0,1}$ given by the harmonic metric $h$.

\begin{proposition}
\label{prop:adapted-metric}
Let $(V,\nabla)$ be the semisimple flat bundle fixed above, with non-expanding cusp monodromy, and let $h$ be a tame harmonic metric adapted to $\mathcal{V}_\bullet$. Write $\nabla=D_h+\Phi+\Phi^{*_h}$ and set $\mathcal{E}=(V,D_h^{0,1})$ and $\mathcal{E}_\bullet=\Xi_h(\mathcal{E})$.
Then there is a constant $C_0>0$ such that:
\begin{enumerate}[label=\textup{(\roman*)}]
\item $|\Phi|_{h,g_{\mathrm{hyp}}}\leqslant C_0$ and $|F_{D_h}|_{h,g_{\mathrm{hyp}}}\leqslant C_0$ on $X$.

\item If $u$ is a nonzero $\nabla$-flat section on a simply connected open set, then $|\mathrm{d}\log\|u\|_h|_{g_{\mathrm{hyp}}}\leqslant 2C_0$.

\item The preceding assertions are preserved under duals and exterior powers. They also apply to a semisimple flat subquotient equipped with its tame harmonic metric.
\end{enumerate}
\end{proposition}

\begin{proof}
The connection $D_h$ is $h$-unitary and has $(0,1)$-part equal to the Dolbeault operator of $\mathcal{E}$, so it is the Chern connection of $(\mathcal{E},h)$. 
In a local holomorphic frame of $\mathcal{E}_0$, write $\Phi=B(q)\,\mathrm{d}q/q$ with $B(q)$ holomorphic. The logarithmic  residue of the Higgs field is $B(0)$. By Simpson's correspondence \cite{Simpson90}, $B(0)$ must be nilpotent. 
Simpson's main estimate \cite[Theorem~1]{Simpson90} yields that the norm of $\Phi$ is bounded from above and below by positive constants. Hence we obtain $|\Phi|_{h,g_{\mathrm{hyp}}}$ is bounded globally because the complement is compact. The Hitchin equation $F_{D_h}+[\Phi,\Phi^{*_h}]=0$ then implies a uniform bound for $F_{D_h}$. Thus \romitem{1} follows.

For \romitem{2}, let $u$ be a flat section and let $\xi$ be a real tangent vector. Since $D_h$ is compatible with $h$, we obtain that
\[\xi\log\|u\|_h=\frac{\operatorname{Re}\langle D_{h,\xi}u,u\rangle_h}{\|u\|_h^2}=-\frac{\operatorname{Re}\langle(\Phi+\Phi^{*_h})(\xi)u,u\rangle_h}{\|u\|_h^2}.\]
Hence $|\xi\log\|u\|_h|\leqslant2|\Phi(\xi)|_h$, and \romitem{2} follows from \romitem{1}.

Finally, duals and exterior powers of a tame harmonic bundle are again tame harmonic bundles, and the bounds above are preserved up to constants depending only on the fixed representation. A semisimple flat subquotient again has non-expanding cusp monodromy and therefore carries its own tame harmonic metric. This proves \romitem{3}.
\end{proof}

For later use, Proposition~\ref{prop:adapted-metric}\romitem{2} implies that if $z,w\in\widetilde X$ and $u$ is flat, then $|\log\|u\|_{h,z}-\log\|u\|_{h,w}|\leqslant C d_{\mathrm{hyp}}(z,w)$. Applying this to $\bigwedge^kV$ gives the same estimate for the logarithm of the product of the first $k$ singular values of flat identification.

\subsection{Lyapunov spectra for the geodesic flow}\label{subsec:Lya}

Let $g_t$ be the unit-speed geodesic flow on the unit tangent bundle $T^1X$. Parallel transport with respect to the flat connection $\nabla$ along the geodesic segment determined by $v\in T^1X$ gives a cocycle $G_t(v)$ on the pullback bundle of $V$ over $T^1X$. Let $\sigma_1(G_t(v))\geqslant\cdots\geqslant\sigma_r(G_t(v))$ be its singular values for $h$. Eskin--Kontsevich--M\"oller--Zorich \cite[Theorem~2.1]{EKMZ18} gives the required integrability for the existence of the corresponding Lyapunov exponents and they do not depend on the chosen admissible metric (see \cite[Definition 2.3]{EKMZ18}). Thus, for $\mu_L$-almost every $v$, the $i$-th Lyapunov exponent
\begin{equation}
\lambda_i^{\mathrm{geo}}(V)=\lim_{t\to\infty}\frac1t\log\sigma_i(G_t(v)),\qquad 1\leqslant i\leqslant r,
\label{eq:geo-lyap}
\end{equation}
exists by Oseledets' theorem \cite{Oseledets68}. We write
\[\boldsymbol{\lambda}^{\mathrm{geo}}(V):=\left(\lambda_1^{\mathrm{geo}}(V),\ldots,\lambda_r^{\mathrm{geo}}(V)\right)\]
for this geodesic Lyapunov spectrum, with multiplicities. In particular, the top exponent of $\bigwedge^kV$ is $\lambda_1^{\mathrm{geo}}(V)+\cdots+\lambda_k^{\mathrm{geo}}(V)$.

Duality gives $\lambda_i^{\mathrm{geo}}(V^\vee)=-\lambda_{r+1-i}^{\mathrm{geo}}(V)$. In addition, the geodesic flow is reversible and Liouville measure is invariant under $v\mapsto-v$, hence $\lambda_i^{\mathrm{geo}}(V)=-\lambda_{r+1-i}^{\mathrm{geo}}(V)$. Therefore $\lambda_1^{\mathrm{geo}}(V^\vee)=\lambda_1^{\mathrm{geo}}(V)$.

\paragraph{Lyapunov spectrum for a semisimple group.}
Suppose more generally that the monodromy factors through a connected semisimple real Lie group $G_{\mathbb R}$ with finite center and without compact factors. Fix a Cartan decomposition $\mathfrak{g}_{\mathbb R}=\mathfrak{k}_{\mathbb{R}}\oplus\mathfrak{p}_{\mathbb{R}}$. Let $K_{\mathbb R}<G_{\mathbb R}$ be the corresponding maximal compact subgroup, choose a maximal abelian subspace $\mathfrak a\subset\mathfrak{p}_{\mathbb{R}}$, and fix a closed positive restricted chamber $\overline{\mathfrak a^+}$. Let $\Delta$ be the simple restricted roots, let $\omega_\alpha$ be the fundamental restricted weights. For $Y,Z\in\mathfrak{a}$, we write $Y\preceq Z$ when $\omega_\alpha(Y)\leqslant\omega_\alpha(Z)$ for every $\alpha\in\Delta$. For $g\in G_{\mathbb R}$, its Cartan decomposition determines a unique $\kappa(g)\in\overline{\mathfrak a^+}$ by $g\in K_{\mathbb R}\exp(\kappa(g))K_{\mathbb R}$. See e.g. \cite[Ch.~\upRoman{9}, Thm.~1.1]{Helgason01}. $\kappa$ is called the \emph{Cartan projection}.

We shall use throughout the following Cartan triangle inequality
\begin{equation}
\kappa(g_2g_1)\preceq\kappa(g_1)+\kappa(g_2).
\label{eq:Cartan-triangle}
\end{equation}
For completeness, let $\sigma$ be the chamber-valued distance on $G_{\mathbb R}/K_{\mathbb R}$ and $o=K_{\mathbb R}$. Then $\sigma(o,go)=\kappa(g)$, while $\sigma(x,z)\preceq\sigma(x,y)+\sigma(y,z)$ by \cite[Remark~3.36\romitem{2}]{KLM09Convex}. Taking $x=o$, $y=g_2o$, and $z=g_2g_1o$ gives \eqref{eq:Cartan-triangle}.

Let $U_t(v)$ denote the $G_{\mathbb R}$-valued flat cocycle over the geodesic flow, written with respect to the compatible tame harmonic $K_{\mathbb R}$-reduction fixed above. The required integrability follows in our setting by applying Proposition~\ref{prop:adapted-metric}\romitem{2} to the adjoint local system: the restricted roots span $\mathfrak a^*$, so the singular-value bounds for the adjoint cocycle control its Cartan projection. The multiplicative ergodic theorem for semisimple Lie-group cocycles gives a unique $\boldsymbol{\lambda}\in\overline{\mathfrak a^+}$ such that
\begin{equation}
\frac1t\kappa(U_t(v))\longrightarrow\boldsymbol{\lambda}
\qquad\text{for }\mu_L\text{-almost every }v.
\label{eq:G-Lyap-spectrum}
\end{equation}
We call $\boldsymbol{\lambda}$ the \emph{Lyapunov spectrum} of the $G_{\mathbb R}$-local system, and write $\lambda_\alpha:=\omega_\alpha(\boldsymbol{\lambda})$. See Kaimanovich \cite{Kaimanovich89} for semisimple groups and Alves--San Martin \cite{AlvesSanMartin13} for a reductive flag-bundle formulation.

In the present setting, the existence of \eqref{eq:G-Lyap-spectrum} can also be seen directly from \eqref{eq:Cartan-triangle}. For each $\alpha\in\Delta$, the cocycle $\omega_\alpha(\kappa(U_t(v)))$ is subadditive. Apply Kingman's subadditive ergodic theorem \cite{Kingman68} to the time-one cocycle. Since the geodesic flow is mixing, its time-one map is ergodic, and the limit at integer times is almost surely constant. Proposition~\ref{prop:adapted-metric}\romitem{2} bounds the Cartan growth on intervals of length at most one, so the same limit holds as $t\to\infty$ through arbitrary real times. Thus $\lim_{t\to\infty}t^{-1}\omega_\alpha(\kappa(U_t(v)))=\lambda_\alpha$ almost everywhere. The fundamental restricted weights form a basis of $\mathfrak a^*$, so the $\lambda_\alpha$ determine a unique $\boldsymbol{\lambda}\in\overline{\mathfrak a^+}$ and give \eqref{eq:G-Lyap-spectrum}. More generally, if $\tau\colon G_{\mathbb R}\to\mathrm{GL}(W)$ has highest restricted weight $\chi_\tau$, then $\lambda_1^{\mathrm{geo}}(W)=\chi_\tau(\boldsymbol{\lambda})$, and $\boldsymbol{\lambda}^{\mathrm{geo}}(W)$ is obtained by evaluating the restricted weights of $W$ on $\boldsymbol{\lambda}$, with multiplicities. In particular, if $G_{\mathbb R}=\mathrm{SL}_r\mathbb R$ and $V=\mathbb R^r$ is the standard representation, then $\boldsymbol{\lambda}=(\lambda_1^{\mathrm{geo}}(V),\ldots,\lambda_r^{\mathrm{geo}}(V))$ and $\omega_k(\boldsymbol{\lambda})=\lambda_1^{\mathrm{geo}}(V)+\cdots+\lambda_k^{\mathrm{geo}}(V)$.

\subsection{Brownian motion and stationary distributions}

For the curvature-($-4$) hyperbolic metric we use $\Delta_{\mathrm{hyp}}=4y^2(\partial_x^2+\partial_y^2)$ and Brownian generator $L_{\mathrm{Br}}=\frac12\Delta_{\mathrm{hyp}}$. This normalization makes the Brownian and unit-speed geodesic Lyapunov exponents agree.

We first recall the probabilistic objects that enter the Brownian definition of the Lyapunov spectrum. By \emph{Brownian motion} on $(X,g_{\mathrm{hyp}})$ we mean the diffusion $Z_t$ with generator $L_{\mathrm{Br}}$. Equivalently, for every $f\in C_c^\infty(X)$,
\[
f(Z_t)-f(Z_0)-\frac12\int_0^t\Delta_{\mathrm{hyp}}f(Z_s)\,\mathrm{d}s
\]
is a local martingale. Since $X$ is complete and has Ricci curvature bounded from below, Brownian motion is non-explosive and is defined for all $t\geqslant0$. See \cite[Chapters~3--4]{Hsu02} for references.

A probability measure $\mu$ on $X$ is called a \emph{stationary distribution} for Brownian motion if, whenever $Z_0$ has law $\mu$, the random variable $Z_t$ has the same law $\mu$ for every $t\geqslant0$. Equivalently, if $H_tf(x)=\mathbf E_x[f(Z_t)]$ is the heat semigroup, then $\int_XH_tf\,\mathrm{d}\mu=\int_Xf\,\mathrm{d}\mu$ for every bounded Borel function $f$. In particular, the normalized hyperbolic area $\mu_X=\operatorname{Area}(X)^{-1}\mathrm{d}A_{\mathrm{hyp}}$ is stationary. Thus saying that Brownian motion is \emph{started in stationarity} means that $Z_0\sim\mu_X$, and then $Z_t\sim\mu_X$ for every $t$.

Parallel transport along $Z_{[0,t]}$ gives a random linear map $P_t\colon V_{Z_0}\to V_{Z_t}$. The norms at the two endpoints are taken with respect to the tame harmonic metric $h$. Similar to the geodesic flow case, Proposition~\ref{prop:adapted-metric}\romitem{2} gives the required integrability for the stationary multiplicative ergodic theorem. And it applies to $P_t$ and its exterior powers, which gives the spectrum $\lambda_1^{\mathrm{Br}}\geqslant\cdots\geqslant\lambda_r^{\mathrm{Br}}$ such that, almost surely,
\[
\lambda_1^{\mathrm{Br}}+\cdots+\lambda_k^{\mathrm{Br}}
=\lim_{t\to\infty}\frac1t\log\|\textstyle\bigwedge^kP_t\|,
\qquad 1\leqslant k\leqslant r.
\]
We write
\[
\boldsymbol{\lambda}^{\mathrm{Br}}(V)
:=
\left(\lambda_1^{\mathrm{Br}}(V),\ldots,\lambda_r^{\mathrm{Br}}(V)\right)
\]
for the Brownian Lyapunov spectrum, with multiplicities, following
\cite[Section~1]{DanielDeroin19}.

\subsection{Flat foliation, leafwise Chern forms, and harmonic measures}
\label{subsec:flat-foliation}

For a right principal $H$-bundle $\mathbb E\to X$ and a left $H$-space $S$, the \emph{associated fibre product} is
\[\mathbb E\times_H S:=(\mathbb E\times S)/{\sim},
\qquad (eh,s)\sim(e,hs).\]
For a representation $\rho\colon\pi_1(X)\to H$, we write
\[
\widetilde X\times_\rho S:=(\widetilde X\times S)/\pi_1(X),
\qquad \gamma\cdot(z,s)=(\gamma z,\rho(\gamma)s).
\]
This is a fibre bundle over $X$ with fibre $S$. 

We now introduce the foliation structure on the projective bundle. Fix a basepoint $x_0\in X$ and let $\rho\colon\pi_1(X,x_0)\to\operatorname{GL}(V_{x_0})$ be the monodromy of $\mathbb{V}$. Parallel transport identifies $\mathbb P(V)$ with the associated fibre product $\widetilde X\times_\rho\mathbb P(V_{x_0})$. Write $p\colon\mathbb P(V)\to X$ for its projection. The \emph{flat foliation} $\mathcal F$ is the foliation obtained from the horizontal slices $\widetilde X\times\{[v]\}$. Equivalently, if $U\subset X$ is a coordinate disk and $V|_U$ is trivialized by flat sections, then $p^{-1}(U)\cong U\times\mathbb P^{r-1}$ and, for a coordinate neighborhood $T\subset\mathbb P^{r-1}$, $U\times T$ is a foliated chart for $\mathcal F$, with local leaves $U\times\{[v]\}$ for $[v]\in T$. Thus $p$ restricts to a local biholomorphism from every leaf to $X$. Each leaf therefore carries the pulled-back complex structure and hyperbolic metric.

For a function or differential form that is smooth along the leaves, write $\partial_{\mathcal F}$ and $\bar\partial_{\mathcal F}$ for the Dolbeault operators in the leaf direction, and set $\mathrm{d}_{\mathcal F}=\partial_{\mathcal F}+\bar\partial_{\mathcal F}$ and $\mathrm{d}_{\mathcal F}^c=(\sqrt{-1}/4\pi)(\bar\partial_{\mathcal F}-\partial_{\mathcal F})$. Hence $\mathrm{d}_{\mathcal F}\mathrm{d}_{\mathcal F}^c=(\sqrt{-1}/2\pi)\partial_{\mathcal F}\bar\partial_{\mathcal F}$. If $(\LL,h_{\LL})$ is a line bundle that is holomorphic along the leaves and $s$ is a nowhere-vanishing leafwise holomorphic frame, its leafwise first Chern form is $c_1^{\mathcal F}(\LL,h_{\LL})=-\mathrm{d}_{\mathcal F}\mathrm{d}_{\mathcal F}^c\log\|s\|_{h_{\LL}}^2$.

\begin{proposition}
\label{prop:leafwise-projective-curvature}
Let $h_{\mathrm{FS}}$ be the Fubini--Study metric induced by $h$ on $\mathcal O_{\mathbb P(V)}(1)$. There is a constant $C_1>0$ such that $|c_1^{\mathcal F}(\mathcal O_{\mathbb P(V)}(1),h_{\mathrm{FS}})|\leqslant C_1p^*\omega_{\mathrm{hyp}}$.
\end{proposition}

\begin{proof}
We fix a foliated chart $U\times T$ and choose a holomorphic coordinate $z=x+\iu y$ on $U$, and trivialize $V|_U$ by a flat frame. Thus $\nabla=\mathrm{d}$, and the metric is represented by a positive Hermitian matrix $H$. Fix the transverse coordinate $[u]\in T$, represented by a constant vector $0\neq u\in\mathbb C^r$, and set
\[
q=u^*Hu=\|u\|_H^2,\qquad f=\log q,
\qquad A=H^{-1}H_z,\qquad B=H^{-1}H_{\bar z},
\]
where $*$ denotes conjugate transpose. The identities $H^*=H$ and $H_z^*=H_{\bar z}$ give $B=H^{-1}A^*H=A^{*_H}$.

Write $D_h=\mathrm{d}+\Gamma$ in this frame. The difference $\nabla-D_h=-\Gamma$ is $h$-self-adjoint on real tangent vectors, while metric compatibility gives
$\mathrm{d}H=\Gamma^*H+H\Gamma=2H\Gamma$. Therefore
\[
D_h=\mathrm{d}+\frac12H^{-1}\mathrm{d}H,
\qquad
\Phi=-\frac12A\,\mathrm{d}z,
\qquad
\Phi^{*_h}=-\frac12B\,\mathrm{d}\bar z.
\]
On endomorphism-valued forms, $D_h^{0,1}$ acts by $\bar\partial+[\frac12B\,\mathrm{d}\bar z,\,\cdot\,]$. Thus harmonicity gives
\[
0=D_h^{0,1}\Phi
=-\frac12\left(\partial_{\bar z}A+\frac12[B,A]\right)
\mathrm{d}\bar z\wedge\mathrm{d}z,
\qquad
\partial_{\bar z}A=\frac12[A,B].
\]
As $u$ is constant, differentiation now yields
\[
q_z=u^*HAu,\qquad
q_{\bar z}=u^*HBu=\overline{q_z},
\]
\[q_{z\bar z}=u^*\left(H_{\bar z}A+H\partial_{\bar z}A\right)u=u^*H\left(BA+\frac12[A,B]\right)u
=\frac12u^*H(AB+BA)u=\frac12\left(\|Au\|_H^2+\|Bu\|_H^2\right).
\]
The last equality uses $B=A^{*_H}$.

The $H$-orthogonal projection onto $\mathbb Cu$ is
$P_uv=u(u^*Hv)/q$. This yields that
\[
\|P_uAu\|_H^2
=\frac{|u^*HAu|^2}{q}
=\frac{|q_z|^2}{q}
=\|P_uBu\|_H^2.
\]
Using the orthogonal decompositions of $Au$ and $Bu$, we obtain that
\[
f_{z\bar z}=\frac{q_{z\bar z}}q-\frac{q_zq_{\bar z}}{q^2}=\frac{\|Au\|_H^2+\|Bu\|_H^2}{2q}-\frac{|q_z|^2}{q^2}=\frac{\|(I-P_u)Au\|_H^2+\|(I-P_u)Bu\|_H^2}{2q}.
\]
Consequently,
\[
0\leqslant f_{z\bar z}
\leqslant\frac12\left(\|A\|_{\mathrm{op},H}^2+
\|B\|_{\mathrm{op},H}^2\right)
=\|A\|_{\mathrm{op},H}^2,
\]
where $\|\cdot\|_{\mathrm{op},H}$ denotes the operator norm for $H$.

On a local leaf transverse to $[u]$, the vector $u$ is a nonvanishing holomorphic frame of $\mathcal O_{\mathbb P(V)}(-1)$ with squared norm $q$. Dualizing and using our normalization of $\mathrm{d}_{\mathcal F}^c$ therefore gives
\[
c_1^{\mathcal F}(\mathcal O_{\mathbb P(V)}(1),h_{\mathrm{FS}})
=\mathrm{d}_{\mathcal F}\mathrm{d}_{\mathcal F}^cf
=\frac{\sqrt{-1}}{2\pi}f_{z\bar z}\,
\mathrm{d}z\wedge\mathrm{d}\bar z.
\]
Finally, write $g_{\mathrm{hyp}}=\sigma^2|\mathrm{d}z|^2$, so that $\omega_{\mathrm{hyp}}=(\sqrt{-1}/2)\sigma^2\mathrm{d}z\wedge\mathrm{d}\bar z$ and $\sigma^{-1}\partial_x$ is a unit real tangent vector. Since $A=-2\Phi(\partial/\partial x)$, Proposition~\ref{prop:adapted-metric}\romitem{1} gives
\[\|A\|_{\mathrm{op},H}
=2\sigma\|\Phi(\sigma^{-1}\cdot\partial/\partial x)\|_{\mathrm{op},h}
\leqslant2C_0\sigma.\]
It follows that
\[
0\leqslant
c_1^{\mathcal F}(\mathcal O_{\mathbb P(V)}(1),h_{\mathrm{FS}})
=\frac{f_{z\bar z}}{\pi\sigma^2}\,p^*\omega_{\mathrm{hyp}}
\leqslant\frac{4C_0^2}{\pi}\,p^*\omega_{\mathrm{hyp}}.
\]
This bound is independent of the chart and of $[u]$, so the assertion holds with $C_1=4C_0^2/\pi$.
\end{proof}

If $y=(x,[v])\in\mathbb P(V)$ and $Z_t$ is Brownian motion started at $x$, parallel transport of $[v]$ along $Z_t$ gives a process $Y_t$ on the leaf through $y$. This is the \emph{leafwise Brownian motion}. Its generator is $L_{\mathcal F}:=\frac12\Delta_{\mathcal F}$, where $\Delta_{\mathcal F}$ is the hyperbolic Laplacian along the leaves.

A nonzero finite positive Borel measure $\nu$ on $\mathbb P(V)$ is a \emph{harmonic measure} in the sense of Garnett \cite{Garnett83} if it is stationary for leafwise Brownian motion. Equivalently, for every compactly supported continuous function $\varphi$ that is $C^2$ along the leaves, one has $\int_{\mathbb P(V)}\Delta_{\mathcal F}\varphi\,\mathrm{d}\nu=0$.

Since $p$ maps leafwise Brownian motion to Brownian motion on $X$, the pushforward $p_*\nu$ is stationary on $X$. By \cite[Proposition~2.6]{DanielDeroin19}, applied to $\nu/\nu(\mathbb P(V))$, this pushforward is a positive multiple of $\mu_X$:
\[
p_*\nu=\nu(\mathbb P(V))\,\mu_X.
\]
Throughout this paper, harmonic measures are normalized by $p_*\nu=\mu_X$. Daniel--Deroin \cite[Theorem~2]{DanielDeroin19} prove the existence of harmonic measure under the bounded-geometry hypotheses. Their construction starts from a probability measure projecting to $\mu_X$ and preserves this projection under the leafwise heat flow. We write $\nu$ for any such harmonic measure, with an invariant support specified when needed.

Following \cite[Section~2.3]{DanielDeroin19}, let $\mathcal W_{\mathcal F}(\mathbb P(V))$ be the space of continuous paths $\gamma\colon[0,+\infty)\to\mathbb P(V)$ contained in a single leaf of $\mathcal F$. For $y\in\mathbb P(V)$, let $\mathbb P_y^{\mathcal F}$ be the law of leafwise Brownian motion started at $y$. The measure $\nu$ determines a probability measure $\bar\nu$ on $\mathcal W_{\mathcal F}(\mathbb P(V))$ by
\[
\int F(\gamma)\,\mathrm{d}\bar\nu(\gamma)
:=\int_{\mathbb P(V)}\left(\int F(\gamma)\,\mathrm{d}\mathbb P_y^{\mathcal F}(\gamma)\right)\mathrm{d}\nu(y)
\]
for every bounded measurable $F$. For $t\geqslant0$, let $\vartheta_t$ denote the time shift on path space, defined by
\[
(\vartheta_t\gamma)(s):=\gamma(t+s),\qquad s\geqslant0.
\]
Since $\nu$ is harmonic, $\bar\nu$ is invariant under each $\vartheta_t$. See \cite[Proposition~2.9]{DanielDeroin19}.

For $\gamma\in\mathcal W_{\mathcal F}(\mathbb P(V))$, choose a nonzero flat lift $s(t)$ of $\gamma(t)$ to the tautological line $\mathcal O_{\mathbb P(V)}(-1)$ along the leaf and set
\[
\overline H_t(\gamma):=\log\frac{\|s(t)\|_h}{\|s(0)\|_h}.
\]
This is independent of the nonzero initial vector and satisfies the additive cocycle identity
\[
\overline H_{t+t'}=\overline H_t+\overline H_{t'}\circ\vartheta_t,
\qquad t,t'\geqslant0.
\]
The integrability needed below is
\[
\sup_{0\leqslant t\leqslant1}|\overline H_t|
\in L^1\left(\mathcal W_{\mathcal F}(\mathbb P(V)),\bar\nu\right).
\]
Indeed, for a lift $\widetilde\gamma$ of $p\circ\gamma$ to $\widetilde X$, Proposition~\ref{prop:adapted-metric}\romitem{2} gives
\[
|\overline H_t(\gamma)|\leqslant
2C_0d_{\mathrm{hyp}}\left(\widetilde\gamma(0),\widetilde\gamma(t)\right).
\]
The proof of \cite[Theorem~1]{DanielDeroin19} bounds the expectation of the maximum of this distance on $[0,1]$, uniformly in the starting point.
This proves the stated integrability and justifies the continuous-time ergodic limit below. In particular, $\overline H_1\in L^1(\bar\nu)$, as also follows from \cite[Lemma~2.10]{DanielDeroin19}, so we may define
\begin{equation}
\lambda(\nu):=\int_{\mathcal W_{\mathcal F}(\mathbb P(V))}\overline H_1\,\mathrm{d}\bar\nu,
\label{eq:lambda-nu-def}
\end{equation}
as in \cite[Definition~2.11]{DanielDeroin19}. By stationarity and the cocycle identity, for every $t>0$, $\lambda(\nu)=t^{-1}\int\overline H_t\,\mathrm{d}\bar\nu$. If $\nu$ is ergodic, Birkhoff's theorem gives $t^{-1}\overline H_t(\gamma)\to\lambda(\nu)$ for $\bar\nu$-almost every leafwise Brownian path \cite[Proposition~2.12]{DanielDeroin19}. 

We next associate to $\nu$ its harmonic current. Let $\Omega^{1,1}_{\mathcal F,c}(\mathbb P(V))$ denote the space of compactly supported continuous leafwise $(1,1)$-forms whose coefficients are smooth along the leaves. Since every leaf has complex dimension one and $p$ is a local isometry on each leaf, every $\alpha\in\Omega^{1,1}_{\mathcal F,c}(\mathbb P(V))$ can be written uniquely as $\alpha=\kappa_\alpha\,p^*\omega_{\mathrm{hyp}}$ for a compactly supported continuous function $\kappa_\alpha$ that is smooth along the leaves. We define
\[\begin{aligned}
 T_\nu\colon\Omega^{1,1}_{\mathcal F,c}(\mathbb P(V))&\longrightarrow\mathbb C\\
 \alpha&\longmapsto\pi\int_{\mathbb P(V)}\kappa_\alpha\,\mathrm{d}\nu.
\end{aligned}\]
For real leafwise $(1,1)$-forms the value is real. The definition extends to any measurable leafwise $(1,1)$-form $\alpha=\kappa_\alpha p^*\omega_{\mathrm{hyp}}$ with $\kappa_\alpha\in L^1(\nu)$. For a compactly supported function $\varphi$ that is $C^2$ along the leaves, our normalization gives $\mathrm{d}_{\mathcal F}\mathrm{d}_{\mathcal F}^c\varphi=(4\pi)^{-1}\Delta_{\mathcal F}\varphi\,p^*\omega_{\mathrm{hyp}}$. Hence
\[
T_\nu(\mathrm{d}_{\mathcal F}\mathrm{d}_{\mathcal F}^c\varphi)
=\frac14\int_{\mathbb P(V)}\Delta_{\mathcal F}\varphi\,\mathrm{d}\nu=0,
\]
where the last equality is precisely the stationarity of $\nu$ for leafwise Brownian motion. Thus $T_\nu$ is $(\mathrm{d}_{\mathcal F}\mathrm{d}_{\mathcal F}^c)$-closed in the leafwise sense. This is the harmonic current associated with $\nu$.

On the base we use the analogous convention $\mathrm{d}^c=(\sqrt{-1}/4\pi)(\bar\partial-\partial)$, so $\mathrm{d}\mathrm{d}^c=(\sqrt{-1}/2\pi)\partial\bar\partial$ and $\mathrm{d}\mathrm{d}^cu=(4\pi)^{-1}\Delta_{\mathrm{hyp}}u\,\mathrm{d}A_{\mathrm{hyp}}$. If $\eta=\kappa_\eta\,\mathrm{d}A_{\mathrm{hyp}}$ is a $(1,1)$-form on $X$, then $p^*\eta=(\kappa_\eta\circ p)p^*\omega_{\mathrm{hyp}}$, and therefore
\[T_\nu(p^*\eta)=\pi\int_{\mathbb P(V)}\kappa_\eta\circ p\,\mathrm{d}\nu=\frac{\pi}{\operatorname{Area}(X)}\int_X\eta,\]
because $p_*\nu=\mu_X$.

With our normalization, Daniel--Deroin \cite[Proposition~3.7]{DanielDeroin19} implies that
\begin{equation}\label{eq:harmonic-current}
\lambda(\nu)
=
T_\nu\!\left(
c_1^{\mathcal F}(\mathcal O_{\mathbb P(V)}(1),h_{\mathrm{FS}})
\right).
\end{equation}
Proposition~\ref{prop:leafwise-projective-curvature} shows that this Chern form has bounded coefficient with respect to $p^*\omega_{\mathrm{hyp}}$, so it lies in the $L^1(\nu)$ domain of the extended current. 

\subsection{Coincidence of spectra}

A representation is \emph{strongly irreducible} if its restriction to every finite-index subgroup is irreducible.
For a flat bundle $W$, a harmonic measure $\nu^{\mathrm{top}}$ on $\mathbb P(W)$ is \emph{Lyapunov-realizing} if $\lambda(\nu^{\mathrm{top}})=\lambda_1^{\mathrm{Br}}(W)$. The following proposition identifies the Brownian and geodesic spectra and realizes their top exponent by such a measure.

\begin{proposition}
\label{prop:lyap-coincidence}
Let $W$ be any flat subquotient or tensor construction of $\mathbb{V}$, equipped with a compatible tame harmonic metric. Then:

\begin{enumerate}[label=\textup{(\roman*)}]
\item the Brownian and geodesic Lyapunov spectra coincide,
\begin{equation}
\lambda_i^{\mathrm{Br}}(W)=\lambda_i^{\mathrm{geo}}(W),\qquad 1\leqslant i\leqslant \operatorname{rank}W.
\label{eq:lyap-coincidence}
\end{equation}
We write $\lambda_i(W)$ for this common value and
\[
\boldsymbol{\lambda}(W)
:=
\left(\lambda_1(W),\ldots,\lambda_{\operatorname{rank}W}(W)\right).
\]
\item for every harmonic probability measure $\nu$ on $\mathbb P(W)$ one has $\lambda(\nu)\leqslant\lambda_1(W)$. Moreover, there exists a harmonic probability measure $\nu^{\mathrm{top}}$ with $p_*\nu^{\mathrm{top}}=\mu_X$ such that $\lambda(\nu^{\mathrm{top}})=\lambda_1(W)$. Hence
\[
\lambda_1(W)
=\lambda(\nu^{\mathrm{top}})
=T_{\nu^{\mathrm{top}}}\!\left(c_1^{\mathcal F}
\left(\mathcal O_{\mathbb P(W)}(1),h_{\mathrm{FS}}\right)\right).
\]
If $W$ is strongly irreducible, every harmonic probability measure on $\mathbb P(W)$ is Lyapunov-realizing. Under the same assumption, every nonempty closed invariant subset $M\subset\mathbb P(W)$ supports a Lyapunov-realizing harmonic probability measure, where invariant means that $M$ is a union of leaves of $\mathcal F$.
\end{enumerate}
\end{proposition}

\begin{proof}
For \romitem{1}, start Brownian motion with initial distribution $\mu_X$ and lift it to $\mathbb H$. Almost surely there is $\xi\in\partial\mathbb H$ such that, if $\gamma_\xi$ is the unit-speed geodesic ray from $Z_0$ to $\xi$, then
\[
d_{\mathrm{hyp}}(Z_t,\gamma_\xi(t))=o(t).
\]
Here the parameter $t$ is the same for Brownian motion and for the geodesic because, for curvature $-4$ and generator $\frac12\Delta_{\mathrm{hyp}}$, the Brownian escape speed is $1$ (see \cite[Appendix~B.3]{FranchiLeJan12} and \cite{Hsu02}).

Let $P_t^{\mathrm{Br}}$ and $P_t^{\mathrm{geo}}$ be parallel transport along the Brownian path and along $\gamma_\xi$, and let $R_t\colon W_{\gamma_\xi(t)}\to W_{Z_t}$ be flat transport between their endpoints. Since the universal cover is simply connected,
\[
P_t^{\mathrm{Br}}=R_tP_t^{\mathrm{geo}}.
\]
By Proposition~\ref{prop:adapted-metric}\romitem{2}, applied to $\bigwedge^kW$, $\|\textstyle\bigwedge^kR_t^{\pm1}\|
\leqslant \mathrm{e}^{C_k d_{\mathrm{hyp}}(Z_t,\gamma_\xi(t))}$.
Therefore
\begin{equation}
\left|
\log\|\textstyle\bigwedge^kP_t^{\mathrm{Br}}\|
-
\log\|\textstyle\bigwedge^kP_t^{\mathrm{geo}}\|
\right|=o(t).
\label{eq:brownian-geodesic-o-t}
\end{equation}
Conditional on $Z_0=x$, the endpoint $\xi$ has visual measure, hence the initial direction $\dot\gamma_\xi(0)$ is uniformly distributed on $T_x^1X$. Since $x$ has law $\mu_X$, this initial direction has normalized Liouville law on $T^1X$. Thus Oseledets' theorem and \eqref{eq:brownian-geodesic-o-t} give, for every $k$,
\[
\sum_{i=1}^k\lambda_i^{\mathrm{Br}}(W)
=
\sum_{i=1}^k\lambda_i^{\mathrm{geo}}(W).
\]
Subtracting consecutive identities proves \romitem{1}.

For \romitem{2}, Daniel--Deroin \cite[Remark~2.17]{DanielDeroin19} associate to every harmonic measure $\nu$ a projective exponent and prove $\lambda(\nu)\leqslant\lambda_1^{\mathrm{Br}}(W)$. We construct a harmonic measure attaining equality.

Let $H$ be the Zariski closure of the monodromy of $W$, and let $H^\circ$ be its identity component. Since $W$ is semisimple, $H$ is reductive. Let $\Gamma_0:=\rho_W^{-1}(H^\circ)$ and let $q\colon X_0\to X$ be the corresponding finite cover. Then $\rho_W(\Gamma_0)$ is Zariski dense in the connected reductive group $H^\circ$. Hence $W_x=E_1\oplus\cdots\oplus E_m$ with each $E_j$ irreducible under $H^\circ$, and therefore
\begin{equation}\label{eq:flat-direct-sum}
q^*W=W_1\oplus\cdots\oplus W_m
\end{equation}
as flat bundles. Now if a finite-index subgroup of $\Gamma_0$ preserved a proper subspace of $E_j$, then its Zariski closure would also preserve it. Then the Zariski closure must be $H^\circ$ since $H^\circ$ is connected. Hence every $W_j$ must be strongly irreducible.

On $W_j$, we have the restricted metric $h_j=(q^*h)|_{W_j}$. Then the top exponent of the direct sum is
\begin{equation}\label{eq:top-direct-sum}
\lambda_1^{\mathrm{Br}}(q^*W)=\max_j\lambda_1^{\mathrm{Br}}(W_j).
\end{equation}
Indeed, if $P_t$ is Brownian parallel transport and $P_{j,t}=P_t|_{W_j}$, then $\|P_{j,t}\|\leqslant\|P_t\|$. Conversely, let $\pi_j$ be the projection onto $W_j$ with respect to \eqref{eq:flat-direct-sum}, $\|P_t\|\leqslant\sum_j\|P_{j,t}\|\,\|\pi_j\|$. Note that $\|\pi_j\|$'s are independent of $t$, so dividing the logarithm by $t$ gives \eqref{eq:top-direct-sum}. Since Brownian motion on $X_0$ projects to Brownian motion on $X$, also
\begin{equation}
\lambda_1^{\mathrm{Br}}(q^*W)=\lambda_1^{\mathrm{Br}}(W).
\label{eq:cover-top-exponent}
\end{equation}

Choose $j$ realizing the maximum in \eqref{eq:top-direct-sum}. By \cite[Theorems~2--3]{DanielDeroin19}, there is a harmonic probability measure $\widetilde\nu_j$ on $\mathbb P(W_j)$, projecting to normalized area on $X_0$, such that
\[
\lambda(\widetilde\nu_j)=\lambda_1^{\mathrm{Br}}(W_j).
\]
Push it forward through $\mathbb P(W_j)\subset\mathbb P(q^*W)\to\mathbb P(W)$ and call the result $\nu^{\mathrm{top}}$. The map intertwines leafwise Brownian motion, so $\nu^{\mathrm{top}}$ is harmonic and projects to $\mu_X$. Because $h_j$ is the restriction of $q^*h$, the projective norm cocycle is unchanged under this map. Hence, using \eqref{eq:top-direct-sum} and \eqref{eq:cover-top-exponent},
\[
\lambda(\nu^{\mathrm{top}})
=
\lambda_1^{\mathrm{Br}}(W_j)
=
\lambda_1^{\mathrm{Br}}(W)
=
\lambda_1(W),
\]
where the last equality is \romitem{1}. The last equality in the statement follows from \eqref{eq:harmonic-current}.

Finally, assume that $W$ is strongly irreducible. The base $X$ is complete, has finite area, and has Ricci curvature bounded from below. Proposition~\ref{prop:adapted-metric}\romitem{1}, applied to $W$, bounds the self-adjoint part of its flat connection. Thus
\cite[Theorem~3]{DanielDeroin19} applies and gives
$\lambda(\nu)=\lambda_1^{\mathrm{Br}}(W)=\lambda_1(W)$ for every harmonic probability measure $\nu$ on $\mathbb P(W)$.
Moreover, Proposition~\ref{prop:leafwise-projective-curvature}, applied to $W$, bounds $\Delta_{\mathcal F}\log\|s\|_h$ uniformly for local nonzero flat sections $s$ of $\mathcal O_{\mathbb P(W)}(-1)$.
This verifies the additional boundedness hypothesis of \cite[Theorem~4]{DanielDeroin19}, which gives a Lyapunov-realizing harmonic probability measure supported in any prescribed nonempty closed invariant subset $M\subset\mathbb P(W)$.
\end{proof}

\begin{convention}
From now on we use $\boldsymbol{\lambda}(W)$ to denote the common Brownian/geodesic Lyapunov spectrum of a linear flat bundle $W$, with multiplicities, and $\lambda_1(W)$ for its top exponent. For a $G_{\mathbb R}$-cocycle, we use $\boldsymbol{\lambda}\in\overline{\mathfrak a^+}$ for its chamber-valued Lyapunov spectrum and $\lambda_\alpha=\omega_\alpha(\boldsymbol{\lambda})$ for its fundamental-weight coordinates.
\end{convention}

\section{Finite-area Poincar\'e--Lelong identity}
\label{sec:PL-identity}

Let $W$ be a semisimple flat bundle with non-expanding cusp monodromy, equipped with a tame harmonic metric $h$ adapted to its Deligne extension $\mathcal W_\bullet$. Let $\mathcal{A}_\bullet\subset\mathcal W_\bullet^\vee$ be a filtered holomorphic line subbundle. Here and below, filtered holomorphic subbundles are saturated and carry the filtration induced from the ambient bundle. Write $\mathcal A$ for its restriction to $X$ and $h_{\mathcal A}$ for the metric induced from $h^\vee$. By the definition of the metric extension, $\Xi_{h_{\mathcal A}}(\mathcal A)=\mathcal A_\bullet$.

We use the parabolic Chern--Weil formula given by \cite{Simpson90}. Since the parabolic structure of $\mathcal{A}$ is induced by $\mathcal{W}_\bullet$, a similar argument of \cite[Lemma~10.5]{Simpson88} implies that the first Chern class $c_1(\mathcal A,h_{\mathcal A})$ is absolutely integrable. And then \cite[Lemma~6.2]{Simpson90} implies that
\begin{equation}
 \int_Xc_1(\mathcal{A},h_{\mathcal{A}})
 =\operatorname{pardeg}(\mathcal{A}_\bullet).
\label{eq:A-CW}
\end{equation}
Furthermore, since $p_*\nu=\mu_X$ and $\operatorname{Area}(X)=\pi\deg\KK_{\overline X}(D)/2$, the definition of $T_\nu$ gives
\[
T_\nu\left(p^*c_1(\mathcal{A},h_{\mathcal{A}})\right)=\frac{2\operatorname{pardeg}(\mathcal{A}_\bullet)}{\deg\KK_{\overline{X}}(D)}.
\]

After flat trivialization over $\widetilde X$, the line is represented by a $\rho_W^\vee$-equivariant holomorphic map
\[
f_{\mathcal{A}}\colon\widetilde X\longrightarrow\mathbb P(W^\vee).
\]
For a projective covector $[a]$, write $\ker[a]:=\ker a$, which does not depend on the nonzero representative $a$.
\begin{assumption}
\label{ass:no-common-kernel}
The projective map $f_{\mathcal{A}}$ satisfies
\[
\bigcap_{z\in\widetilde X}\ker f_{\mathcal{A}}(z)=\{0\}.
\]
\end{assumption}
Assumption~\ref{ass:no-common-kernel} is automatic when $W$ is irreducible.

Let $p\colon\mathbb P(W)\to X$ be the flat projective bundle. Evaluation defines a leafwise holomorphic section $s_{\mathcal{A}}$ of
\[
\mathcal N_{\mathcal{A}}
=\mathcal O_{\mathbb P(W)}(1)\otimes p^*\mathcal{A}^\vee,
\]
whose zero divisor is
\[
\Sigma_{\mathcal{A}}=\{(z,[u]):u\in\ker f_{\mathcal{A}}(z)\}.
\]
For nonzero representatives $a\in\mathcal{A}_z$ and $u\in W_z$, with $[a]=f_{\mathcal{A}}(z)$, its squared norm is
\[
\|s_{\mathcal{A}}\|^2
=\frac{|a(u)|^2}{\|a\|^2\|u\|^2}.
\]
Thus $0\leqslant \|s_{\mathcal{A}}\|^2\leqslant1$, and $s_{\mathcal{A}}$ vanishes precisely on $\Sigma_{\mathcal{A}}$. Moreover, we have the equality of leafwise Chern forms with respect to suitable metrics:
\begin{equation}
c_1^{\mathcal F}(\mathcal N_{\mathcal{A}},h_{\mathrm{FS}}\otimes p^*h_{\mathcal{A}}^\vee)
=c_1^{\mathcal F}(\mathcal O_{\mathbb P(W)}(1),h_{\mathrm{FS}})-p^*c_1(\mathcal{A},h_{\mathcal{A}}).
\label{eq:evaluation-c1}
\end{equation}
Proposition~\ref{prop:leafwise-projective-curvature} and \eqref{eq:A-CW} show that the coefficient of this form with respect to $p^*\omega_{\mathrm{hyp}}$ is $L^1(\nu)$ for every normalized harmonic measure $\nu$.

Assumption~\ref{ass:no-common-kernel} says exactly that no flat leaf is contained in $\Sigma_{\mathcal{A}}$. Indeed, such a leaf would be defined by a nonzero flat vector $u$ lying in $\ker f_{\mathcal{A}}(z)$ for every $z$.
Under this assumption, in a foliated chart the restriction of $s_{\mathcal{A}}$ to each local leaf has discrete zeros, which implies that
\begin{equation}
\nu(\Sigma_{\mathcal{A}})=0.
\label{eq:zero-divisor-measure}
\end{equation}

\subsection{Two lemmata}

\begin{lemma}
\label{lem:bounded-potential-stationarity}
Let $\nu$ be a normalized harmonic measure. If $\phi$ is bounded and leafwise $C^2$, and if $L_{\mathcal F}\phi\geqslant-g$ for some nonnegative $g\in L^1(\nu)$, then $L_{\mathcal F}\phi\in L^1(\nu)$ and
\begin{equation}
\int L_{\mathcal F}\phi\,\mathrm{d}\nu=0.
\label{eq:bounded-stationarity}
\end{equation}
\end{lemma}

\begin{proof}
Set $h=L_{\mathcal F}\phi$. Let $(Y_t)$ be leafwise Brownian motion started with law $\nu$. Since $\nu$ is harmonic, $Y_t$ has law $\nu$ for every $t$, and $h^-\leqslant g$ gives
\[
\mathbf E\!\int_0^t h^-(Y_s)\,\mathrm{d}s
=t\int h^-\,\mathrm{d}\nu
\leqslant t\int g\,\mathrm{d}\nu<\infty.
\]
Localize Dynkin's formula by stopping times $\tau_n\uparrow\infty$ for which the martingale part is integrable. Then
\[
\mathbf E\!\int_0^{t\wedge\tau_n}h^+(Y_s)\,\mathrm{d}s
\leqslant 2\|\phi\|_\infty+t\int g\,\mathrm{d}\nu.
\]
Monotone convergence and stationarity therefore imply $h^+\in L^1(\nu)$, so $h\in L^1(\nu)$. We can now remove the localization in Dynkin's formula by boundedness of $\phi$ and integrability of $h$, which give
\[
0=\mathbf E[\phi(Y_t)-\phi(Y_0)]
=t\int h\,\mathrm{d}\nu,
\]
which proves \eqref{eq:bounded-stationarity}.
\end{proof}

\begin{lemma}
\label{lem:regularized-PL}
Suppose that Assumption~\ref{ass:no-common-kernel} holds.
In this lemma, we use $D_h$ to denote the leafwise Chern connection of $\mathcal N_{\mathcal{A}}$ with its induced Hermitian metric. For $\varepsilon>0$ set $\phi_\varepsilon=\log(\|s_{\mathcal{A}}\|^2+\varepsilon)$. Then
\begin{equation}
\begin{aligned}
\mathrm{d}_{\mathcal F}\mathrm{d}^c_{\mathcal F}\phi_\varepsilon
&=-\frac{\|s_{\mathcal{A}}\|^2}{\|s_{\mathcal{A}}\|^2+\varepsilon}
 c_1^{\mathcal F}(\mathcal N_{\mathcal{A}},h_{\mathrm{FS}}\otimes p^*h_{\mathcal{A}}^\vee)
+\Theta_\varepsilon,\\
\Theta_\varepsilon
&=\frac{\sqrt{-1}}{2\pi}
\frac{\varepsilon}{(\|s_{\mathcal{A}}\|^2+\varepsilon)^2}
\langle D_h^{1,0}s_{\mathcal{A}},D_h^{1,0}s_{\mathcal{A}}\rangle
\geqslant0,
\end{aligned}
\label{eq:regularized-PL}
\end{equation}
where the bracket denotes the scalar $(1,1)$-form induced by the Hermitian metric on $\mathcal N_{\mathcal{A}}$. If locally we can express $s_{\mathcal{A}}=fe$ in a holomorphic frame $e$ on a local leaf, then, for every compactly supported continuous test function $\chi$ on that local leaf,
\begin{equation}
\lim_{\varepsilon\downarrow0}\int\chi\,\Theta_\varepsilon
=
\sum_{f(z)=0}\operatorname{ord}_z(f)\,\chi(z).
\label{eq:leafwise-zero-convergence}
\end{equation}
\end{lemma}

\begin{proof}
In a foliated chart, choose a holomorphic frame $e$ of $\mathcal N_{\mathcal{A}}$ on a local leaf, and write $s_{\mathcal{A}}=fe$, $\|e\|^2=\mathrm{e}^{-\psi}$, and $\|s_{\mathcal{A}}\|^2=|f|^2\mathrm{e}^{-\psi}$. If $\alpha=\partial f-f\,\partial\psi$, then $D_h^{1,0}s_{\mathcal{A}}=\alpha\otimes e$ and a direct calculation gives
\[
\partial\bar\partial\log(\|s_{\mathcal{A}}\|^2+\varepsilon)
=-\frac{\|s_{\mathcal{A}}\|^2}{\|s_{\mathcal{A}}\|^2+\varepsilon}\,\partial\bar\partial\psi
+\frac{\varepsilon}{(\|s_{\mathcal{A}}\|^2+\varepsilon)^2}
 \mathrm{e}^{-\psi}\alpha\wedge\bar\alpha.
\]
Since $c_1^{\mathcal F}(\mathcal N_{\mathcal{A}})=(\sqrt{-1}/2\pi)\partial\bar\partial\psi$, this follows from \eqref{eq:regularized-PL}. Note that its second term is positive.

By Assumption~\ref{ass:no-common-kernel}, $f$ is not identically zero on this local leaf, so $\log \|s_{\mathcal{A}}\|^2\in L^1_{\mathrm{loc}}$. Moreover $\log(\|s_{\mathcal{A}}\|^2+\varepsilon)\to\log \|s_{\mathcal{A}}\|^2$ in $L^1_{\mathrm{loc}}$, while $\frac{\|s_{\mathcal{A}}\|^2}{\|s_{\mathcal{A}}\|^2+\varepsilon}\to1$ in $L^1_{\mathrm{loc}}$. Hence, after rearranging \eqref{eq:regularized-PL},
\[
\Theta_\varepsilon\longrightarrow
\mathrm{d}_{\mathcal F}\mathrm{d}^c_{\mathcal F}\log \|s_{\mathcal{A}}\|^2
+c_1^{\mathcal F}(\mathcal N_{\mathcal{A}},h_{\mathrm{FS}}\otimes p^*h_{\mathcal{A}}^\vee)
=[Z(f)]
\]
as currents, by the Poincar\'e--Lelong formula. Since the $\Theta_\varepsilon$ are positive measures, current convergence is weak convergence against $C_c^0$ test functions, which is exactly \eqref{eq:leafwise-zero-convergence}.
\end{proof}

Under Assumption~\ref{ass:no-common-kernel}, the negative part of $L_{\mathcal F}\phi_\varepsilon$ is dominated by an $L^1(\nu)$ function, so \eqref{eq:bounded-stationarity} gives
\[
T_\nu(\Theta_\varepsilon)
=T_\nu\!\left(
\frac{\|s_{\mathcal{A}}\|^2}{\|s_{\mathcal{A}}\|^2+\varepsilon}c_1^{\mathcal F}(\mathcal N_{\mathcal{A}},h_{\mathrm{FS}}\otimes p^*h_{\mathcal{A}}^\vee)
\right).
\]
Using \eqref{eq:zero-divisor-measure} and dominated convergence, define
\begin{equation}
\delta_{\mathcal{A}}(\nu)
:=\lim_{\varepsilon\downarrow0}T_\nu(\Theta_\varepsilon)
=T_\nu\!\left(c_1^{\mathcal F}(\mathcal N_{\mathcal{A}},h_{\mathrm{FS}}\otimes p^*h_{\mathcal{A}}^\vee)\right)
\geqslant0.
\label{eq:PL-positive-term}
\end{equation}

\subsection{Gap between Lyapunov exponents and Harder--Narasimhan slopes}

\begin{theorem}
\label{thm:main-A}
Under Assumption~\ref{ass:no-common-kernel}, every normalized harmonic probability measure $\nu$ on $\mathbb P(W)$ satisfies
\[
\lambda(\nu)
-\frac{2\operatorname{pardeg}(\mathcal{A}_\bullet)}{\deg\KK_{\overline{X}}(D)}
=\delta_{\mathcal{A}}(\nu)\geqslant0.
\]
\end{theorem}

\begin{proof}
Pair \eqref{eq:evaluation-c1} with $T_\nu$. The tautological term is $\lambda(\nu)$ by \eqref{eq:harmonic-current}, and the term contributed by $p^*c_1(\mathcal{A},h_{\mathcal{A}})$ is $2\operatorname{pardeg}(\mathcal{A}_\bullet)/\deg\KK_{\overline{X}}(D)$ by \eqref{eq:A-CW} and $p_*\nu=\mu_X$. Now the statement follows from \eqref{eq:PL-positive-term}.
\end{proof}

We recover the following consequence of the lower bound of Eskin--Kontsevich--M\"oller--Zorich \cite[Theorem~4.1]{EKMZ18} in the spirit of \cite{DanielDeroin19}.

\begin{corollary}
\label{cor:HN-Lyap}
Let $(V,\nabla)$ be a semisimple flat bundle of rank $r$ with non-expanding cusp monodromy, and let $\mathcal V_\bullet$ be the canonical filtered Deligne extension of its de Rham holomorphic bundle $\mathcal V=(V,\nabla^{0,1})$. Then
\[
\frac{2}{\deg\KK_{\overline{X}}(D)}
\sum_{i=1}^k\mu_i(\mathcal V_\bullet)
\leqslant\sum_{i=1}^k\lambda_i(V),
\qquad 1\leqslant k\leqslant r,
\]
with equality for $k=r$. Equivalently,
\[
\frac{2}{\deg\KK_{\overline{X}}(D)}
\boldsymbol{\mu}_{\mathrm{HN}}(\mathcal V_\bullet)
\preceq\boldsymbol{\lambda}(V),
\]
where $\preceq$ denotes dominance of decreasing vectors: all partial sums are bounded above and the total sums agree.
\end{corollary}

\begin{proof}
Let $\mathcal F_\bullet$ be an HN step of rank $k$ and put $\mathcal A_\bullet=\det\mathcal F_\bullet$. The smallest flat subbundle $U\subset\bigwedge^kV$ containing $\mathcal A$ is semisimple. Its fibers are spanned by the flat translates of $\mathcal A$, so $W=U^\vee$ satisfies Assumption~\ref{ass:no-common-kernel}.
Exactness of filtered Deligne extensions preserves the induced filtration $\mathcal A_\bullet\subset\mathcal U_\bullet$. Equip $U$ with an adapted tame harmonic metric and $\mathcal A$ with the induced metric. Thus $\mathcal A_\bullet\subset\mathcal W_\bullet^\vee$ is a filtered line subbundle to which Theorem~\ref{thm:main-A} applies.

Applying Theorem~\ref{thm:main-A} to a Lyapunov-realizing harmonic measure on $\mathbb P(W)$, supplied by Proposition~\ref{prop:lyap-coincidence}\romitem{2}, we obtain that
\[
\frac{2\operatorname{pardeg}\mathcal F_\bullet}{\deg\KK_{\overline{X}}(D)}\leqslant\lambda_1(U^\vee)
=\lambda_1(U)\leqslant\lambda_1\!\left(\bigwedge^kV\right)=\sum_{i=1}^k\lambda_i(V).
\]
Here the first equality uses duality and geodesic reversibility. Since $\operatorname{pardeg}\mathcal F_\bullet=\sum_{i=1}^k\mu_i(\mathcal V_\bullet)$ at an HN rank, this proves the required inequalities there. Other remaining ranks then follow directly by affinity.
Finally, $\operatorname{pardeg}\mathcal V_\bullet=0$ by the degree--residue formula, and $\sum_i\lambda_i(V)=0$ by reversibility, giving equality at full rank.
\end{proof}

\begin{corollary}\label{cor:equality-all-harmonic}
	Suppose that Assumption~\ref{ass:no-common-kernel} holds. If
	\[\lambda_1(W)=\frac{2\operatorname{pardeg}(\mathcal{A}_\bullet)}{\deg\KK_{\overline{X}}(D)},\]
	then every normalized harmonic measure $\nu$ on $\mathbb P(W)$ is Lyapunov-realizing and
	\[\lambda(\nu)=\lambda_1(W),\qquad\delta_{\mathcal{A}}(\nu)=0.\]
\end{corollary}

\begin{proof}
	Theorem~\ref{thm:main-A} gives $\lambda(\nu)\geqslant\lambda_1(W)$, while \cite[Remark~2.17]{DanielDeroin19} gives $\lambda(\nu)\leqslant\lambda_1(W)$.
\end{proof}

\subsection{Criterion of equality}

\begin{theorem}\label{thm:strict-zero-divisor}
	Suppose that Assumption~\ref{ass:no-common-kernel} holds. For every normalized harmonic measure $\nu$,
	\begin{equation}\label{eq:strict-zero-divisor}
		\operatorname{supp}\nu\cap\Sigma_{\mathcal{A}}\neq\varnothing\quad\Longrightarrow\quad\delta_{\mathcal{A}}(\nu)>0.
	\end{equation}
	Hence $\delta_{\mathcal{A}}(\nu)=0$ implies $\operatorname{supp}\nu\cap\Sigma_{\mathcal{A}}=\varnothing$.
\end{theorem}

\begin{proof}
	Choose $y_0\in\operatorname{supp}\nu\cap\Sigma_{\mathcal{A}}$ and a relatively compact foliated chart $B\cong\mathbb D_R\times T$ with $y_0=(0,t_0)$. After a leafwise holomorphic trivialization of $\mathcal N_{\mathcal{A}}$, write $s_{\mathcal{A}}=f_t(z)e$, where $f_t$ is holomorphic in $z$ and depends continuously on $(z,t)$.
	Since no leaf is contained in $\Sigma_{\mathcal{A}}$, $f_{t_0}\not\equiv0$, while $f_{t_0}(0)=0$. Hence one may choose $0<r<r_1<R$ such that $f_{t_0}$ has no zero on $\partial\mathbb D_r$ and has $m\geqslant1$ zeros in $\mathbb D_r$, counted with multiplicity.
	By Rouch\'e's theorem, after shrinking to a transverse neighborhood $U\ni t_0$, every $f_t$, $t\in U$, has exactly $m$ zeros in $\mathbb D_r$.
	
	By the local disintegration theorem for harmonic measures \cite{Garnett83}, on $\mathbb D_{r_1}\times U$ we may write
	\[\mathrm{d}\nu=h_t(z)\,\mathrm{d}A_{\mathrm{hyp}}(z)\,\mathrm{d}\eta(t),\]
	with $h_t$ nonnegative harmonic for $\eta$-almost every $t$. Harnack's inequality on the fixed pair $\mathbb D_r\Subset\mathbb D_{r_1}\Subset\mathbb D_R$ gives a constant $c>0$, independent of $t$, such that
	\[h_t(z)\geqslant c\int_{\mathbb D_{r_1}}h_t\,\mathrm{d}A_{\mathrm{hyp}}\qquad (z\in\mathbb D_r)\]
	for $\eta$-almost every $t\in U$.
	
	Choose $\chi\in C_c(\mathbb D_{r_1})$ with $0\leqslant\chi\leqslant1$ and $\chi\equiv1$ on $\mathbb D_r$. For $\eta$-almost every $t\in U$, Lemma~\ref{lem:regularized-PL}, applied to the continuous test function $\chi h_t$, gives
	\[\lim_{\varepsilon\downarrow0}\int_{\mathbb D_{r_1}}\chi h_t\,\Theta_{\varepsilon,t}=\sum_{\substack{z\in\mathbb D_{r_1}\\ f_t(z)=0}}\operatorname{ord}_z(f_t)\,\chi(z)h_t(z)\geqslant mc\int_{\mathbb D_{r_1}}h_t\,\mathrm{d}A_{\mathrm{hyp}}.\]
Define the localized mass
\[I_\varepsilon:=\pi\int_U\left(\int_{\mathbb D_{r_1}}\chi h_t\,\Theta_{\varepsilon,t}\right)\mathrm{d}\eta(t).\]
Since $\Theta_\varepsilon\geqslant0$, one has $0\leqslant I_\varepsilon\leqslant T_\nu(\Theta_\varepsilon)$. Therefore, by Fatou's lemma,
\[\begin{aligned}
	\delta_{\mathcal{A}}(\nu)
	&\geqslant\liminf_{\varepsilon\downarrow0}I_\varepsilon\\
	&\geqslant\pi\int_U\liminf_{\varepsilon\downarrow0}\left(\int_{\mathbb D_{r_1}}\chi h_t\, \Theta_{\varepsilon,t}\right)\mathrm{d}\eta(t)\\
	&\geqslant\pi mc\int_U\int_{\mathbb D_{r_1}}h_t\, \mathrm{d}A_{\mathrm{hyp}}\,\mathrm{d}\eta(t)\\
	&=\pi mc\,\nu(\mathbb D_{r_1}\times U)>0.
\end{aligned}\]
The last inequality holds because $\mathbb D_{r_1}\times U$ is a neighborhood of $y_0\in\operatorname{supp}\nu$. This proves \eqref{eq:strict-zero-divisor}.
\end{proof}

Under Assumption~\ref{ass:no-common-kernel}, the preceding two results give
\[
\lambda_1(W)
=\frac{2\operatorname{pardeg}(\mathcal{A}_\bullet)}{\deg\KK_{\overline{X}}(D)}
\quad\Longrightarrow\quad
\operatorname{supp}\nu\cap\Sigma_{\mathcal{A}}=\varnothing
\quad\text{for every normalized harmonic }\nu.
\]

For comparison with the compact theory, put $\mathcal F=\ker(\mathcal W\to\mathcal{A}^\vee)$. Then $\mathbb P(\mathcal F)=\Sigma_{\mathcal{A}}$ and
\[
\mathcal O_{\mathbb P(\mathcal W)}([\mathbb P(\mathcal F)])
\cong
\mathcal O_{\mathbb P(\mathcal W)}(1)\otimes p^*\mathcal{A}^\vee
=\mathcal N_{\mathcal{A}}.
\]
Thus, on a compact base, $\delta_{\mathcal{A}}(\nu)$ is exactly the dynamical degree of Daniel--Deroin. Their dynamical-degree formula and geometric-intersection theorem give both nonnegativity and
\[
\delta_{\mathcal{A}}(\nu)=0
\iff
\operatorname{supp}\nu\cap\mathbb P(\mathcal F)=\varnothing
\]
when no leaf germ is contained in the divisor, see \cite[\S3.2--3.3]{DanielDeroin19}. Costantini's compact asymptotic degree agrees with the same dynamical degree whenever it is defined for the corresponding developing datum, see \cite{Costantini20}. Theorem \ref{thm:main-A} is the finite-area parabolic analytic version of this identity. On a noncompact base, its geometric-intersection interpretation requires cusp control \cite[\S3.5]{DanielDeroin19}.

\subsection{Principal Harder--Narasimhan reduction}
In this subsection we assume that $X=\overline{X}$ is compact, so $D=\varnothing$, and that the monodromy factors through a connected semisimple group $G_{\mathbb R}$ with complexification $G$.
Denote the holomorphic de Rham principal $G$-bundle by $\mathbb E$. For a $G$-representation $V$, write $\mathbb E(V):=\mathbb E\times_G V$ for its associated holomorphic vector bundle.

Choose a Cartan of $\mathfrak{g}_{\mathbb R}$ containing $\mathfrak a$ and stable under the Cartan involution, with positive complex roots compatible with the chosen positive restricted roots. This determines a maximal torus $T\subset G$ and an inclusion $\iota\colon\mathfrak a\hookrightarrow\mathfrak h_{\mathbb R}$ into the real span of the complex coroots.

By \cite[Theorem~4.1]{BiswasHolla04}, $\mathbb E$ has a unique Harder--Narasimhan reduction. The dominant type $\boldsymbol{\mu}_{\mathrm{HN}}(\mathbb E)\in\mathfrak h_{\mathbb R}$ of the holomorphic HN reduction is determined by the degrees of its characters and has zero pairing with the simple roots of its Levi factor.
In the vector-bundle case, this is the ordered list of HN slopes, repeated with multiplicity, as recalled in Section~\ref{sec:setup}.

\begin{proposition}
\label{prop:principal-HN-Lyap}
For compact $X$,
\begin{equation}
\frac{2}{\deg\KK_X}\,\boldsymbol{\mu}_{\mathrm{HN}}(\mathbb E)
\preceq
\iota(\boldsymbol{\lambda})
\label{eq:principal-HN-Lyap}
\end{equation}
in the dominant chamber of the complexified group. For split $G_{\mathbb R}$, the two chambers identify via $\iota$.
\end{proposition}

\begin{proof}
Write $\mu=\boldsymbol{\mu}_{\mathrm{HN}}(\mathbb E)$ and let $\mathbb E_P$ be this reduction to the standard parabolic $P$. The roots of its Levi factor satisfy $\alpha(\mu)=0$, whereas those of its unipotent radical satisfy $\alpha(\mu)>0$.
For every algebraic character $\chi$ of $P$, our convention for the HN type is
\[
\deg(\mathbb E_P\times_P\mathbb C_\chi)=\chi(\mu),
\]
where $\mathbb C_\chi$ is the one-dimensional $P$-module with character $\chi$.

Fix a fundamental dominant weight $\varpi$. Since the character lattice of $T$ has finite index in the weight lattice, choose $m\geqslant1$ such that $m\varpi$ is a character of $T$, and let $\tau\colon G\to\mathrm{GL}(V)$ be the irreducible representation with highest weight $m\varpi$. Every weight $\beta$ of $V$ satisfies $m\varpi-\beta=\sum_\alpha n_\alpha\alpha$, with $n_\alpha\geqslant0$ and $\alpha$ ranging over the simple complex roots. Since $\mu$ is dominant, the largest value of $\beta(\mu)$ is $m\varpi(\mu)$. Define
\[
F=\bigoplus_{\beta(\mu)=m\varpi(\mu)}V_\beta,
\qquad k=\dim F.
\]
A root operator for $\alpha$ maps $V_\beta$ to $V_{\beta+\alpha}$. Thus the Levi preserves $F$, while every root operator of the unipotent radical annihilates $F$, since it would strictly increase the maximal $\mu$-weight. Since $P$ is generated by $T$ and these root subgroups, $F$ is $P$-invariant. Consequently $\mathcal F=\mathbb E_P\times_P F\subset\mathbb E(V)$ is a holomorphic subbundle of rank $k$.

Let $\chi$ be the character of $P$ acting on $\det F$.
Its restriction to $T$ is the sum of the weights of $F$, counted with multiplicity. Since $\det\mathcal F=\mathbb E_P\times_P\mathbb C_\chi$, we obtain
\[
\deg\mathcal F=\chi(\mu)=\sum_{\beta(\mu)=m\varpi(\mu)}(\dim V_\beta)\,\beta(\mu)=km\varpi(\mu).
\]
On the compact curve $X$, any smooth Hermitian metric makes the flat cocycle on $\mathbb E(V)$ integrable. Applying \cite[Theorem~4.1]{EKMZ18} to $\mathcal F$ gives
\[
\frac{2m}{\deg\KK_X}\,\varpi(\mu)\leqslant\frac1k\sum_{j=1}^k\lambda_j(V)\leqslant\lambda_1(V).
\]

To identify the last term, choose a Hermitian norm compatible with the Cartan involution, so that $\tau(K_{\mathbb R})$ is unitary and $\mathrm{d}\tau(H)$ is self-adjoint for $H\in\mathfrak a$.
Compatibility of the positive systems gives $\iota(\overline{\mathfrak a^+})\subset\overline{\mathfrak h_{\mathbb R}^+}$, where the latter is the closed complex dominant chamber. For $g\in G_{\mathbb R}$, the Cartan decomposition and the weight decomposition therefore imply
\[
\log\|\tau(g)\|_{\mathrm{op}}
=\max_\beta\beta\left(\iota(\kappa(g))\right)
=m\varpi\left(\iota(\kappa(g))\right).
\]
Applying this identity to the flat cocycle and using \eqref{eq:G-Lyap-spectrum} yields $\lambda_1(V)=m\varpi\left(\iota(\boldsymbol{\lambda})\right)$.
Dividing the preceding inequality by $m$ proves the required inequality for each fundamental dominant weight. These inequalities characterize dominance order, and hence give \eqref{eq:principal-HN-Lyap}.
\end{proof}

\begin{remark}
We restrict Proposition~\ref{prop:principal-HN-Lyap} to closed surfaces so that ordinary holomorphic principal HN theory applies. For punctured surfaces, there should be an analogous formulation using a principal HN reduction compatible with the boundary filtrations and defined using parabolic degrees. For the existence and uniqueness of the principal HN reduction when the parabolic weights are rational, see \cite[\S6]{DeyParthasarathi05} and \cite[Remark~6.3.8]{BalajiSeshadri15}; for parahoric torsors with arbitrary real weights and their stability theory, see \cite[\S\S6--7]{BalajiBiswasPandey17}.
\end{remark}

\section{Transversality and Hermitian domains}
\label{sec:hermitian-transversality}

Let $G_{\mathbb R}$ be a connected semisimple real Lie group with finite center and without compact factors. Throughout this section, let \[\eta\colon G_{\mathbb R}\longrightarrow G_{\mathbb C}\] be its universal complexification. The group $G_{\mathbb C}$ is connected complex semisimple, with Lie algebra $\operatorname{Lie}(G_{\mathbb R})\otimes_{\mathbb R}\mathbb C$. 
Every finite-dimensional complex representation of $G_{\mathbb R}$ extends uniquely to a holomorphic, hence algebraic, representation of $G_{\mathbb C}$. We pull back $G_{\mathbb C}$-actions along $\eta$. A subgroup $\Gamma<G_{\mathbb R}$ is called \emph{Zariski dense} if $\eta(\Gamma)$ is Zariski dense over $\mathbb C$ in $G_{\mathbb C}$. 

A \emph{lattice} in $G_{\mathbb R}$ is a discrete subgroup of finite covolume. For the thinness specializations, by an \emph{arithmetic lattice} $\Gamma_{\mathrm{ar}}<G_{\mathbb R}$ we mean, up to finite central isogeny, a lattice arising from the integral points of a $\mathbb Q$-form: there are a connected semisimple algebraic group $\mathbf G$ over $\mathbb Q$ and a homomorphism $\iota\colon\mathbf G(\mathbb R)^\circ\to G_{\mathbb R}$ with compact kernel and finite-index image such that $\Gamma_{\mathrm{ar}}$ is commensurable with $\iota(\mathbf G(\mathbb Z)\cap\mathbf G(\mathbb R)^\circ)$. A Zariski-dense subgroup of an arithmetic lattice is called \emph{thin} if it has infinite index.
A subgroup of a lattice is itself a lattice if and only if its index is finite. Thus, under arithmetic containment and Zariski density, being a nonlattice is equivalent to thinness.

Denote by $\rho\colon\pi_1(X)\to G_{\mathbb R}$ the monodromy and $\Gamma:=\rho(\pi_1(X))$. For a representation on $W$, write $\mathcal W$ for the associated complex flat bundle, using $W_{\mathbb C}=W\otimes_{\mathbb R}\mathbb C$ when $W$ is real.

Fix a maximal split Cartan and a positive system, with simple restricted roots $\Delta$ and fundamental restricted weights $\omega_\alpha$. Let $\varnothing\neq\theta\subset\Delta$. By \cite[Lemma~3.2]{GGKW17}, applied to the adjoint group and pulled back to $G_{\mathbb R}$, there is an integer $N\geqslant1$ and an irreducible proximal real linear representation $\tau\colon G_{\mathbb R}\to\mathrm{GL}(W)$ with highest restricted weight
\[\chi=N\sum_{\alpha\in\theta}\omega_\alpha.\]
Its highest restricted-weight space is one-dimensional, and the support of $\chi$ in the fundamental weights is exactly $\theta$. Thus $\tau$ is $\theta$-proximal in the sense of \cite[Definition~3.1]{GGKW17}. Moreover, $W_{\mathbb C}$ is irreducible, otherwise $W$ would admit an invariant complex structure preserving its one-dimensional highest restricted-weight space. We use the same symbol $\tau$ for its algebraic extension to $G_{\mathbb C}$ on $W_{\mathbb C}$.

Let $0\neq v_+\in W$ span the highest restricted-weight line. By \cite[Proposition~3.3(a)]{GGKW17}, the stabilizer of the line $\mathbb Rv_+$ is $P_\theta$, and \cite[Proposition~3.3(b)]{GGKW17} gives the equivariant embedding
\[\begin{aligned}
	\iota^+\colon G_{\mathbb R}/P_\theta&\longrightarrow\mathbb P(W),\\
	gP_\theta&\longmapsto[\tau(g)v_+].
\end{aligned}\]
Following the terminology used in \cite[Section~2.3]{BCLS15}, we call $\tau$ a \emph{Pl\"ucker representation} for $P_\theta$.

Let $P_{\theta,\mathbb C}:=\operatorname{Stab}_{G_{\mathbb C}}([v_+])$. We denote the complex flag orbit and its real flag orbit by
\[\mathcal F_{\mathbb C}(\tau):=G_{\mathbb C}\cdot[v_+]\cong G_{\mathbb C}/P_{\theta,\mathbb C}\subset\mathbb P(W_{\mathbb C}),
\qquad
\mathcal F_{\mathbb R}(\tau):=\iota^+(G_{\mathbb R}/P_\theta)\subset\mathcal F_{\mathbb C}(\tau).\]
Choose $0\neq\alpha_-\in W_{\mathbb C}^\vee$ that vanishes on every restricted-weight space other than $\mathbb Cv_+$. The stabilizer $P_{\theta,\mathbb C}^{-}$ of $[\alpha_-]$ is opposite to $P_{\theta,\mathbb C}$. The dual flag orbit is
\[\mathcal F_{\mathbb C}^{\vee}(\tau)
=G_{\mathbb C}\cdot[\alpha_-]
\cong G_{\mathbb C}/P_{\theta,\mathbb C}^{-}
\subset\mathbb P(W_{\mathbb C}^{\vee}).\]
For $[v]\in\mathcal F_{\mathbb C}(\tau)$ and $[\alpha]\in\mathcal F_{\mathbb C}^{\vee}(\tau)$, we say that $[v]$ and $[\alpha]$ are \emph{transverse} if the pair lies in the unique open $G_{\mathbb C}$-orbit in $\mathcal F_{\mathbb C}(\tau)\times\mathcal F_{\mathbb C}^{\vee}(\tau)$. Equivalently, their stabilizer parabolics are opposite. In the above Pl\"ucker embeddings this is exactly the projective non-vanishing condition $\alpha(v)\neq0$. See \cite[Proposition~3.3]{GGKW17}.

For an arbitrary holomorphic line subbundle $\mathcal{A}\subset\mathcal W^\vee$, pulling back to $\widetilde X$ and using the flat trivialization gives a holomorphic projective map
\[\begin{aligned}
f_{\mathcal{A}}\colon\widetilde X&\longrightarrow\mathbb P(W_{\mathbb C}^\vee),\\
z&\longmapsto[\mathcal{A}_z].
\end{aligned}\]
We call this map \emph{flag-valued} (for $\tau$) if $f_{\mathcal{A}}(\widetilde X)\subset \mathcal F_{\mathbb C}^{\vee}(\tau)$.
In this case we regard $f_{\mathcal{A}}$ as a map to $\mathcal F_{\mathbb C}^{\vee}(\tau)$.

For arbitrary $[v]\in\mathbb P(W_{\mathbb C})$ and $[\alpha]\in\mathbb P(W_{\mathbb C}^\vee)$, we use the same transversality criterion $\alpha(v)\neq0$. Following \cite{AlessandriniSaracco11}, for a subset $S\subset\mathbb P(W_{\mathbb C})$, we say its \emph{dual complement} is
\[S^*:=\{[\alpha]\in\mathbb P(W_{\mathbb C}^\vee):\mathbb P(\ker\alpha)\cap S=\varnothing\}.\]
In particular, set
\[\Omega_{\mathrm{tr}}(\tau):=\mathcal F_{\mathbb R}(\tau)^*.\]
Thus $[\alpha]\in\Omega_{\mathrm{tr}}(\tau)$ exactly when it is transverse to every point of $\mathcal F_{\mathbb R}(\tau)$. These definitions are independent of the choice of nonzero representatives. For any complex representation $W$, we use the same notation $S^*\subset\mathbb P(W^\vee)$ for the dual complement of $S\subset\mathbb P(W)$.

\subsection{Transversality and lattice monodromy}

\begin{theorem}
\label{thm:contact-limit-set}
Let $W$ be an irreducible complex representation of $G_{\mathbb C}$. Assume that $\Gamma$ is Zariski dense, and let $\mathcal{A}_\bullet\subset\mathcal W_\bullet^\vee$ be a filtered holomorphic line subbundle, with restriction $\mathcal A$ to $X$. Suppose that
\[\lambda_1(W)=\frac{2\operatorname{pardeg}\mathcal{A}_\bullet}{\deg\KK_{\overline{X}}(D)}.\]
Then, for every nonempty minimal closed $\Gamma$-invariant subset $\Lambda\subset\mathbb P(W)$,
\[f_{\mathcal{A}}(\widetilde X)\subset\Lambda^*.\]
Equivalently, $v\notin\ker f_{\mathcal{A}}(z)$ for every $z\in\widetilde X$ and $[v]\in\Lambda$.
\end{theorem}

\begin{proof}
Since $G_{\mathbb C}$ is connected and $\eta(\Gamma)$ is Zariski dense, the image of every finite-index subgroup of $\Gamma$ is Zariski dense. Irreducibility of $W$ therefore makes its restriction to $\Gamma$ strongly irreducible. The common kernel $\bigcap_{z\in\widetilde X}\ker f_{\mathcal{A}}(z)$ is a proper $\Gamma$-invariant complex subspace of $W$, so it vanishes. Thus Assumption~\ref{ass:no-common-kernel} holds.

Put $M_\Lambda:=\widetilde X\times_\rho\Lambda$. Proposition~\ref{prop:lyap-coincidence}\romitem{2} supplies a harmonic probability measure $\nu$ supported on $M_\Lambda$, with $p_*\nu=\mu_X$ and $\lambda(\nu)=\lambda_1(W)$. Its support is closed and saturated. Since $\mu_X$ has full support and the projective fibers are compact, its support meets every fiber.
The fiber over a base point is a nonempty closed $\Gamma$-invariant subset of $\Lambda$, so minimality gives $\operatorname{supp}\nu=M_\Lambda$. Theorem~\ref{thm:main-A} and the assumed equality imply $\delta_{\mathcal{A}}(\nu)=0$. Theorem~\ref{thm:strict-zero-divisor} then gives $M_\Lambda\cap\Sigma_{\mathcal{A}}=\varnothing$, which is precisely $f_{\mathcal{A}}(\widetilde X)\subset\Lambda^*$.
\end{proof}

\begin{proposition}\label{prop:lattice-minimality}
Let $W$ be an irreducible complex representation of $G_{\mathbb C}$ and let $\mathcal S\subset\mathbb P(W)$ be a compact $G_{\mathbb R}$-orbit. If $\Gamma<G_{\mathbb R}$ is a lattice, every $\Gamma$-orbit in $\mathcal S$ is dense. In particular, every lattice in $G_{\mathbb R}$ acts minimally on $\mathcal F_{\mathbb R}(\tau)$.
\end{proposition}

\begin{proof}
Let $G_{\mathbb R}=K_{\mathbb R}AN$ be the Iwasawa decomposition, where $A=\exp(\mathfrak{a})$. Choose an $H\in\mathfrak a^+$ such that different restricted weights of $W$ valued on $H$ differently. and let $\chi$ be the restricted weight such that $\chi(H)$ is maximal.
The complex span of $\mathcal S$ is $G_{\mathbb R}$-invariant, hence it is $G_{\mathbb C}$-invariant. Hence it equals $W$ by irreducibility. Therefore, there exists a $[v]\in\mathcal S$ whose component $v_\chi$ in the restricted-weight space $W_\chi$ is nonzero. Then the weight decomposition gives
\[\lim_{t\to+\infty}\exp(tH)\cdot[v]=[v_\chi].\]
Since $\mathcal S$ is compact, $[v_\chi]\in\mathcal S$. Note that the group $A$ acts by scalars on $W_\chi$ since $\chi(H)$ is maximal. Hence $N$ fixes $W_\chi$ pointwise. Therefore
\[AN\subset Q:=\operatorname{Stab}_{G_{\mathbb R}}([v_\chi]).\]

So $\mathcal{S}\cong G_{\mathbb R}/Q$. The closed subgroup $Q$ contains the maximal connected split solvable subgroup $AN$. Since $G_{\mathbb R}$ has no compact factors, the conditions in \cite[Corollary~1.7]{GorodnikOh07} are satisfied. Hence $\Gamma$-orbit in $G_{\mathbb R}/Q$ equidistributed with respect to its $K_{\mathbb R}$-invariant probability measure. This measure has full support because $G_{\mathbb R}=K_{\mathbb R}Q$. As a result, every such orbit is dense.
\end{proof}

\begin{corollary}\label{cor:lattice-transversality}
Let $W$ be an irreducible complex representation of $G_{\mathbb C}$ and let $\mathcal S\subset\mathbb P(W)$ be a compact $G_{\mathbb R}$-orbit. Let $\mathcal{A}_\bullet\subset\mathcal W_\bullet^\vee$ be a filtered holomorphic line subbundle, with restriction $\mathcal A$ to $X$. Suppose that
\[\lambda_1(W)=\frac{2\operatorname{pardeg}\mathcal{A}_\bullet}{\deg\KK_{\overline{X}}(D)}.\]
Then either
\[\Gamma\text{ is not a lattice},
\qquad\text{or}\qquad
f_{\mathcal{A}}(\widetilde X)\subset\mathcal S^*.\]
For a real Pl\"ucker representation $\tau$, apply this to its complexification and $\mathcal S=\mathcal F_{\mathbb R}(\tau)$, then $\mathcal S^*=\Omega_{\mathrm{tr}}(\tau)$.
\end{corollary}

\begin{proof}
If $\Gamma$ is a lattice, Borel density \cite{Borel60} gives Zariski density, and Proposition~\ref{prop:lattice-minimality} makes $\mathcal S$ a minimal $\Gamma$-invariant subset of $\mathbb P(W)$.
Theorem~\ref{thm:contact-limit-set}, applied with $\Lambda=\mathcal S$, gives $f_{\mathcal{A}}(\widetilde X)\subset\mathcal S^*$.
\end{proof}

We use the Cartan-projection characterization \cite{bochi2019anosov,kapovich2018morse} of Anosov subgroups. A finitely generated subgroup $\Gamma<G_{\mathbb R}$ is \emph{$P_\theta$-Anosov} if there are constants $c>0$ and $C\geqslant0$ such that
\[
\alpha\left(\kappa(\gamma)\right)
\geqslant c|\gamma|_\Gamma-C
\qquad\forall\gamma\in\Gamma,\ \alpha\in\theta.
\]
Here $|\cdot|_\Gamma$ is word length with respect to a finite symmetric generating set, and $\kappa\colon G_{\mathbb R}\to\overline{\mathfrak a^+}$ is the Cartan projection, characterized by $g\in K_{\mathbb R}\exp(\kappa(g))K_{\mathbb R}$.
The condition is independent of the generating set and is equivalent to the ordinary dynamical definition for the inclusion $\Gamma\hookrightarrow G_{\mathbb R}$.
Such a group is word hyperbolic and admits a continuous $\Gamma$-equivariant boundary map $\xi\colon\partial_\infty\Gamma\to G_{\mathbb R}/P_\theta$, where $\partial_\infty\Gamma$ denotes its Gromov boundary.

\begin{corollary}\label{cor:anosov-opposition}
Let $\tau\colon G_{\mathbb R}\to\mathrm{GL}(W)$ be the real Pl\"ucker representation fixed above. Assume the hypotheses of Theorem~\ref{thm:contact-limit-set} for its complexification $W_{\mathbb C}$, and suppose that $\Gamma<G_{\mathbb R}$ is $P_\theta$-Anosov. Let
\[\xi\colon\partial_\infty\Gamma\longrightarrow G_{\mathbb R}/P_\theta\]
be its boundary map. Then \[f_{\mathcal{A}}(\widetilde X)\subset\left(\iota^+\left(\xi(\partial_\infty\Gamma)\right)\right)^*,\] where $\iota^+$ denotes the Pl\"ucker embedding. If $f_{\mathcal{A}}$ is furthermore flag-valued, then this means that $f_{\mathcal{A}}(z)$ is transverse to $\xi(t)$ for any $z\in\widetilde X$, and $t\in\partial_\infty\Gamma$.
\end{corollary}

\begin{proof}
Zariski density makes $\Gamma$ non-elementary, so its boundary action is minimal. Thus Theorem~\ref{thm:contact-limit-set} applies.
\end{proof}

\subsection{Hermitian transversality}\label{subsec:cominuscule-opposition}

Assume from now on that $G_{\mathbb R}$ is simple Hermitian, and let $\mathcal D=G_{\mathbb R}/K_{\mathbb R}$ be its irreducible Hermitian symmetric space. A maximal parabolic subgroup of $G_{\mathbb C}$ is called \emph{cominuscule} if its unipotent radical is abelian. Or equivalently, the simple root omitted from the root system of its Levi factor has coefficient $1$ in the highest root, see e.g.~\cite{Viviani14}. The complex structure of $\mathcal D$ determines such a parabolic $P\subset G_{\mathbb C}$. The quotient $\check{\mathcal D}=G_{\mathbb C}/P$ is the compact dual of $\mathcal D$, and the Borel embedding $gK_{\mathbb R}\mapsto\eta(g)P$ identifies $\mathcal D$ biholomorphically with an open $G_{\mathbb R}$-orbit in $\check{\mathcal D}$. Denote by
\[
\check{\mathcal D}^{\mathrm{opp}}:=G_{\mathbb C}/P^{\mathrm{opp}}
\]
the opposite flag variety, choosing $P^{\mathrm{opp}}$ to be the conjugate of $P$ with respect to the real form $\mathfrak{g}_{\mathbb R}$.
Define
\[
\mathcal D^{\mathrm{opp}}
:=\eta(G_{\mathbb R})\cdot P^{\mathrm{opp}}
\subset\check{\mathcal D}^{\mathrm{opp}}.
\]
This open orbit is identified with $\mathcal D$ endowed with the opposite complex structure by $gK_{\mathbb R}\mapsto\eta(g)P^{\mathrm{opp}}$. When $P^{\mathrm{opp}}$ is conjugate to $P$, we identify the two flag varieties by their parabolic stabilizers and regard $\mathcal D^{\mathrm{opp}}$ as an open orbit in $\check{\mathcal D}$.

Let $\omega$ be the fundamental weight associated with the omitted simple root of $P$. This weight is also called \emph{cominuscule}. Choose a positive integer $m$ sufficiently divisible that $m\omega$ defines an algebraic $G_{\mathbb C}$-representation
\[
V:=V_{m\omega},\qquad W:=V^\vee.
\]
We regard $V$ and $W$ as complex representations of $G_{\mathbb R}$ by pullback along $\eta$. Their projective embeddings and associated flat bundles are taken over $\mathbb C$.

The $P$-invariant line in $V$ and the $P^{\mathrm{opp}}$-invariant line in $W$ give projective embeddings
\[
\check{\mathcal D}\hookrightarrow\mathbb P(V),
\qquad
\check{\mathcal D}^{\mathrm{opp}}\hookrightarrow\mathbb P(W).
\]
A holomorphic line subbundle $\mathcal{A}\subset\mathcal W^\vee$ gives the projective map
\[
f_{\mathcal{A}}\colon\widetilde X\longrightarrow\mathbb P(W^\vee)=\mathbb P(V).
\]
For the representation $W=V_{m\omega}^{\vee}$, the map $f_{\mathcal{A}}$ is \emph{flag-valued} precisely when
$f_{\mathcal{A}}(\widetilde X)\subset\check{\mathcal D}$.

Let $\mathcal S^{\mathrm{opp}}\subset\check{\mathcal D}^{\mathrm{opp}}$ be the unique closed $G_{\mathbb R}$-orbit. It is equivariantly isomorphic to the Shilov boundary, hence to $G_{\mathbb R}/P_{\mathcal S}$ for the corresponding real parabolic $P_{\mathcal S}$. The complex representation $W$ embeds this compact homogeneous orbit in $\mathbb P(W)$.

Set
\begin{equation}\label{eq:omega-opp-def}
\Omega_{\mathrm{opp}}:=
\{x\in \check{\mathcal D}:x\text{ is transverse to every }s\in \mathcal S^{\mathrm{opp}}\}.
\end{equation}

\begin{theorem}\label{thm:universal-hermitian-avoidance}
Let $\mathcal D\subset\check{\mathcal D}$ be the Borel-embedded irreducible Hermitian symmetric domain and let $\Omega_{\mathrm{opp}}$ be defined by \eqref{eq:omega-opp-def}. Then
\[\Omega_{\mathrm{opp}}=
\begin{cases}
\mathcal D,&\mathcal D\text{ is not of tube type},\\[1mm]
\mathcal D\sqcup\mathcal D^{\mathrm{opp}},&\mathcal D\text{ is of tube type},
\end{cases}.\]
\end{theorem}

\begin{proof}
Fix $s\in\mathcal S^{\mathrm{opp}}$. The set of points of $\check{\mathcal D}$ transverse to $s$ is an opposite Bruhat cell, hence an affine space. Since $\Omega_{\mathrm{opp}}$ is contained in this cell, it contains no positive-dimensional compact connected complex submanifold: the affine coordinate functions would all be constant on such a submanifold. Moreover, $\Omega_{\mathrm{opp}}$ is $G_{\mathbb R}$-invariant and open, because non-transversality is closed and $\mathcal S^{\mathrm{opp}}$ is compact.

Let $\mathcal O$ be an open $G_{\mathbb R}$-orbit meeting $\Omega_{\mathrm{opp}}$. Invariance gives $\mathcal O\subset\Omega_{\mathrm{opp}}$. Let $K_{\mathbb C}$ be the complex Zariski closure of $\eta(K_{\mathbb R})$ in $G_{\mathbb C}$.
By \cite[Lemma~1.3]{HuckleberryWolf02}, every open real-group orbit in a complex flag manifold contains a compact complex $K_{\mathbb C}$-orbit, called its base cycle. Here that cycle must be a point. Thus $\mathcal O$ contains a $K_{\mathbb C}$-fixed point.

We identify these fixed points. Write the Harish--Chandra decomposition
\[
\mathfrak{g}_{\mathbb C}
=\mathfrak{p}^-\oplus\mathfrak{k}_{\mathbb C}\oplus\mathfrak{p}^+,
\]
choosing signs so that the stabilizer in $G_{\mathbb C}$ of the base point of $\mathcal D$ has Lie algebra $\mathfrak{k}_{\mathbb C}\oplus\mathfrak{p}^+$ (e.g.~\cite[Theorem~2.22]{Viviani14}).
The Lie algebra $\mathfrak q_x$ of the stabilizer of a $K_{\mathbb C}$-fixed point $x\in\check{\mathcal D}$ contains $\mathfrak{k}_{\mathbb C}$.
The modules $\mathfrak{p}^\pm$ are irreducible, with different central characters. Since $\mathfrak q_x$ has the dimension of the original parabolic, it follows that
\[\mathfrak q_x=\mathfrak{k}_{\mathbb C}\oplus\mathfrak{p}^+\quad\text{or}\quad\mathfrak{k}_{\mathbb C}\oplus\mathfrak{p}^-.\]
Parabolic subgroups are connected and self-normalizing, so each of these stabilizer Lie algebras determines at most one point. The first gives the base point of $\mathcal D$. For the second possibility, the candidate stabilizer is the parabolic opposite to $P$. Since the stabilizers of points of $\check{\mathcal D}=G_{\mathbb C}/P$ are exactly the conjugates $gPg^{-1}$ of $P$, this second possibility occurs precisely when $P$ is conjugate to its opposite.

This conjugacy condition is equivalent to the opposition involution fixing the Dynkin-diagram vertex corresponding to the simple root omitted by $P$, or equivalently its fundamental weight $\omega$. Since opposition sends a highest weight to the highest weight of the dual module, this means that the irreducible $\mathfrak{g}_{\mathbb C}$-module of highest weight $\omega$ is self-dual. By \cite[Proposition~2.35]{ChaputMaubon19}, this holds precisely when $\mathcal D$ is of tube type. In that case the second point lies in $\mathcal D^{\mathrm{opp}}$. These two open orbits are distinct, since each open orbit has a unique compact $K_{\mathbb C}$-orbit \cite[\S1]{HuckleberryWolf02}. Thus the only open orbits that can meet $\Omega_{\mathrm{opp}}$ are $\mathcal D$ and, in tube type, $\mathcal D^{\mathrm{opp}}$.

We next prove their inclusion and exclude their boundaries. Conjugation of parabolic stabilizers with respect to $\mathfrak{g}_{\mathbb R}$ gives a $G_{\mathbb R}$-equivariant antiholomorphic map $\sigma\colon\check{\mathcal D}\to\check{\mathcal D}^{\mathrm{opp}}$. It sends the Shilov boundary $\mathcal S$ of $\mathcal D$ onto $\mathcal S^{\mathrm{opp}}$.

Use the Harish--Chandra chart $j\colon\mathfrak{p}^-\to\check{\mathcal D}$, $j(z)=\exp(z)P$. Its image consists exactly of the flags transverse to the base point $P^{\mathrm{opp}}\in\check{\mathcal D}^{\mathrm{opp}}$ and contains $\overline{\mathcal D}$ in the bounded realization. Write $\overline z$ for conjugation with respect to $\mathfrak{g}_{\mathbb R}$. For $x,y\in\mathfrak{p}^-$, the pair $(x,y)$ is \emph{quasi-invertible} in the sense of \cite{ClercNeeb06} if and only if $\exp(-\overline y)\cdot j(x)\in j(\mathfrak{p}^-)$. Since $\sigma(j(y))=\exp(\overline y)P^{\mathrm{opp}}$, invariance of transversality shows that this is equivalent to $j(x)$ being transverse to $\sigma(j(y))$.

By \cite[Theorem~2.6]{ClercNeeb06}, a pair $(x,y)\in\overline{\mathcal D}\times\overline{\mathcal D}$ is quasi-invertible exactly when no proper face of the closed bounded domain contains both points. An interior point lies in no proper face. Thus, for every $x\in\mathcal D$ and $s\in\mathcal S$, the pair $(x,s)$ is quasi-invertible, so $x$ is transverse to $\sigma(s)$. Since $\sigma(\mathcal S)=\mathcal S^{\mathrm{opp}}$, each $x\in\mathcal D$ is transverse to every point of $\mathcal S^{\mathrm{opp}}$, hence belongs to $\Omega_{\mathrm{opp}}$ by its definition. Therefore $\mathcal D\subset\Omega_{\mathrm{opp}}$.

Conversely, every $x\in\partial\mathcal D$ lies in a proper closed face, which contains an extreme point $s$. The extreme points are the Shilov points \cite[Proposition~2.1]{ClercNeeb06}. Hence $x$ is not transverse to $\sigma(s)\in\mathcal S^{\mathrm{opp}}$, and $\partial\mathcal D\cap\Omega_{\mathrm{opp}}=\varnothing$.

In tube type the two flag types coincide, so $\sigma$ is an involution of the compact dual exchanging $\mathcal D$ and $\mathcal D^{\mathrm{opp}}$.
It also acts on the opposite flag variety, preserves transversality, and preserves its unique closed real orbit $\mathcal S^{\mathrm{opp}}$.
It therefore preserves $\Omega_{\mathrm{opp}}$, giving
\[
\mathcal D^{\mathrm{opp}}\subset\Omega_{\mathrm{opp}},\qquad
\partial\mathcal D^{\mathrm{opp}}\cap\Omega_{\mathrm{opp}}=\varnothing.
\]

Finally, there are finitely many $G_{\mathbb R}$-orbits in the compact dual, and their open orbits form a dense subset (see \cite[Theorem~3.8]{WolfZierau97}). Since $\Omega_{\mathrm{opp}}$ is open, each of its points is a limit of points in the permitted open orbits. The boundary exclusions show that it belongs to one of them, proving the asserted equality.
\end{proof}

\subsection{Nonlattice-transversality alternative for Hermitian groups}

When $\Gamma$ is a lattice, the equality in the theorem below implies that every $f_{\mathcal{A}}(z)$ is transverse to every point of $\mathcal S^{\mathrm{opp}}$. Thus $f_{\mathcal{A}}(\widetilde X)\subset\Omega_{\mathrm{opp}}$.
The classification of $\Omega_{\mathrm{opp}}$ and connectedness of $\widetilde X$ then place the image in one Hermitian component.

\begin{theorem}\label{thm:main-B}
Let $G_{\mathbb R}$ be a simple Hermitian real group with finite center and without compact factors, and let $W=V_{m\omega}^{\vee}$ be the dual of a positive algebraic cominuscule Cartan power of $G_{\mathbb C}$. If there is a filtered holomorphic line subbundle $\mathcal{A}_\bullet\subset\mathcal W_\bullet^\vee$ whose projective map is flag-valued, i.e.
\[f_{\mathcal{A}}\colon\widetilde X\longrightarrow \check{\mathcal D},\]
and
\[\lambda_1(W)=\frac{2\operatorname{pardeg}\mathcal{A}_\bullet}{\deg\KK_{\overline{X}}(D)},\]
then either $\Gamma$ is not a lattice in $G_{\mathbb R}$, or the flag map takes values in a single Hermitian component: $\mathcal D$ in non-tube type, and either $\mathcal D$ or $\mathcal D^{\mathrm{opp}}$ in tube type.
\end{theorem}

\begin{proof}
Assume that $\Gamma$ is a lattice. The representation $W$ is irreducible, and $\mathcal S^{\mathrm{opp}}$ is a compact $G_{\mathbb R}$-orbit.
Corollary~\ref{cor:lattice-transversality} and the flag-valued hypothesis therefore give
\[
f_{\mathcal{A}}(\widetilde X)\subset\check{\mathcal D}\cap(\mathcal S^{\mathrm{opp}})^*=\Omega_{\mathrm{opp}}.
\]
Theorem~\ref{thm:universal-hermitian-avoidance} identifies $\Omega_{\mathrm{opp}}$ with $\mathcal D$ in non-tube type and $\mathcal D\sqcup\mathcal D^{\mathrm{opp}}$ in tube type. Since $\widetilde X$ is connected and $f_{\mathcal{A}}$ is continuous, its image is contained in one of these Hermitian components.
\end{proof}

\subsection{Classical Hermitian models}
\label{subsec:classical-hermitian-models}

In this subsection, we consider Theorem~\ref{thm:universal-hermitian-avoidance} in the following classical models. A holomorphic equivariant map into one of these components pulls back the homogeneous Hodge filtration to the same flat local system, giving the induced Hermitian variation. See \cite{Helgason01,GriffithsSchmid69} for references.

\begin{example}[Type A\upRoman{3}]\label{ex:hermitian-aiii}
Let $G_{\mathbb R}=\mathrm{SU}_{p,q}$, $1\leqslant p\leqslant q$, and let $h$ be the preserved Hermitian form of signature $(p,q)$ on the standard complex representation $V$. The compact dual is $\mathrm{Gr}_p(V)$, and the real flag orbit is the orbit of maximal
$h$-isotropic $p$-planes. The Hermitian component is
\[
\mathcal D^+=\{E\in\mathrm{Gr}_p(V):h|_E>0\}.
\]
If $p=q$, there is also the opposite component
\[
\mathcal D^-=\{E\in\mathrm{Gr}_p(V):h|_E<0\}.
\]
For $p<q$ only $\mathcal D^+$ occurs in the transversality domain.

If $f\colon\widetilde X\to \mathcal D^\varepsilon$ is the equivariant holomorphic map and $\mathcal{E}\subset\mathcal{V}$ is the corresponding rank-$p$ holomorphic subbundle, choose $\varepsilon\in\{\pm1\}$ so that $\varepsilon h|_{\mathcal{E}}>0$. Then, as smooth complex vector bundles,
\[
\mathcal{V}=\mathcal{E}\oplus\mathcal{E}^{\perp_h},\qquad
\mathcal{V}^{1,0}_{\mathrm{ind}}=\mathcal{E},\qquad
\mathcal{V}^{0,1}_{\mathrm{ind}}=\mathcal{E}^{\perp_h},
\]
and, with polarization $\varepsilon h$, this is the induced polarized complex variation of Hodge structure of weight one and Hodge numbers $(p,q)$.
\end{example}

\begin{example}[Type C\upRoman{1}]
\label{ex:hermitian-ci}
Let $G_{\mathbb R}=\mathrm{Sp}_{2n}\mathbb R$, let $Q$ be the preserved real symplectic form, and put
\[
H(v,w)=\sqrt{-1}\,Q(v,\overline w)
\qquad (v,w\in V_{\mathbb C}).
\]
The compact dual is the Lagrangian Grassmannian $\mathrm{LGr}_n(V_{\mathbb C})$. Its two Hermitian components are
\[
\mathcal D^+=\{E:H|_E>0\},\qquad
\mathcal D^-=\{E:H|_E<0\}.
\]
We call $\mathcal D^+$ and $\mathcal D^-$ the \emph{Siegel components}. They are the two open $\mathrm{Sp}_{2n}\mathbb R$-orbits of complex Lagrangian subspaces on which $H$ is definite.

If $f\colon\widetilde X\to \mathcal D^\varepsilon$ is equivariant and holomorphic and $\mathcal{E}\subset\mathcal{V}$ is the corresponding holomorphic Lagrangian subbundle, choose $\varepsilon\in\{\pm1\}$ so that $\varepsilon H|_{\mathcal{E}}>0$. Then, as smooth complex vector bundles,
\[
\mathcal{V}=\mathcal{E}\oplus\overline{\mathcal{E}},\qquad
\mathcal{V}^{1,0}_{\mathrm{ind}}=\mathcal{E},\qquad
\mathcal{V}^{0,1}_{\mathrm{ind}}=\overline{\mathcal{E}},
\]
so $F^1_{\mathrm{ind}}=\mathcal{E}$ and, with polarization $\varepsilon Q$, one obtains the induced polarized real variation of Hodge structure of weight one and Hodge numbers $(n,n)$ on the same flat real local system.
\end{example}

\begin{example}[Type D\upRoman{3}]
\label{ex:hermitian-diii}
Let $G_{\mathbb R}=\mathrm{SO}_{2n}^{*}$. On the standard complex representation $V$, let $Q$ be the preserved nondegenerate complex symmetric form and let $j$ be the preserved anti-linear quaternionic structure, so that $j^2=-1$ and $Q(jv,jw)=-\overline{Q(v,w)}$. Put
\[h(v,w)=Q(v,jw).\]

Let $\mathrm{OGr}_n(V,Q)$ denote the variety of $n$-dimensional complex subspaces $E\subset V$ with $Q|_{E\times E}=0$. These are maximal $Q$-isotropic subspaces. It has two connected components, denoted $\mathrm{OGr}_n^\delta(V,Q)$ for $\delta\in\{+,-\}$. Choose $\delta$ so that this component contains an $h$-positive plane, and take it as the compact dual. Its positive Hermitian domain is
\[\mathcal D^+=\{E:h|_E>0\}.\]
When $n$ is even, the same compact-dual component also contains
\[\mathcal D^-=\{E:h|_E<0\}.\]
Type D\upRoman{3} is of tube type exactly when $n$ is even. Hence Theorem~\ref{thm:universal-hermitian-avoidance} gives $\Omega_{\mathrm{opp}}=\mathcal D^+$ for odd $n$ and $\Omega_{\mathrm{opp}}=\mathcal D^+\sqcup\mathcal D^-$ for even $n$.
For odd $n$, the $h$-negative planes lie in the other component of the maximal isotropic Grassmannian.

If $f\colon\widetilde X\to \mathcal D^\varepsilon$ is equivariant and holomorphic and $\mathcal{E}$ is the corresponding maximal $Q$-isotropic subbundle, choose $\varepsilon\in\{\pm1\}$ so that $\varepsilon h|_{\mathcal{E}}>0$. Then, as smooth complex vector bundles,
\[
\mathcal{V}=\mathcal{E}\oplus j\mathcal{E},\qquad
\mathcal{V}^{1,0}_{\mathrm{ind}}=\mathcal{E},\qquad
\mathcal{V}^{0,1}_{\mathrm{ind}}=j\mathcal{E},
\]
and the flat Hermitian form $\varepsilon h$ polarizes the induced complex variation of Hodge structure of weight one.

\end{example}

\begin{example}[Type BD\upRoman{1}]
\label{ex:hermitian-iv}
Let $G_{\mathbb R}=\mathrm{SO}_{2,q}^{0}$, $q\geqslant3$, and let $Q$ have signature $(2,q)$ on the standard real representation $V_{\mathbb R}$. The compact dual is the isotropic quadric
\[\mathcal Q=\{[u]\in\mathbb P(V_{\mathbb C}):Q(u,u)=0\}.\]
The transversality domain is
\[\mathcal D^+\sqcup \mathcal D^-=\{[u]\in\mathcal Q:Q(u,\overline u)>0\},\]
whose two components are distinguished by the orientation of the positive real two-plane spanned by $\operatorname{Re}u$ and $\operatorname{Im}u$.

If $f\colon\widetilde X\to \mathcal D^\varepsilon$ is equivariant and holomorphic and $\LL\subset\mathcal{V}$ is the resulting isotropic line bundle, then the induced filtration is
\[
F^2_{\mathrm{ind}}=\LL,\qquad
F^1_{\mathrm{ind}}=\LL^{\perp_Q},\qquad
F^0_{\mathrm{ind}}=\mathcal{V}.
\]
After the standard overall sign normalization of $Q$, this is the induced polarized real variation of Hodge structure of weight two and Hodge numbers $(1,q,1)$. Its only nontrivial Griffiths-transversality condition follows from $Q(s,s)=0$ for a local generator $s$ of $\LL$, since
\[Q(\nabla s,s)=0
\quad\Longrightarrow\quad
\nabla F^2_{\mathrm{ind}}\subset
F^1_{\mathrm{ind}}\otimes\Omega_X^1.\]
\end{example}

\section{Variations of Hodge structure and hypergeometric applications}\label{sec:CY-applications}

In this section, we discuss more about applying Theorem~\ref{thm:main-B} to symplectic and special unitary flat bundles, and then to their Hodge filtrations. 

We first recall uniqueness of the Hodge grading on an irreducible complex local system. An \emph{integral translation of the Hodge grading} is a reindexing by one integer $r\in\mathbb Z$: the translated Hodge summands and filtration are
\[
\widetilde{\mathcal{E}}^{p}=\mathcal{E}^{p+r},\qquad
\widetilde F^p=F^{p+r}.
\]
For a real or rational variation, both Hodge indices are translated by the same integer, so $\widetilde V^{p,q}=V^{p+r,q+r}$. The weight changes from $w$ to $w-2r$, and with our convention the polarization changes from $Q$ to $(-1)^rQ$. The underlying flat local system is unchanged.

\begin{lemma}\label{lem:grading-uniqueness}
Let $\mathbb{V}$ be an irreducible complex local system on $X$ with non-expanding cusp monodromy. Any two polarized complex variations of Hodge structure on $\mathbb{V}$ differ by a unique integral translation of the Hodge grading. Writing $\mathcal{E}_i^p$ for their Hodge summands and $F_i^p$ for their Hodge filtrations, there is a unique $m\in\mathbb Z$ such that
\[
\mathcal{E}_1^p=\mathcal{E}_2^{p+m},\qquad
F_1^p=F_2^{p+m}\qquad(p\in\mathbb Z).
\]
Equivalently, if $g_{i,t}$ denotes the $\mathbb C^*$-grading action of the $i$-th variation on the common harmonic bundle, with $g_{i,t}|_{\mathcal{E}_i^p}=t^p\operatorname{id}$, then $g_{2,t}=t^m g_{1,t}$.
\end{lemma}

\begin{proof}

	The two variations give tame harmonic metrics adapted to the same canonical parabolic structure. Since $\mathbb{V}$ is irreducible, uniqueness of tame harmonic metric (see \cite{Simpson90,Mochizuki07}) identifies these metrics up to a positive scalar and hence identifies the corresponding parabolic Higgs bundle $(\mathcal{E},\bar\partial_{\mathcal{E}},\Phi)$. This Higgs bundle is stable and therefore simple. The grading actions preserve the parabolic filtration and satisfy $g_t\Phi g_t^{-1}=t^{-1}\Phi$. If $g_t$ and $\widetilde g_t$ are the two grading actions, then $a_t=\widetilde g_tg_t^{-1}$ satisfies $[a_t,\Phi]=0$. Hence $a_t=\chi(t)\operatorname{id}$. Since both grading actions are algebraic, $\chi$ is an algebraic character of $\mathbb C^*$, so $\chi(t)=t^m$ for a unique $m\in\mathbb Z$. Comparing weight spaces gives the asserted translation.
\end{proof}

\begin{corollary}\label{cor:real-grading-uniqueness}
	Assume that the two variations in Lemma~\ref{lem:grading-uniqueness} are defined on the same real, respectively rational, local system whose complexification is irreducible, and have weights $w_1,w_2$. Then
	\[w_2-w_1=2m,\qquad V_1^{p,q}=V_2^{p+m,q+m}.\]
	Thus the two Hodge decompositions differ by a simultaneous integral translation of both Hodge indices.
\end{corollary}

\begin{proof}
	Let $\sigma$ be the common real structure. From $\mathcal{E}_1^p=\mathcal{E}_2^{p+m}$ and $\sigma(\mathcal{E}_i^p)=\mathcal{E}_i^{w_i-p}$ we obtain, for every nonzero Hodge summand,
	\[\mathcal{E}_2^{w_2-p-m}=\mathcal{E}_1^{w_1-p}=\mathcal{E}_2^{w_1-p+m}.\]
	Since the Hodge decomposition is direct, $w_2-w_1=2m$, and then $V_1^{p,q}=V_2^{p+m,q+m}$. The rational case follows after extension of scalars to $\mathbb R$.
\end{proof}

\subsection{Nonlattice-transversality alternative for VHS}

For a real or complex vector space $U$ equipped with a bilinear or Hermitian form $Q$, write $\mathrm{IGr}_k(U,Q)$ for the Grassmannian of $Q$-isotropic $k$-planes. Let $(V_{\mathbb{R}},Q)$ be symplectic of dimension $2n$. In this subsection, for $1\leqslant k\leqslant n$, let $W_k=\bigwedge\nolimits_0^kV_{\mathbb R}$ be the primitive fundamental representation, which is the kernel of contraction with respect to $Q$. We write $W_{k,\mathbb C}=W_k\otimes_{\mathbb R}\mathbb C$.
The invariant pairing
\[
B_k(v_1\wedge\cdots\wedge v_k,w_1\wedge\cdots\wedge w_k)
=\det\left(Q(v_i,w_j)\right)
\]
restricts nondegenerately to $W_k$. We use it to identify $W_k$ with $W_k^\vee$, and likewise their complexifications.

\begin{theorem}\label{thm:symplectic-contact}
	Let $\mathbb{V}$ be a semisimple symplectic flat bundle of rank $2n$, with non-expanding cusp monodromy $\Gamma<\mathrm{Sp}_{2n}\mathbb R$, associated de Rham holomorphic bundle $\mathcal{V}=(V_{\mathbb C},\nabla^{0,1})$. Let $\mathcal{E}_{k,\bullet}\subset\mathcal V_\bullet$ be a filtered holomorphic subbundle of rank $1\leqslant k\leqslant n$ whose restriction $\mathcal E_k$ to $X$ is isotropic. Suppose that
	\begin{equation}\label{eq:LyaHNsymplectic}
		\lambda_1(\mathbb{V})+\cdots+\lambda_k(\mathbb{V})=\frac{2\operatorname{pardeg}\mathcal{E}_{k,\bullet}}{\deg\KK_{\overline{X}}(D)}.
	\end{equation}
	Then the following hold:
	\begin{itemize}
		\item If $k<n$, then $\Gamma$ is not a lattice. In particular, this applies when $\mathcal{E}_{k,\bullet}$ is an isotropic Harder--Narasimhan step.
		
		\item If $k=n$, then either $\Gamma$ is not a lattice, or the Lagrangian flag map takes values in a Siegel component. The latter condition is equivalent to $\mathcal{E}_n$ being the $F^1$-bundle of the polarized real weight-one variation of type C\upRoman{1} described in Example~\ref{ex:hermitian-ci}, on the same flat real local system.
	\end{itemize}
	In particular, if $\mathbb{V}$ underlies a polarized real variation of Hodge structure of odd weight $2m+1$, then $F^{m+1}$ is Lagrangian and extends as a filtered holomorphic subbundle of $\mathcal V_\bullet$. If \eqref{eq:LyaHNsymplectic} holds for $\mathcal{E}_n=F^{m+1}$ and $\Gamma$ is a lattice, its Hodge grading is an integral translation of that of a polarized real weight-one variation of type C\upRoman{1}. Thus $\Gamma$ is not a lattice if no such translation exists, in particular when the Hodge decomposition has more than two nonzero summands.
\end{theorem}

\begin{proof}
	Since $\lambda_1(W_k)=\lambda_1(\mathbb{V})+\cdots+\lambda_k(\mathbb{V})$ and $\operatorname{pardeg}\det\mathcal{E}_{k,\bullet}=\operatorname{pardeg}\mathcal{E}_{k,\bullet}$, \eqref{eq:LyaHNsymplectic} is the Lyapunov--HN equality for $\mathcal{A}=B_k(\cdot,\det\mathcal{E}_k)\subset\mathcal W_k^\vee$. 
	
	Suppose first that $k<n$. If $\Gamma$ were a lattice, Corollary~\ref{cor:lattice-transversality} would force every complex isotropic $k$-plane $U$ in the map $f_{\mathcal{A}}$ to be transverse to every real isotropic $k$-plane. But
	\[R_U=\{v\in V_{\mathbb R}:Q(v,U)=0\}\]
	has real dimension at least $2n-2k>0$. Choosing $0\neq p_1\in R_U$, extending $\mathbb Rp_1$ to a real Lagrangian, and choosing a real isotropic $k$-plane $P=\operatorname{span}_{\mathbb R}\{p_1,\ldots,p_k\}$ inside it, we get
	\[B_k(\det P,\det U)=\det\left(Q(p_i,u_j)\right)=0\]
	for any basis $u_1,\ldots,u_k$ of $U$, because the first row vanishes, contradiction.
	
	Now let $k=n$. The Pl\"ucker line belongs to the type-C\upRoman{1} representation $V_{\omega_n}$. If $\Gamma$ were a lattice, Theorem~\ref{thm:main-B} places the Lagrangian flag map in a Siegel component. By the type-C\upRoman{1} description in Example~\ref{ex:hermitian-ci}, this is equivalent to $\mathcal{E}_n$ being the $F^1$-bundle of the corresponding polarized real weight-one variation on the same flat real local system.
	
	Finally, for a polarized odd-weight variation, the results follows directly from above discussions and Corollary~\ref{cor:real-grading-uniqueness}.
\end{proof}

We have a similar result for $\mathrm{SU}_{m,n}$-bundles.

\begin{theorem}\label{thm:unitary-contact}
	Let $\mathbb{V}$ be a semisimple flat $\mathrm{SU}_{m,n}$-bundle, where $m\geqslant n\geqslant1$, with non-expanding cusp monodromy $\Gamma<\mathrm{SU}_{m,n}$, flat Hermitian form $H$ of signature $(m,n)$, associated de Rham holomorphic bundle $\mathcal{V}=(V,\nabla^{0,1})$.
	Let $\mathcal{E}_{k,\bullet}\subset\mathcal V_\bullet$ be a filtered holomorphic subbundle of rank $1\leqslant k<m+n$, and write $\mathcal E_k$ for its restriction to $X$. Suppose that
	\begin{equation}\label{eq:LyaHNunitary}
		\lambda_1(\mathbb{V})+\cdots+\lambda_k(\mathbb{V})=\frac{2\operatorname{pardeg}\mathcal{E}_{k,\bullet}}{\deg\KK_{\overline{X}}(D)}.
	\end{equation}
	Then the following hold:
	\begin{itemize}
		\item If $k\notin\{m,n\}$, then $\Gamma$ is not a lattice.
		
		\item If $k\in\{m,n\}$, then either $\Gamma$ is not a lattice, or $H|_{\mathcal{E}_k}$ is definite everywhere. In the latter case, if $m>n$, its sign is negative for $k=n$ and positive for $k=m$. Moreover, $\mathcal{E}_k$ is the $F^1$-bundle of a polarized complex weight-one variation of type A\upRoman{3} described in Example~\ref{ex:hermitian-aiii}, on the same flat complex local system, with the appropriate sign of $H$.
	\end{itemize}
	In particular, if $\mathbb{V}$ underlies a polarized complex variation of Hodge structure and \eqref{eq:LyaHNunitary} holds for $\mathcal{E}_k=F^p$, then, if $\Gamma$ is a lattice, its Hodge grading is an integral translation of that of a polarized complex weight-one variation of type A\upRoman{3}. Thus $\Gamma$ is not a lattice if no such translation exists, in particular when the Hodge decomposition has more than two nonzero summands.
\end{theorem}

\begin{proof}
	Set \[W:=\left(\bigwedge^kV\right)^\vee,\qquad \mathcal{A}:=\det\mathcal{E}_k\subset\mathcal W^\vee.\]
	Then $W$ is an irreducible representation for $\mathrm{SL}(V)$ and $\lambda_1(W)=\sum_{i=1}^k\lambda_i(\mathbb{V})$. Thus \eqref{eq:LyaHNunitary} is equivalent to the Lyapunov--HN equality for $\mathcal{A}$.
	
	Suppose that $k\leqslant m$ and $\Gamma$ is a lattice. Let $P$ be an $H$-isotropic $k$-planes if $k\leqslant n$ or nonnegative $k$-planes with $n$-dimensional null space if $n<k\leqslant m$. For a basis $p_1,\ldots,p_k$ of $P$, define $\ell_P\in W$ by
	\[\ell_P(v_1\wedge\cdots\wedge v_k)=\det\left(H(v_i,p_j)\right).\]
	This construction is $\mathrm{SU}_{m,n}$-equivariant. By Witt's extension theorem, see e.g. \cite[Theorem~7.4]{Taylor92}, any two such planes are related by an $H$-isometry. Without loss of generality, we can assume that this isometry has determinant $1$. Thus these lines form an $\mathrm{SU}_{m,n}$-orbit $\mathcal S\subset\mathbb P(W)$, which is compact. By Corollary~\ref{cor:lattice-transversality}, we obtain that $f_{\mathcal{A}}(z)\in\mathcal S^*$.
	
	Let $E:=\EE_k(z)$ for some $z\in\widetilde{X}$. Choose an orthogonal decomposition $V=V_+\oplus V_-$ with respect to $H$, with $H$ positive on $V_+$ and negative on $V_-$. Then $\dim E^{\perp_H}=m+n-k$. By $\dim(A\cap B)\geqslant\dim A+\dim B-\dim V$, we obtain that
	\[\dim(E^{\perp_H}\cap V_+)\geqslant m-k,\qquad\dim(E^{\perp_H}\cap V_-)\geqslant n-k.\]
	If $k<n$, both intersections are nonzero, so $E^{\perp_H}$ contains a nonzero null vector $p$. Extend $\mathbb Cp$ to an $H$-isotropic $k$-plane $P$. If $n<k<m$, the inequality above gives a positive vector $p\in E^{\perp_H}$. Extend $\mathbb Cp$ to a nonnegative $k$-plane $P$ with $n$-dimensional null space. In either case, 
	\[\det\left(H(u_i,p_j)\right)=0,\]
	for bases $u_1,\ldots,u_k$ of $E$ and $p_1=p,p_2,\ldots,p_k$ of $P$, because $H(u_i,p)\equiv 0$, contradiction. For $k=n$ or $k=m$, the statements follow directly from Theorem~\ref{thm:main-B} and Example~\ref{ex:hermitian-aiii}.
	
	For $k>m$, the same argument applied to the annihilator of $\mathcal{E}_k$ in $\mathcal{V}^\vee$, of rank $m+n-k<n$, with induced filtered extension $(\mathcal V_\bullet/\mathcal E_{k,\bullet})^\vee$, shows that $\Gamma$ is not a lattice. Finally, Zariski density implied by Borel density makes $\mathbb{V}$ irreducible, and Lemma~\ref{lem:grading-uniqueness} gives the assertion for Hodge filtrations.
\end{proof}

\subsection{Weight-three variations and Assumption A}\label{subsec:weight-three}

Let $\mathbb{V}$ be a polarized complex weight-three variation with Hodge numbers $(m-1,1,1,n-1)$, where $m,n\geqslant2$, over the complete finite-area hyperbolic surface $X=\overline{X}\setminus D$. Write $h$ for its Hodge metric and $H$ for its flat Hermitian polarization, normalized so that $H|_{V^{p,q}}=(-1)^qh$. Thus $H$ has signature $(m,n)$, and $F^2=V^{3,0}\oplus V^{2,1}$ has rank $m$ and signature $(m-1,1)$. We assume that the monodromy lies in $\mathrm{SU}_{m,n}$. For each Hodge bundle, write $V^{p,q}_\bullet=\Xi_h(V^{p,q})$ for its metric extension with respect to the Hodge holomorphic structure $D_h^{0,1}$. We use the same symbol $\Phi_{p,q}$ for its logarithmic extension, with values in $\KK_{\overline X}(D)$, which carries parabolic weight zero. 

Below, without loss of generality, we assume that $m\geqslant n$. Otherwise we consider the dual bundle. Put
\[W_n:=\bigwedge^n V,\qquad\mathrm{IGr}_n(V,H)=\{L\in\mathrm{Gr}_n(V):H|_L=0\}.\]
The compact flag manifold $\mathrm{IGr}_n(V,H)$ parametrizes the maximal $H$-isotropic complex subspaces of $V$. The two lowest Hodge pieces of $W_n$ are
\[W_n^{1,3n-1}=\det V^{0,3}\otimes V^{1,2},\qquad W_n^{2,3n-2}=\det V^{0,3}\otimes V^{2,1},\]
both lines, and the component of the Higgs field between them is
\begin{equation}\label{eq:signature-one-induced-current}
	\left.\Phi_{W_n}\right|_{W_n^{2,3n-2}}=\mathrm{id}_{\det V^{0,3}}\otimes\Phi_{2,1}:\ W_n^{2,3n-2}\longrightarrow W_n^{1,3n-1}\otimes\KK_X.
\end{equation}

The following theorem generalizes the definition of Assumption~A and corresponding results in \cite{Filip21}. 

\begin{theorem}\label{thm:main-C}
Let $\mathbb{V}$ be as above. Suppose that the logarithmic extension
\[\Phi_{2,1}\colon V^{2,1}_\bullet\longrightarrow V^{1,2}_\bullet\otimes\KK_{\overline X}(D)\]
is an isomorphism of filtered line bundles on $\overline X$. Then there is a proper closed set $\Lambda\subset\mathrm{IGr}_n(V,H)$, which is monodromy-invariant, such that $F^2(z)\cap L=0$ for every $L\in\Lambda$ and every $z\in\widetilde X$. Moreover,
\begin{equation}\label{eq:signature-one-contact-current} 
	\lambda_1(\mathbb{V})+\cdots+\lambda_m(\mathbb{V})=\frac{2\operatorname{pardeg}F^2}{\deg\KK_{\overline{X}}(D)}.
\end{equation}

Furthermore, the monodromy $\Gamma$ is not a lattice in $\mathrm{SU}_{m,n}$. If $m=n$ and the variation is the complexification of a real polarized variation, then $\Gamma<\mathrm{Sp}_{2n}\mathbb R$ is not a lattice in $\mathrm{Sp}_{2n}\mathbb R$ either.
In either case, if $\Gamma$ is contained in an arithmetic lattice $\Gamma_{\mathrm{ar}}$ of the ambient group, then $[\Gamma_{\mathrm{ar}}:\Gamma]=\infty$. If $\Gamma$ is also Zariski dense in that group, it is thin.
\end{theorem}

\begin{proof} 
	We first give the uniform lower bound of $\Phi_{2,1}$, which follows the idea of \cite[Proposition~2.1.11]{Filip21} and \cite[Lemma 3.13]{Zhang25}. At a cusp with coordinate $z$ centered at the cusp, we set $y=-\log|z|$ and choose holomorphic frames $e_{2,1}$ and $e_{1,2}$ of the extension of $V^{2,1}$ and $V^{1,2}$ respectively such that
	\[\Phi_{2,1}(e_{2,1})=e_{1,2}\,\mathrm{d}z/z.\]
	Since $\Phi_{2,1}$ is still an isomorphism after extension, this gives equal parabolic weights on $V^{2,1}$ and $V^{1,2}$, and the residue of $\Phi_{2,1}$ is nonzero. Then the comparison between the harmonic metric and Simpson's model metric \cite[\S5]{Simpson90} gives
	\[\frac{\|e_{1,2}\|_h^2}{\|e_{2,1}\|_h^2}\asymp y^{-2},\qquad\left\|\frac{\mathrm{d}z}{z}\right\|_{\mathrm{hyp}}^2\asymp y^2.\]
	Thus $\|\Phi_{2,1}\|_{h,\mathrm{hyp}}^2\asymp1$ near each cusp. Therefore, the uniform lower bound follows by the smoothness on the complement compact part. Hence there are $0<\delta\leqslant M$ such that $\delta\leqslant\|\Phi_{2,1}\|_{h,\mathrm{hyp}}\leqslant M$ on $X$.
	
	Trivialize the flat bundle over $\widetilde X$. For $L\in\mathrm{IGr}_n(V,H)$, choose $0\neq w_L\in\bigwedge^nL$ and put $f_L=\|w_L^{1,3n-1}\|_h^2$. In the Hodge splitting, the quotient map $V\to V/F^2(z)$ is projection onto $V^{1,2}(z)\oplus V^{0,3}(z)$, whose determinant is $W_n^{1,3n-1}(z)$. Its top exterior power therefore sends $w_L$ to $w_L^{1,3n-1}$. Since both $L$ and the quotient have dimension $n$, this image vanishes exactly when the restriction to $L$ has nonzero kernel, namely $F^2(z)\cap L$. Thus 
	\[f_L(z)=0\quad\Longleftrightarrow\quad F^2(z)\cap L\neq0.\]
	Set $L_0:=V^{0,3}\oplus V^{2,1}$, $P_0:=V^{3,0}\oplus V^{1,2}$. Then $H|_{L_0}=-h$ and $H|_{P_0}=h$. Since $H|_L=0$, the projection $L\to L_0$ with respect to the decomposition $L_0\oplus P_0$ is injective. Hence it is an isomorphism. Thus $L$ is the graph of a linear map $U\colon L_0\to P_0$. Since $L$ is isotropic, we have the identity $0=-h(v,w)+h(Uv,Uw)$, which shows that $U^*U=I_n$. 

	Choose $h$-orthonormal bases $e_1,\ldots,e_n$ of $L_0$ and $f_1,\ldots,f_m$ of $P_0$, with $e_1,\ldots,e_{n-1}\in V^{0,3}$, $e_n\in V^{2,1}$, $f_1,\ldots,f_{m-1}\in V^{3,0}$, and $f_m\in V^{1,2}$. Write $Ue_j=\sum_{i=1}^m U_{ij}f_i$. The isomorphism $v\mapsto v+Uv$ from $L_0$ to $L$ sends the first basis to $e_j+Ue_j$. Since $L_0$ and $P_0$ are $h$-orthogonal and $U^*U=I_n$,
	\[h(e_j+Ue_j,e_k+Ue_k)=h(e_j,e_k)+h(Ue_j,Ue_k)=2\delta_{jk}.\]
	Write $w_L=a(e_1+Ue_1)\wedge\cdots\wedge(e_n+Ue_n)$, with $a\neq0$.
	Therefore,
	\[
		w_L^{1,3n-1}=aU_{mn}\,e_1\wedge\cdots\wedge e_{n-1}\wedge f_m,\qquad w_L^{2,3n-2}=a\,e_1\wedge\cdots\wedge e_n.\]
	Since $|U_{mn}|^2\leqslant\sum_{i=1}^m|U_{in}|^2=1$ and $\|w_L\|_h^2=|a|^2\det(2I_n)=2^n|a|^2$, we obtain that
	\[\|w_L^{1,3n-1}\|_h\leqslant\|w_L^{2,3n-2}\|_h,\qquad\|w_L\|_h=2^{n/2}\|w_L^{2,3n-2}\|_h.\]
	The flatness of $w_L$ gives
	\[D_h^{0,1}w_L^{1,3n-1}=0,\qquad D_h^{1,0}w_L^{1,3n-1}=-\Phi_{W_n}w_L^{2,3n-2}.\]
	These identities show that $w_L^{1,3n-1}$ is nonzero and holomorphic. Thus its zeros, and hence those of $f_L$, are discrete. The same estimates give $\|\nabla f_L\|\geqslant c f_L$, where $c>0$ depends only on $\delta$ and the fixed norm conventions, so $f_L$ has no critical points away from its zero set. At a zero $p$, holomorphicity and the second identity imply that $D_{h,\xi}w_L^{1,3n-1}\neq0$ for every nonzero real tangent vector $\xi\in T_p\widetilde X$. Therefore,
	\[\operatorname{Hess}_p f_L(\xi,\xi)=2\left\|D_{h,\xi}w_L^{1,3n-1}\right\|_h^2>0\]
	by metric compatibility.
	Thus $f_L$ is a Morse function. Applying \cite[Lemma~4.1 and Theorem~4.3(2)]{Zhang25} on the complete surface $\widetilde X$, we obtain that $f_L$ has at most one zero. See also \cite[Proposition~2.2.4 and Theorem~2.2.8]{Filip21}.
	
	Following the argument in \cite[Theorem~4.1]{LiZhang26}, consider
	\[\mathcal I=\left\{(L,z)\in\mathrm{IGr}_n(V,H)\times\widetilde X:f_L(z)=0\right\},\]
	equipped with the projection $\pi\colon\mathcal I\to\mathrm{IGr}_n(V,H)$ onto the first factor. Locally choose $w_L$ smoothly in $L$ and write $w_L^{1,3n-1}(z)=b(L,z)e(z)$ in a smooth $h$-unit frame. Then $\mathcal I$ is defined by $b=0$, and at each of its points
	\[\operatorname{Hess}_z f_L(\xi,\xi)=2|\mathrm{d}_z b(\xi)|^2.\]
	Therefore, $\mathrm{d}_z b\colon T_z\widetilde X\to\mathbb C$ is injective, hence an isomorphism. The implicit function theorem therefore makes $\mathcal I$ a smooth submanifold, with
	\[T_{(L,z)}\mathcal I=\left\{(\eta,\xi):\mathrm{d}_L b(\eta)+\mathrm{d}_z b(\xi)=0\right\}.\]
	The differential $\mathrm{d}\pi(\eta,\xi)=\eta$ is an isomorphism since $\mathrm{d}_zb$ is an isomorphism. Thus $\pi$ is a local diffeomorphism. Note that $\pi$ is injective because the zero of $f_L$ is unique. Thus $\pi$ is a diffeomorphism onto its open image $\Omega$. Let $\pi^{-1}(L)=(L,p_L)$, then $L\mapsto p_L$ is smooth. Now $\Omega$ is monodromy-invariant and the map $L\mapsto p_L$ is monodromy-equivariant. For every $z$, an $H$-isotropic line $\ell\subset F^2(z)$ exists because its signature is $(m-1,1)$. Extending $\ell$ to a maximal $H$-isotropic subspace gives $L\in\mathrm{IGr}_n(V,H)$ with $f_L(z)=0$. Thus $L\mapsto p_L$ maps $\Omega$ onto $\widetilde X$. Since $\widetilde X$ is noncompact, $\Omega$ cannot be the whole compact flag manifold $\mathrm{IGr}_n(V,H)$. A minimal nonempty closed invariant subset $\Lambda$ of its complement has the required transversality.

	For the Lyapunov--HN equality, we follow the idea of \cite[Remark~2.2.10 \& Lemmata~2.2.12--2.2.13]{Filip21}. Choose $o$ over a basepoint where the Oseledets theorem for $W_n$ and Birkhoff's theorem applied to the curvature density $c_1(\det F^2,h)/\mathrm dA_{\mathrm{hyp}}$ hold in almost every direction. We equip $\mathrm{IGr}_n(V,H)$ with a smooth positive volume measure. Note that the span of the complex lines $\bigwedge^nL$, for $L\in\mathrm{IGr}_n(V,H)$ is invariant under $\mathrm{SU}_{m,n}$, whose complex Zariski closure is $\mathrm{SL}_{m+n}\mathbb{C}$. Since $\bigwedge^n\mathbb C^{m+n}$ is an irreducible representation of $\mathrm{SL}_{m+n}\mathbb{C}$, this span is $W_n$. Consequently, for every proper complex linear subspace $A\subset W_n$, the set
	\[\left\{L\in\mathrm{IGr}_n(V,H):\bigwedge^nL\subset A\right\}\]
	is a proper real-analytic subset of the connected flag manifold $\mathrm{IGr}_n(V,H)$, and hence has measure zero.

	For almost every direction $\theta$ from $o$, the nonzero vectors with growth rate strictly smaller than $\lambda_1(W_n)$, together with zero, form a proper complex Oseledets subspace $A_\theta$. Hence for each such $\theta$, almost every $L$ satisfies $w_L\notin A_\theta$ and therefore realizes $\lambda_1(W_n)$ along that ray. Then, Fubini's theorem implies that almost every $L$ realizes this exponent along almost every ray from $o$. Since the map $\Omega\to\widetilde X,\ L\mapsto p_L$ is smooth and surjective, $\Omega\setminus\{L\in\Omega:p_L=o\}$ is nonempty and open. Thus, this open set has positive measure and we can choose such an $L\in\Omega$ with $p_L\neq o$. 
	We rescale $w_L$ so that $\|w_L(o)\|_h=1$. Then $\|w_L\|_{\det F^2}=\|w_L^{1,3n-1}\|_h=f_L^{1/2}$.
	Set $\varphi_L:=\log\|w_L\|_{\det F^2}$ away from $p_L$.
	The flatness equation and the equality $\|w_L\|_h=2^{n/2}\|w_L^{2,3n-2}\|_h$ give the pointwise identity
	\[\frac{\|w_L\|_h}{\|w_L\|_{\det F^2}}
	=\frac{2^{n/2+1}}{\|\Phi_{2,1}\|_{h,\mathrm{hyp}}}\|\partial\varphi_L\|_{\mathrm{hyp}}
	\leqslant\frac{2^{n/2+1}}{\delta}\|\nabla\varphi_L\|_{\mathrm{hyp}}.\]
	
	Near its unique minima $p_L$, the function $f_L$ satisfies $\|\nabla f_L\|\geqslant c_L f_L^{1/2}$ for some $c_L>0$ on a neighborhood $U$ since the minima is nondegenerate. By \cite[Theorem~3.1(2)]{LiZhang26}, any sequence with $f_L(z_j)\to0$ converges to $p_L$. Hence there is an $\eta_L>0$ such that $\{f_L<\eta_L\}\subset U$. Otherwise such a sequence exists outside $U$. On $\{f_L\geqslant\eta_L\}$, the previous estimate implies that
	\[\|\nabla f_L\|\geqslant c f_L\geqslant c\sqrt{\eta_L}\,f_L^{1/2}.\]
	Replacing $c_L$ by $\min\{c_L,c\sqrt{\eta_L}\}$ therefore makes $\|\nabla f_L\|\geqslant c_L f_L^{1/2}$ global. Note that the constant $c_L$ may depend on $L$.
	
	The exponential growth estimate \cite[Lemma~3.2]{LiZhang26}, see also \cite[Lemma~4.4]{Zhang25}, therefore gives
	\[d_{\mathrm{hyp}}(p_L,z)\leqslant A_L+c^{-1}\log(1+f_L(z))\]
	for some $A_L\geqslant0$. This constant may depend on the chosen $L$, the normalization at $o$, and the fixed variation. 
	
	Put
	\[r_L=A_L+c^{-1}\log2,\qquad R_L=d_{\mathrm{hyp}}(o,p_L)+r_L+3.\]
	When $d_{\mathrm{hyp}}(p_L,y)\geqslant r_L$, we obtain that $f_L(y)\geqslant1$. Let $C_0$ be the constant in Proposition~\ref{prop:adapted-metric} for $V$, and set $B=2nC_0$. Proposition~\ref{prop:adapted-metric}\romitem{2} gives 
	\[|\mathrm{d}\log\|w_L\|_h|_{\mathrm{hyp}}\leqslant B,\qquad \log\|w_L(y)\|_h\leqslant B\,d_{\mathrm{hyp}}(o,y).\]
	Here $B$ depends only on the fixed variation and $n$, and is independent of $L$ and $o$. If $d_{\mathrm{hyp}}(o,z)=R\geqslant R_L$, then every $y\in B_2(z)$ satisfies $d_{\mathrm{hyp}}(p_L,y)>r_L$. Thus this ball avoids $p_L$, and
	\[0\leqslant\varphi_L(y)\leqslant\log\|w_L(y)\|_h\leqslant B(R+2),\qquad y\in B_2(z).\]
	All balls here are taken in $\widetilde X$ for the lifted hyperbolic metric.
	
	Note that the Hitchin equation gives the Chern curvature $-\Psi\wedge\Psi^{*_h}$ on $W_n^{1,3n-1}$, where $\Psi=\Phi_{W_n}|_{W_n^{2,3n-2}}$. By \eqref{eq:signature-one-induced-current}, $\|\Psi\|_{h,\mathrm{hyp}}=\|\Phi_{2,1}\|_{h,\mathrm{hyp}}\leqslant M$.
	Consequently,
	\[K:=2\pi\sup_X\left|\frac{c_1(W_n^{1,3n-1},h)}{\mathrm{d}A_{\mathrm{hyp}}}\right|\leqslant C_{\mathrm{curv}}M^2,\]
	where $C_{\mathrm{curv}}$ is a numerical constant determined only by the norm conventions. Away from $p_L$, the nonvanishing holomorphic section $w_L^{1,3n-1}$ gives
	\[\Delta_{\mathrm{hyp}}\varphi_L=-2\pi\frac{c_1(W_n^{1,3n-1},h)}{\mathrm{d}A_{\mathrm{hyp}}},
	\qquad |\Delta_{\mathrm{hyp}}\varphi_L|\leqslant K.\]
	In particular, this bound is independent of $L$, $o$, $z$, and $R$.
	
	The interior gradient estimate for the Laplace--Beltrami operator on Riemannian manifolds, applied to $B_1(z)\subset B_2(z)$ (see, e.g., \cite[Theorem~2.4]{ChengThalmaierThompson18}), gives
	\[\|\nabla u(z)\|_{\mathrm{hyp}}\leqslant C_{\mathrm P}\left(\|u\|_{L^\infty(B_2(z))}+\|\Delta_{\mathrm{hyp}}u\|_{L^\infty(B_2(z))}\right).\]
	Here $C_{\mathrm P}$ depends only on the dimension $2$, the fixed radii $1$ and $2$, and the Ricci curvature lower bound $-4$. It is independent of $L$, $o$, $z$, and $R$, including over the cusps. 
	
	Applying this estimate to $\varphi_L$ and using the norm comparison above gives 
	\[1\leqslant\frac{\|w_L(z)\|_h}{\|w_L(z)\|_{\det F^2}}
	\leqslant\frac{2^{n/2+1}C_{\mathrm P}}{\delta}\left(B(R+2)+K\right)
	\leqslant C(1+R),\]
	for $d_{\mathrm{hyp}}(o,z)=R\geqslant R_L$, where one may take $C=\max\{1,2^{n/2+1}C_{\mathrm P}\delta^{-1}(2B+C_{\mathrm{curv}}M^2)\}$. This constant depends only on $n$, $\delta$, $M$, $C_0$, and the fixed norm conventions. 
	
	Writing $m_R$ for angular average on this circle, we obtain that
	\begin{equation}\label{eq:signature-one-tempered-comparison}
		0\leqslant m_R\left(\log\|w_L\|_h-\varphi_L\right)
		\leqslant\log C+\log(1+R)=O(\log R),\qquad R\geqslant R_L.
	\end{equation}
	Therefore,
	\begin{equation}\label{eq:Lyapunov-circular-average}
		\lambda_1(W_n)
	=\lim_{R\to\infty}\frac1R m_R\log\|w_L\|_h
	=\lim_{R\to\infty}\frac1R m_R\varphi_L.
	\end{equation}

	Let $m_L$ be the vanishing order of $w_L^{1,3n-1}$ at $p_L$. The Poincar\'e--Lelong formula (see, e.g., \cite[Chapter~V, equation~(13.2)]{Demailly12}) gives
	\[\Delta_{\mathrm{hyp}}\varphi_L\,\mathrm{d}A_{\mathrm{hyp}}=2\pi c_1(\det F^2,h)+2\pi m_L\delta_{p_L}.\]
	The Green formula used in the proof of \cite[Theorem~4.1]{EKMZ18} is
	\[\frac{\mathrm{d}}{\mathrm{d}t}m_tu=\frac{\tanh t}{2\operatorname{Area}(B_t(o))}\int_{B_t(o)}\Delta_{\mathrm{hyp}}u\,\mathrm{d}A_{\mathrm{hyp}}=\frac{1}{\pi\sinh(2t)}\int_{B_t(o)}\Delta_{\mathrm{hyp}}u\,\mathrm{d}A_{\mathrm{hyp}},\]
	since $\operatorname{Area}(B_t(o))=\pi\sinh^2t$ for curvature $-4$. Applying this formula to $\varphi_L$ we obtain that
	\[\frac{\mathrm{d}}{\mathrm{d}t}m_t\varphi_L=\frac{2}{\sinh(2t)}\left(\int_{B_t(o)}c_1(\det F^2,h)+m_L\mathbf1_{\{t>d_{\mathrm{hyp}}(o,p_L)\}}\right)\] 
	for almost every $t>0$. Integrating from $0$ to $R>d_{\mathrm{hyp}}(o,p_L)$, and using $\mathrm{d}A_{\mathrm{hyp}}=\tfrac12\sinh(2r)\,\mathrm{d}r\,\mathrm{d}\theta$ in polar coordinates about $o$ yields that
	\begin{equation}\label{eq:signature-one-circle-green}
		m_R\varphi_L-\varphi_L(o)
		=2\pi\int_0^R k_R(r)m_r\!\left(\frac{c_1(\det F^2,h)}{\mathrm{d}A_{\mathrm{hyp}}}\right)\,\mathrm{d}r+m_L\log\frac{\tanh R}{\tanh d_{\mathrm{hyp}}(o,p_L)},
	\end{equation}
	where
	\[k_R(r)=\sinh(2r)\int_r^R\frac{\mathrm{d}t}{\sinh(2t)},\qquad 0\leqslant r\leqslant R,\]
	with $k_R(0)=0$. In particular, the atom contributes exactly
	\[0\leqslant2m_L\int_{d_{\mathrm{hyp}}(o,p_L)}^R\frac{\mathrm{d}t}{\sinh(2t)}=m_L\log\frac{\tanh R}{\tanh d_{\mathrm{hyp}}(o,p_L)}\leqslant m_L\log\coth d_{\mathrm{hyp}}(o,p_L).\]
	The last inequality follows from $m_L\geqslant0$ and $\tanh R<1$.
	
	For $0<r\leqslant t$, the inequality $\sinh(2t)\geqslant\mathrm e^{2(t-r)}\sinh(2r)$ gives $0\leqslant k_R(r)\leqslant1/2$. And we have the integration $\int_0^R k_R(r)\mathrm{d}r=(\log\cosh R)/2$.
	Consequently, these yield that
	\[\int_0^R\left|k_R(r)-\frac12\right|\,\mathrm{d}r=\frac{R-\log\cosh R}{2}\leqslant\frac{\log2}{2}.\]
	The curvature bound above gives $\left|c_1(\det F^2,h)/\mathrm{d}A_{\mathrm{hyp}}\right|\leqslant K/(2\pi)$, so \eqref{eq:signature-one-circle-green} implies
	\begin{equation}\label{eq:error}
	\left|m_R\varphi_L-\varphi_L(o)-\pi\int_0^R m_r\!\left(\frac{c_1(\det F^2,h)}{\mathrm{d}A_{\mathrm{hyp}}}\right)\,\mathrm{d}r\right|
	\leqslant\frac{K}{2}\log2+m_L\log\coth d_{\mathrm{hyp}}(o,p_L).
	\end{equation}
	
	Birkhoff's theorem implies that, for $\sigma_o$-almost every $v\in T_o^1\widetilde X$, the ray $r\mapsto\exp_o(rv)$ satisfies
	\begin{equation}\label{eq:Birkhoff}
	\lim_{R\to\infty}\frac1R\int_0^R\frac{c_1(\det F^2,h)}{\mathrm{d}A_{\mathrm{hyp}}}\left(\exp_o(rv)\right)\,\mathrm{d}r=\int_X\frac{c_1(\det F^2,h)}{\mathrm{d}A_{\mathrm{hyp}}}\,\mathrm{d}\mu_X=\frac1{\operatorname{Area}(X)}\int_Xc_1(\det F^2,h).
	\end{equation}
	Since $m_r$ is normalized angular average, Fubini's theorem gives, for every $R>0$,
	\[\frac1R\int_0^R m_r\!\left(\frac{c_1(\det F^2,h)}{\mathrm{d}A_{\mathrm{hyp}}}\right)\,\mathrm{d}r=\int_{T_o^1\widetilde X}\left[\frac1R\int_0^R\frac{c_1(\det F^2,h)}{\mathrm{d}A_{\mathrm{hyp}}}\left(\exp_o(rv)\right)\,\mathrm{d}r\right]\mathrm{d}\sigma_o(v).\]
	As $\sigma_o$ is a probability measure, dominated convergence and \eqref{eq:Birkhoff} yield
	\begin{equation}\label{eq:circular-average-limit}
	\lim_{R\to\infty}\frac1R\int_0^R m_r\!\left(\frac{c_1(\det F^2,h)}{\mathrm{d}A_{\mathrm{hyp}}}\right)\,\mathrm{d}r
	=\frac1{\operatorname{Area}(X)}\int_X c_1(\det F^2,h).
	\end{equation}
	Therefore, \eqref{eq:Lyapunov-circular-average}, \eqref{eq:error}, and \eqref{eq:circular-average-limit} yield that
	\[\lambda_1(\mathbb{V})+\cdots+\lambda_n(\mathbb{V})=\frac{\pi}{\operatorname{Area}(X)}\int_X c_1(\det F^2,h)=\frac{2\operatorname{pardeg}F^2}{\deg\KK_{\overline{X}}(D)},\]
	where the last equality is given by the parabolic Chern--Weil formula. Since the monodromy lies in $\mathrm{SU}_{m,n}$ with $m\geqslant n$, its middle Lyapunov exponents $\lambda_{n+1},\ldots,\lambda_m$ vanish when $m>n$. Thus
	\[\sum_{i=1}^m\lambda_i(\mathbb V)=\sum_{i=1}^n\lambda_i(\mathbb V),\]
	which proves \eqref{eq:signature-one-contact-current} with $\operatorname{rank}F^2=m$.
	
	To prove the nonlattice statement, apply Theorem~\ref{thm:unitary-contact} with $\mathcal{E}_m=F^2$, using \eqref{eq:signature-one-contact-current}. Since $H|_{F^2}$ has signature $(m-1,1)$, the definite-plane exception is impossible, so $\Gamma$ is not a lattice.
	
	In the real symplectic case, restrict the previous construction to $\mathrm{LGr}_n(V_{\mathbb R})$. Viewing each real Lagrangian through its complexification, these are $H$-isotropic complex $n$-planes, so the same uniqueness and nondegeneracy arguments apply. And all other statements follow directly.
\end{proof}

\begin{remark}
	The dependence of the constant $A_L$ on $L$ is analogous to \cite[Lemma~2.2.13]{Filip21}, where uniformity is asserted on compact families whose associated functions attain their minima. 
\end{remark}

\subsection{Geometric consequences}\label{subsec:weight-three-geometry}

The construction and norm estimates above also give geometric consequences, following the methods of \cite[\S3.4 and Theorem~4.1]{LiZhang26}.

\begin{corollary}\label{cor:weight-three-geometry}
	Under the hypotheses of Theorem~\ref{thm:main-C}, let $\rho\colon\pi_1(X)\to\mathrm{SU}_{m,n}$ be the monodromy and $\Gamma=\rho(\pi_1(X))$. Then $\rho$ is discrete and faithful. Moreover:
	\begin{enumerate}[label=\textup{(\roman*)}]
		\item The proper nonempty open set
		\[\Omega=\left\{L\in\mathrm{IGr}_n(V,H):F^2(z)\cap L\neq0\text{ for some }z\in\widetilde X\right\}\]
		is a domain of discontinuity: $\Gamma$ acts freely and properly discontinuously on $\Omega$. The map sending $L$ to the unique zero $p_L$ of $f_L$ induces a smooth fibre bundle $\Omega/\Gamma\to X$ with compact fibres. In particular, $\Omega/\Gamma$ is compact if and only if $X$ is compact.
		
		\item Let $u\colon\widetilde X\to\mathcal S=\mathrm{SU}_{m,n}/\mathrm S(\mathrm U_m\times\mathrm U_n)$ be the equivariant harmonic map defined by the Hodge metric, and fix an invariant Riemannian metric on $\mathcal S$. Then $u$ is a bi-Lipschitz embedding: there are constants $c>0$ and $C>0$ such that
		\[c\,d_{\mathrm{hyp}}(x,y)
		\leqslant d_{\mathcal S}(u(x),u(y))
		\leqslant C\,d_{\mathrm{hyp}}(x,y)
		\qquad(x,y\in\widetilde X).\]
		These constants depend only on $m,n$, the resulting positive lower bound for $\|\Phi_{2,1}\|_{h,\mathrm{hyp}}$, the bound for $\|\Phi\|_{h,\mathrm{hyp}}$, and the fixed metric conventions.
		If $X$ is compact, then the orbit map of $\rho$ is a quasi-isometric embedding with respect to any word metric on $\pi_1(X)$.
	\end{enumerate}
\end{corollary}

\begin{proof}
	The proof of Theorem~\ref{thm:main-C} identifies the manifold $\mathcal I$ diffeomorphically with $\Omega$ and gives the smooth, equivariant, surjective map
	\[\begin{aligned}
		p\colon\Omega&\longrightarrow\widetilde X\\
		L&\longmapsto p_L.
	\end{aligned}\]
	For every compact $A\subset\widetilde X$, the set $\mathcal I\cap(\mathrm{IGr}_n(V,H)\times A)$ is closed in a compact space. Its image is $p^{-1}(A)$, so $p$ is proper. Locally choose smooth transformations $g(z)$ in $\mathrm{SU}_{m,n}$, taking the polarized Hodge decomposition at a fixed point $z_0$ to that at $z$. Then $(L,z)\mapsto(g(z)^{-1}L,z)$ identifies $\mathcal{I}$ locally with a product over the coordinate disk. These identifications provide smooth lifts of every base tangent vector, so $p$ is a submersion and a smooth fibre bundle with compact fibres.
	
	If $A\subset\Omega$ is compact and $\rho(\gamma)A\cap A\neq\varnothing$, then $\gamma p(A)\cap p(A)\neq\varnothing$. Only finitely many deck transformations have this property, so the action on $\Omega$ is properly discontinuous. A fixed point of $\rho(\gamma)$ in $\Omega$ would give a fixed point of $\gamma$ in $\widetilde X$, hence $\gamma=1$. This also proves faithfulness. Proper discontinuity on the nonempty open set $\Omega$ implies that $\Gamma$ is discrete: a sequence of distinct elements tending to the identity would move a fixed point of $\Omega$ into a compact neighbourhood of that point. The quotient map $\Omega/\Gamma\to X$ inherits the local bundle description and is proper and surjective. This proves \romitem{1}.
	
	For \romitem{2}, fix any $x\in\widetilde X$, and let $P=V^{0,3}(x)\oplus V^{2,1}(x)$ be regarded as a flat complex $n$-plane. Choose a flat $w\in\det P$ with $\|w(x)\|_h=1$, and set $f(y)=\|w^{1,3n-1}(y)\|_h^2$. The form $H$ is negative definite on $P$. At any $y$, projection of $P$ onto $L_0(y)=V^{0,3}(y)\oplus V^{2,1}(y)$ is an isomorphism, so $P$ is the graph of a strict contraction $U\colon L_0(y)\to P_0(y)$. In the adapted orthonormal basis used in the proof of Theorem~\ref{thm:main-C}, we know $w=a(e_1+Ue_1)\wedge\cdots\wedge(e_n+Ue_n)$.
	Flatness of $H$ and the normalization at $x$ give
	\[|a|^2\det(I-U^*U)=1,\qquad |a|\geqslant1.\]
	The same computation gives
	\[f(y)=|aU_{mn}|^2\leqslant |a|^2,\qquad \|w^{2,3n-2}(y)\|_h=|a|.\]
	Since all singular values of $U$ are smaller than $1$,
	\[\det(I-U^*U)\leqslant1-\|U\|_{\mathrm{op}}^2\leqslant1-|U_{mn}|^2.\]
	The determinant identity therefore gives $|a|^2\geqslant1+f$. The flatness identities now implies that 
	\begin{equation}\label{eq:lower-bound}
	\|\nabla f\|_{\mathrm{hyp}}
	\geqslant\varepsilon |a|\sqrt f
	\geqslant\varepsilon\sqrt{f(1+f)},
	\end{equation}
	where $\varepsilon>0$ depends on the lower bound for $\|\Phi_{2,1}\|_{h,\mathrm{hyp}}$.
	At a zero, the same Hessian calculation as before is positive definite. Thus $f$ is a nonnegative Morse function with $f(x)=0$. Then \eqref{eq:lower-bound} implies $\|\nabla f\|\geqslant\varepsilon f$ and $\|\nabla f\|\geqslant\varepsilon\sqrt f$. By \cite[Theorem~3.1]{LiZhang26}, the negative gradient trajectory from any $y$ converges to the unique zero $x$. Along this trajectory, $-\mathrm{d}f/\mathrm{d}t=\|\nabla f\|^2$.
	This gives
	\[d_{\mathrm{hyp}}(x,y)
	\leqslant\frac1\varepsilon\int_0^{f(y)}\frac{\mathrm{d}s}{\sqrt{s(1+s)}}
	=\frac2\varepsilon\operatorname{arcsinh}\sqrt{f(y)}.\]
	We compare this bound with the symmetric-space distance using the $KAK$ argument of \cite[Theorem~5.1.4]{Filip21}, \cite[\S5.2]{Zhang25}, and \cite[Theorem~3.9]{LiZhang26}.
	Let $K_\mathbb{R}\cong\mathrm S(\mathrm U_m\times\mathrm U_n)$ be the stabilizer of $u(y)=L_0(y)$, and choose $g\in\mathrm{SU}_{m,n}$ with $gL_0(y)=P$. Write its Cartan decomposition as
	\[g=k_1a_tk_2,\qquad k_1,k_2\in K_\mathbb{R}.\]
	$a_t$ can be expressed as follows: choose orthonormal bases $e_1,\ldots,e_n$ of $L_0(y)$ and $f_1,\ldots,f_m$ of $P_0(y)$ by
	\[a_te_j=\cosh t_j\,e_j+\sinh t_j\,f_j,
	\qquad
	a_tf_j=\sinh t_j\,e_j+\cosh t_j\,f_j,\]
	for $1\leqslant j\leqslant n$, where $t_1\geqslant\cdots\geqslant t_n\geqslant0$. It fixes $f_{n+1},\ldots,f_m$. The vectors $e_j+f_j$ and $e_j-f_j$ are eigenvectors with eigenvalues $\mathrm e^{t_j}$ and $\mathrm e^{-t_j}$, respectively.
	Put $w_0=e_1\wedge\cdots\wedge e_n$. Since $g$ preserves $H$, both $(\bigwedge^n g)w_0$ and $w$ have unit norm for $\det(-H|_P)$, so they differ by a scalar of modulus one. The compact factors preserve $h_y$, and $k_2$ acts on $\det L_0(y)$ by a scalar of modulus one.
	Thus,
	\[\|w(y)\|_h^2
	=\left\|\bigwedge_{j=1}^n(\cosh t_j\,e_j+\sinh t_j\,f_j)\right\|_h^2
	=\prod_{j=1}^n\cosh(2t_j).\]
	Now $|a|^2\geqslant1+f(y)$ implies that
	\[1+2f(y)\leqslant |a|^2+f(y)\leqslant\|w(y)\|_h^2=\prod_{j=1}^n\cosh(2t_j)\leqslant\cosh\!\left(2\sum_{j=1}^n t_j\right),\]
	where the last inequality follows from $\cosh(s+t)\geqslant\cosh s\cosh t$ for $s,t\geqslant0$. Therefore, $\operatorname{arcsinh}\sqrt{f(y)}\leqslant\sum_jt_j$. Note that The invariant distance between $u(x)$ and $u(y)$ is a fixed positive multiple of $(\sum_jt_j^2)^{1/2}$ since $g\cdot u(y)=u(x)$. Hence Cauchy--Schwarz inequality shows that
	\[\sum_{j=1}^n t_j\leqslant C_{m,n} d_{\mathcal S}(u(x),u(y)),\]
	where $C_{m,n}>0$ depends only on $m,n$ and the fixed invariant metric on $\mathcal S$. Combining the estimates yields that
	\[\frac\varepsilon2 d_{\mathrm{hyp}}(x,y)
	\leqslant\operatorname{arcsinh}\sqrt{f(y)}
	\leqslant\sum_{j=1}^n t_j
	\leqslant C_{m,n} d_{\mathcal S}(u(x),u(y)).\]
	Thus the lower part of bi-Lipschitz holds. The upper bound follows directly from the bound of the Higgs field, see Proposition~\ref{prop:adapted-metric}\romitem{1} gives the upper Lipschitz bound. Hence $u$ is proper and injective. $u$ is also an immersion because $\Phi_{2,1}$ is nowhere-vanishing. So $u$ is a bi-Lipschitz embedding. When $X$ is compact, the orbit map $\pi_1(X)\to\widetilde X$ is a quasi-isometry, and composition with $u$ proves the quasi-isometric embedding.
\end{proof}

If $m=n$ and the variation arises from a real polarized variation, then the harmonic map takes values in the totally geodesic subspace $\mathrm{Sp}_{2n}\mathbb R/\mathrm U_n$. Thus the statement about bi-Lipschitz embedding also holds for this symmetric space.

\subsection{Orthogonal variations and Assumption B}\label{subsec:orthogonal-B}

Let $\mathbb{V}$ be a polarized real weight-four variation on $X=\overline{X}\setminus D$, with the tame hypotheses of Section~\ref{sec:setup} and Hodge numbers
\[\left(h^{4,0},h^{3,1},h^{2,2},h^{1,3},h^{0,4}\right)
=(n,1,1,1,n),\qquad n\geqslant1.\]
Let $Q$ be the invariant real quadratic form with signature $(2,2n+1)$. Its positive real subspace is $(V^{3,1}\oplus V^{1,3})\cap V_{\mathbb R}$. We assume that the monodromy takes values in $\mathrm{SO}_{2,2n+1}$ and that the
logarithmic extension
\[
\Phi_{2,2}\colon V^{2,2}_\bullet\longrightarrow
V^{1,3}_\bullet\otimes\KK_{\overline X}(D)
\]
is an isomorphism of filtered line bundles on $\overline X$. This extends the rank-five setting of Assumption~B in \cite[Definition~5.4.1]{Filip21}. Since $V^{2,2}$ and $V^{1,3}$ are line bundles, a similar argument in the proof of Theorem~\ref{thm:main-C} gives
\[
\|\Phi_{2,2}\|_{h,\mathrm{hyp}}\geqslant\delta>0
\quad\text{on }X.
\]
We show directly how the estimates of \cite{LiZhang26} apply in this generality.

We have 
\[
\|\Phi_{3,1}\|_{h,\mathrm{hyp}}=\|\Phi_{2,2}\|_{h,\mathrm{hyp}},
\qquad
\|\Phi_{1,3}\|_{h,\mathrm{hyp}}=\|\Phi_{4,0}\|_{h,\mathrm{hyp}}.
\]
by the real structure. We consider the cyclic ordering
\[
(E_0,E_1,E_2,E_3)
=\left(V^{2,2},V^{1,3},V^{0,4}\oplus V^{4,0},V^{3,1}\right).
\]
The ranks are $(1,1,2n,1)$ and the composition $E_1\to E_2\otimes\KK_X\to E_3\otimes\KK_X^2$ vanishes. Note that the metric
\[
g=\|\Phi_{2,2}\|_{h,\mathrm{hyp}}^2g_{\mathrm{hyp}}
\geqslant\delta^2g_{\mathrm{hyp}}
\]
is complete. Hence 
\[
r:=\frac{\|\Phi_{4,0}\|_{h,\mathrm{hyp}}^2}
{\|\Phi_{2,2}\|_{h,\mathrm{hyp}}^2}.
\]
satisfies that
\[
K_g=2r-2\geqslant-2,\qquad\frac14\Delta_g\log r\geqslant3r-2
\]
on $\{r>0\}$ by the computations in \cite[proof of Lemma~2.5]{LiZhang26}. The identity $\Delta_g\log r=(\Delta_g r)/r-\|\mathrm{d}r\|_g^2/r^2$ therefore implies
\[
\frac14\Delta_g r\geqslant r(3r-2).
\]
The inequality above holds everywhere by continuity. Moreover, $r$ is bounded above by Proposition~\ref{prop:adapted-metric}\romitem{1} and $\|\Phi_{2,2}\|_{h,\mathrm{hyp}}\geqslant\delta$. Therefore, by the Omori--Yau maximum principle \cite{omori1967isometric,yau1975harmonic} with respect to the complete metric $g$ with $K_g\geqslant-2$, we obtain that
\[
r\leqslant\frac23,
\qquad
\|\Phi_{2,2}\|_{h,\mathrm{hyp}}-\|\Phi_{4,0}\|_{h,\mathrm{hyp}}
\geqslant\left(1-\sqrt{2/3}\right)\delta>0.
\]
This is the gap required in \cite[Assumption~3.1]{LiZhang26}.

For a real $Q$-isotropic two-plane $S$, extend a nonzero $w_S\in\bigwedge^2S$ flatly and put
\[
f_S=\left\|\operatorname{pr}_{V^{2,2}\wedge V^{3,1}}w_S\right\|_h^2.
\]
Since $S_{\mathbb C}$ is isotropic for $Q(\cdot,\overline{\cdot})$, the gap and \cite[Theorem~3.1 \& Lemmata~3.6--3.7]{LiZhang26} show that $f_S$ has at most one zero, necessarily nondegenerate.

\begin{corollary}\label{cor:orthogonal-B-thinness}
	Under the preceding hypotheses, the monodromy image $\Gamma$ is not a lattice in $\mathrm{SO}(Q)$. If it is contained in an arithmetic lattice $\Gamma_{\mathrm{ar}}$ of $\mathrm{SO}(Q)$, then $[\Gamma_{\mathrm{ar}}:\Gamma]=\infty$; if it is also Zariski dense in $\mathrm{SO}(Q)$, it is thin.
\end{corollary}

\begin{proof}
	As in the proof of Theorem~\ref{thm:main-C}, consider the real flag manifold $\mathrm{IGr}_2(V_{\mathbb R},Q)$ and its subset
	\[\Omega=\left\{S\in\mathrm{IGr}_2(V_{\mathbb R},Q):f_S(z)=0\text{ for some }z\in\widetilde X\right\}.\]
	Uniqueness defines an equivariant map $p\colon\Omega\to\widetilde X$, sending $S$ to its zero. Writing the projection defining $f_S$ in a local unit frame, nondegeneracy and the implicit function theorem show that $\Omega$ is open and $p$ is smooth, exactly as in the proof of Theorem~\ref{thm:main-C}.
	
	For every $z$, choose unit vectors $e_0\in V^{0,4}(z)$ and $e_3\in V^{3,1}(z)$. The real span of $\operatorname{Re}(e_0+e_3)$ and $\operatorname{Im}(e_0+e_3)$ is $Q$-isotropic and has $f_S(z)=0$. Thus $p$ is surjective. Since $\widetilde X$ is noncompact and the flag manifold is compact, $\Omega$ is proper. Its complement is therefore a proper closed $\Gamma$-invariant set.
	
	The Pl\"ucker embedding realizes this flag manifold as a compact orbit in the projectivized irreducible adjoint representation $\bigwedge^2V_{\mathbb C}$. If $\Gamma$ were a lattice in $\mathrm{SO}(Q)$, its intersection with the identity component would also be a lattice. Proposition~\ref{prop:lattice-minimality}, applied to this intersection, would make its action on the flag manifold minimal, a contradiction. Thus $\Gamma$ is not a lattice.
	If it is contained in an arithmetic lattice, it has infinite index; if it is also Zariski dense in $\mathrm{SO}(Q)$, it is thin.
\end{proof}

Let $F^3=V^{4,0}\oplus V^{3,1}$. The following proposition is stated in \cite[\S5.4.4]{Filip21} using geometric structure studied in \cite{collier2019geometry}. Here we give a proof similar to Theorem~\ref{thm:main-C}.

\begin{proposition}\label{prop:orthogonal-B-equality}
	Under the preceding hypotheses, if $n=1$, then
	\[\lambda_1(\mathbb{V})+\lambda_2(\mathbb{V})=\frac{2\operatorname{pardeg}F^3}{\deg\KK_{\overline{X}}(D)}.\]
\end{proposition}

\begin{proof}
	Let $S$ be a real $Q$-isotropic two-plane and extend a nonzero $w\in\bigwedge^2S$ flatly over $\widetilde X$. Write $w^{p,q}$ for its Hodge components. The two lowest Hodge summands of $\bigwedge^2V$ are the line bundles
	\[\left(\bigwedge^2V\right)^{1,7}=V^{0,4}\wedge V^{1,3},\qquad
	\left(\bigwedge^2V\right)^{2,6}=V^{0,4}\wedge V^{2,2}.\]
	Then we set $f_S=\|\operatorname{pr}_{V^{2,2}\wedge V^{3,1}}w\|_h^2$.
	
	As in Theorem~\ref{thm:main-C}, projection onto the positive subspace $V^{3,1}\oplus V^{1,3}$ identifies $S_{\mathbb C}$ with the graph of a linear map
	\[U\colon V^{3,1}\oplus V^{1,3}\longrightarrow
	V^{0,4}\oplus V^{2,2}\oplus V^{4,0}.\]
	Isotropy for $Q(\cdot,\overline{\cdot})$ gives $U^*U=I_2$, and the real structure implies that $U$ commutes with conjugation.
	Choose unit vectors $e_p\in V^{p,4-p}$ with $\overline e_p=e_{4-p}$, and write $(U_{ij})$ for the matrix of $U$ in the ordered bases $(e_3,e_1)$ and $(e_0,e_2,e_4)$.
	Then
	\[(U_{12},U_{22},U_{32})=(\overline{U_{31}},\overline{U_{21}},\overline{U_{11}}).\]
	The Gram matrix of $(e_3+Ue_3,e_1+Ue_1)$ is $I_2+U^*U=2I_2$. Hence, writing $w=c(e_3+Ue_3)\wedge(e_1+Ue_1)$, we have $\|w\|_h^2=4|c|^2$. Expanding the exterior product gives
	\[w^{1,7}=cU_{11}\,e_0\wedge e_1,\qquad w^{2,6}=c(U_{11}\overline{U_{21}}-U_{21}\overline{U_{31}})\,e_0\wedge e_2,\qquad f_S=|c|^2|U_{21}|^2.\]
	The identities $U^*U=I_2$ give
	\[|U_{11}|^2+|U_{21}|^2+|U_{31}|^2=1,\qquad 2U_{11}U_{31}+U_{21}^2=0.\]
	Taking moduli in the second identity and substituting into the first yields $|U_{11}|+|U_{31}|=1$, and consequently $|U_{21}|^2=2|U_{11}|(1-|U_{11}|)$. The same identities give
	\[|U_{11}\overline{U_{21}}-U_{21}\overline{U_{31}}|^2=|U_{21}|^2\left(|U_{11}|^2+|U_{21}|^2+|U_{31}|^2\right)=|U_{21}|^2.\]
	Thus $\|w^{2,6}\|_h^2=f_S$. The pairing induced by $Q$ identifies $\det F^3$ with a line in $(\bigwedge^2V)^\vee$, and its evaluation seminorm satisfies $\|w\|_{\det F^3}=\|w^{1,7}\|_h$. We obtain
	\begin{equation}\label{eq:orthogonal-B-real-comparison}
		\|w\|_h\|w\|_{\det F^3}=2\|w\|_{\det F^3}^2+f_S.
	\end{equation}
	The lowest component $w^{1,7}$ is holomorphic in $(\det F^3)^\vee$.
	
	For $S\in\Omega$, write $p_S=p(S)$ for the unique zero of $f_S$. The planes constructed in Corollary~\ref{cor:orthogonal-B-thinness} satisfy $w^{1,7}(p_S)\neq0$, so this condition defines a nonempty open subset of $\Omega$. For these $S$, \eqref{eq:orthogonal-B-real-comparison} and uniqueness of $p_S$ show that $w^{1,7}$ is nowhere zero. Moreover, \cite[Lemma~3.6]{LiZhang26} gives $\|\mathrm{d}f_S\|_{\mathrm{hyp}}\geqslant\varepsilon f_S$.
	Since $p_S$ is a nondegenerate zero, a similar argument to that in the proof of Theorem~\ref{thm:main-C} gives $\|\mathrm{d}f_S\|_{\mathrm{hyp}}\geqslant c_S\sqrt{f_S}$ globally, with $c_S>0$ allowed to depend on $S$. This is the same dependence on the chosen plane discussed in Theorem~\ref{thm:main-C}, analogous to the uniformity on compact families in \cite[Lemma~2.2.13]{Filip21}.
	Thus \cite[Lemma~3.2]{LiZhang26} applies to this isotropic plane, and $f_S$ grows exponentially away from $p_S$.
	The Pl\"ucker lines of this open family span $\bigwedge^2V$. A similar argument to that in the proof of Theorem~\ref{thm:main-C} therefore selects $S$ and a basepoint $o$ such that $w$ realizes $\lambda_1(\mathbb{V})+\lambda_2(\mathbb{V})$ along almost every ray from $o$. Rescale $w$ so that $\|w(o)\|_h=1$, and put $\varphi=\log\|w\|_{\det F^3}$. Proposition~\ref{prop:adapted-metric}\romitem{2}, applied to $\bigwedge^2V$, gives a constant $B>0$, independent of $S$, such that
	\[\log\|w(y)\|_h\leqslant B\,d_{\mathrm{hyp}}(o,y).\]
	Exponential growth gives $r_S>0$ such that $f_S(y)\geqslant1$ whenever $d_{\mathrm{hyp}}(p_S,y)\geqslant r_S$. If $d_{\mathrm{hyp}}(o,z)=R\geqslant d_{\mathrm{hyp}}(o,p_S)+r_S+2$, then every $y\in B_2(z)$ satisfies $f_S(y)\geqslant1$ and $d_{\mathrm{hyp}}(o,y)\leqslant R+2$.
	The inequalities $f_S\leqslant\|w\|_h\|w\|_{\det F^3}$ and $\|w\|_{\det F^3}\leqslant\|w\|_h$ therefore give
	\[\mathrm{e}^{-B(R+2)}\leqslant\|w(y)\|_h^{-1}
	\leqslant\|w(y)\|_{\det F^3}\leqslant\|w(y)\|_h
	\leqslant \mathrm{e}^{B(R+2)}.\]
	Thus $\|\varphi\|_{L^\infty(B_2(z))}\leqslant B(R+2)$. Since its Laplacian is bounded, the interior gradient estimate gives $\|\mathrm{d}\varphi(z)\|_{\mathrm{hyp}}\leqslant C_1(1+R)$, with $C_1$ independent of $z$.

	Flatness gives $D_h^{1,0}w^{1,7}=-(\mathrm{id}\wedge\Phi_{2,2})w^{2,6}$. Since $w^{1,7}$ is a nowhere-zero holomorphic frame and $\|w^{1,7}\|_h^2=e^{2\varphi}$, its Chern connection also satisfies $D_h^{1,0}w^{1,7}=2\partial\varphi\otimes w^{1,7}$.
	Therefore,
\[
4\|w\|_{\det F^3}^2\|\partial\varphi\|_{\mathrm{hyp}}^2
=\|D_h^{1,0}w^{1,7}\|_{h,\mathrm{hyp}}^2
=\|\Phi_{2,2}\|_{h,\mathrm{hyp}}^2\|w^{2,6}\|_h^2
=\|\Phi_{2,2}\|_{h,\mathrm{hyp}}^2f_S.
\]
Dividing \eqref{eq:orthogonal-B-real-comparison} by $\|w\|_{\det F^3}^2$ and using $\|\Phi_{2,2}\|_{h,\mathrm{hyp}}\geqslant\delta$ therefore yields
\[2\leqslant\frac{\|w(z)\|_h}{\|w(z)\|_{\det F^3}}=2+\frac{4}{\|\Phi_{2,2}(z)\|_{h,\mathrm{hyp}}^2}\|\partial\varphi(z)\|_{\mathrm{hyp}}^2\leqslant2+\frac{4C_1^2}{\delta^2}(1+R)^2\leqslant C(1+R)^2,\]
where $C=2+4C_1^2/\delta^2$. Taking logarithms gives
\[
\log2\leqslant\log\|w(z)\|_h-\varphi(z)
\leqslant\log C+2\log(1+R)=O(\log R).
\]
The angular Green formula and parabolic Chern--Weil argument from Theorem~\ref{thm:main-C} now yield
\[
\lambda_1(\mathbb{V})+\lambda_2(\mathbb{V})
=\frac{\pi}{\operatorname{Area}(X)}\int_Xc_1(\det F^3,h)
=\frac{2\operatorname{pardeg}F^3}{\deg\KK_{\overline{X}}(D)}.
\]
\end{proof}

\begin{proposition}
\label{prop:orthogonal-B-strictness}
Under the preceding hypotheses, suppose that $n>1$ and that the monodromy is Zariski dense in $\mathrm{SO}(Q)$. Then
\[
\sum_{j=1}^{n+1}\lambda_j(\mathbb{V})
>\frac{2\operatorname{pardeg}F^3}{\deg\KK_{\overline{X}}(D)}.
\]
\end{proposition}

\begin{proof}
Let $W=\bigwedge^{n+1}V_{\mathbb C}$. The pairing induced by $Q$ identifies $\mathcal A=\det F^3$ with a holomorphic line subbundle of $\mathcal W^\vee$ with the induced parabolic structure. The representation $\bigwedge^{n+1}\mathbb C^{2n+3}$ of $\mathrm{SO}_{2n+3}\mathbb C$ is irreducible, so Zariski density makes $W$ strongly irreducible. The common kernel of the equivariant family $\mathcal A_z$ is a proper invariant subspace and therefore vanishes. Thus Assumption~\ref{ass:no-common-kernel} holds. Moreover,
\[
\lambda_1(W)
=\sum_{j=1}^{n+1}\lambda_j(\mathbb V).
\]
By Theorems~\ref{thm:main-A} and \ref{thm:strict-zero-divisor} and Proposition~\ref{prop:lyap-coincidence}\romitem{2}, it suffices to construct a normalized harmonic measure $\widehat\nu$ on $\mathbb P(W)$ such that $\operatorname{supp}\widehat\nu\cap\Sigma_{\mathcal A}
\neq\varnothing$.

Let $\mathcal{B}\to X$ be the flat bundle of real $Q$-isotropic two-planes. Its Pl\"ucker image is a nonempty closed invariant subset of $\mathbb P(\bigwedge^2V_{\mathbb C})$. The adjoint representation $\bigwedge^2V_{\mathbb C}$ is strongly irreducible by Zariski density, so Proposition~\ref{prop:lyap-coincidence}\romitem{2} supplies a normalized harmonic probability measure $\nu$ on $\mathcal{B}$.

For $(x,S)\in\mathcal B$, the quotient $S^{\perp_Q}/S$ is negative definite of real dimension $2n-1$. Set
\[
C_S=\mathrm{IGr}_{n-1}
\left((S^{\perp_Q}/S)_{\mathbb C},Q\right).
\]
Taking inverse images in $S^{\perp_Q}_{\mathbb C}$ identifies $C_S$ with the variety of maximal complex isotropic $(n+1)$-planes containing $S_{\mathbb C}$.
We regard it as a projective subvariety of $\mathbb P(W_x)$ through the Pl\"ucker embedding. Its complex dimension is
\[
\dim_{\mathbb C}C_S=\frac{n(n-1)}2>0.
\]

The compact orthogonal group of $S^{\perp_Q}/S$ acts
transitively on $C_S$, giving a canonical invariant probability measure $m_S$ with full support, which depends continuously on $(x,S)$. These measures are preserved by flat transport, which induces isometries of the quotients.

For a bounded continuous function $f$ on $\mathbb P(W)$, define
\[
Tf(x,S)=\int_{C_S}f(x,\ell)\,\mathrm{d}m_S(\ell),
\qquad
\int f\,\mathrm{d}\widehat\nu
=\int_{\mathcal B}Tf\,\mathrm{d}\nu.
\]
This defines a probability measure with $p_*\widehat\nu=\mu_X$. Since $\nu$ is harmonic, $\widehat\nu$ is also harmonic.

We next verify that $C_S\subset\operatorname{supp}\widehat\nu$ when $(x,S)\in\operatorname{supp}\nu$. Fix $\ell\in C_S$ and an open neighborhood $U$ of $(x,\ell)$ in $\mathbb P(W)$. Choose a nonnegative continuous function $f$, compactly supported in $U$, with $f(x,\ell)>0$.
Full support of $m_S$ gives $Tf(x,S)>0$. By continuity, $Tf$ is bounded below by a positive constant on a neighborhood of $(x,S)$.
This neighborhood has positive $\nu$-measure, so
\[
\int f\,\mathrm{d}\widehat\nu
=\int Tf\,\mathrm{d}\nu>0.
\]
Consequently $\widehat\nu(U)>0$, proving the support inclusion.

Every positive-dimensional projective variety meets every projective hyperplane. Therefore, for any $(x,S)\in\operatorname{supp}\nu$, $C_S\cap\mathbb P(\ker\mathcal A_x)\neq\varnothing$.
Together with the support inclusion, this gives $\operatorname{supp}\widehat\nu\cap\Sigma_{\mathcal A}
\neq\varnothing$.
Theorem~\ref{thm:strict-zero-divisor} yields $\delta_{\mathcal A}(\widehat\nu)>0$. Theorem~\ref{thm:main-A} and Proposition~\ref{prop:lyap-coincidence}\romitem{2} therefore give
\[
\frac{2\operatorname{pardeg}F^3}{\deg\KK_{\overline X}(D)}=\lambda(\widehat\nu)-\delta_{\mathcal A}(\widehat\nu)<\lambda(\widehat\nu)\leqslant\lambda_1(W)=\sum_{j=1}^{n+1}\lambda_j(\mathbb V).
\]
This proves the proposition.
\end{proof}

\begin{remark}
 
For $n=1$, the variety $C_S$ in the proof is a point, so this obstruction disappears. For $n>1$, Corollary~\ref{cor:orthogonal-B-thinness} and Proposition~\ref{prop:orthogonal-B-strictness} show that Zariski-dense monodromy is not a lattice and the Lyapunov--HN inequality is strict. If the monodromy is contained in an arithmetic lattice, it is therefore thin.

\end{remark}

\subsection{C-32 in \texorpdfstring{$\mathrm{Sp}_6\mathbb R$}{Sp(6,R)}}\label{subsec:C32-example}

Let $A$ and $B$ be the companion matrices of
\[
f(t)=(t-1)^4(t^2-t+1),\qquad g(t)=t^6+1.
\]
We define C-32 to be the rank-six hypergeometric local system on $X=\mathbb P^1\setminus\{0,1,\infty\}$ with local monodromies $T_0=B^{-1}$, $T_1=A^{-1}B$, and $T_\infty=A$.
This is case C-32 in the degree-six symplectic census \cite{BDSS21,BDN25}. The eigenvalues of $T_\infty$ and $T_0$ are $\mathrm e^{2\pi i\alpha_j}$ and $\mathrm e^{-2\pi i\beta_j}$, respectively, where
\begin{equation}
\label{eq:C32-parameters}
\alpha=\left(0,0,0,0,\frac16,\frac56\right),\qquad
\beta=\left(\frac1{12},\frac14,\frac5{12},
\frac7{12},\frac34,\frac{11}{12}\right).
\end{equation}
Let $z$ be the standard coordinate on $\mathbb P^1$. Applying Fedorov's Hodge-number formula in the coordinate $z^{-1}$, the counts are
\[
\left(\#\{i:\alpha_i<\beta_k\}-k\right)_{k=1}^6=(3,3,2,1,0,0).
\]
The multiplicity of $p$ is $h^{p,3-p}$, so the Hodge numbers are $(2,1,1,2)$ \cite[Theorem~1]{Fedorov18}.
We describe the local filtered spaces used in the calculation. Fix a puncture $x$, and identify the flat fibres on the universal cover of a punctured disc with a fixed vector space $E_x$.
The local monodromy acts on this space by $T_x$. Write $T_x=T_{x,s}\exp N_x$, where $T_{x,s}$ is diagonalizable, $N_x$ is nilpotent, and they commute. Here $T_{x,s}$ has finite order. Choose a limiting Hodge filtration $F^\bullet E_x$, given by Schmid's theorem after a finite cover making the monodromy unipotent \cite{Schmid73}.
It is preserved by $T_{x,s}$, so every generalized eigenspace
\[
E_{x,\lambda}:=\ker\left((T_x-\lambda I)^6\right)
\]
inherits a filtration.

\begin{remark}
\label{rem:C32-hodge-indices}
We use the local invariants $\nu$ of \cite[Definition~2.14 and Remark~2.15]{Fedorov18} in Fedorov's normalization, whose Hodge indices are one higher than our weight-three indices. Thus $\nu^{p+1}_{a,\ell}$ is the dimension of the $p$-th primitive Hodge quotient in our grading. In particular, the index $4$ below corresponds to Hodge index $3$.
The parameter $a$ records the eigenvalue $\mathrm e^{2\pi ia}$, with parabolic weight $-a$ modulo $1$.
\end{remark}

At $0$, each space $E_{0,\mathrm e^{-2\pi i\beta_k}}$ is a line because the $\beta_k$ are distinct. With the indexing convention of Remark~\ref{rem:C32-hodge-indices}, we have
\[
\dim\operatorname{Gr}_F^pE_{0,\mathrm e^{-2\pi i\beta_k}}=\nu^{p+1}_{-\beta_k,0}.
\]
Applying \cite[Theorem~3]{Fedorov18} in the coordinate $z^{-1}$, this dimension is equal to one precisely when
\[
p=\#\{i:\alpha_i<\beta_k\}
 -\#\{i:\beta_i<\beta_k\}-1
 =\#\{i:\alpha_i<\beta_k\}-k,
\]
and zero otherwise. The last equality uses $\#\{i:\beta_i<\beta_k\}=k-1$. Explicitly,
\[
\begin{array}{c|cccccc}
k&1&2&3&4&5&6\\\hline
\beta_k&\frac1{12}&\frac14&\frac5{12}&\frac7{12}&\frac34&\frac{11}{12}\\[1mm]
\#\{i:\alpha_i<\beta_k\}&4&5&5&5&5&6\\
p&3&3&2&1&0&0.
\end{array}
\]
At infinity, let $a\in\{0,1/6,5/6\}$ and put $E=E_{\infty,\mathrm e^{2\pi ia}}$. The restriction $T_\infty|_E$ has a single Jordan block because $T_\infty=A$ is a companion matrix. Hence $E/N_\infty E$ is a line, and
\[
\dim\operatorname{Gr}_F^p(E/N_\infty E)=\nu^{p+1}_{a,\dim E-1}.
\]
Applying \cite[Theorem~3]{Fedorov18} in the coordinate $z^{-1}$, this dimension is equal to one precisely when
\[
p=\#\{i:\alpha_i\leqslant a\}-\#\{i:\beta_i\leqslant a\}-1,
\]
and zero otherwise. Explicitly,
\[
\begin{array}{c|ccc}
a&0&\frac16&\frac56\\\hline
\dim E&4&1&1\\
\#\{i:\alpha_i\leqslant a\}&4&5&6\\
\#\{i:\beta_i\leqslant a\}&0&1&5\\
p&3&3&0.
\end{array}
\]
Note that $\dim\operatorname{Gr}_F^pE_\infty=h^{p,3-p}$ for every $p$.
Since $T_{\infty,s}$ preserves $F^pE_\infty$, its eigenspace projections, which are polynomials in $T_{\infty,s}$, preserve it as well. Consequently,
\[
F^pE_\infty
=\bigoplus_{a\in\{0,1/6,5/6\}}
 F^pE_{\infty,\mathrm e^{2\pi ia}}.
\]
For $a=1/6$ and $5/6$, these eigenspaces are lines and $N_\infty$ vanishes on them, so the preceding computation gives their only nonzero graded quotients at $p=3$ and $p=0$, respectively. Taking graded quotients in this direct sum yields
\[
\left(\dim\operatorname{Gr}_F^pE_{\infty,1}\right)_{p=0}^3
=(2,1,1,2)-(0,0,0,1)-(1,0,0,0)
=(1,1,1,1).
\]
Parabolic weights are the negatives of monodromy exponents modulo $1$, giving
\[
\begin{array}{c|c|c}
&0&\infty\\\hline
V^{3,0}&\{\frac1{12},\frac14\}&\{0,\frac56\}\\[1mm]
V^{2,1}&\{\frac5{12}\}&\{0\}\\[1mm]
V^{1,2}&\{\frac7{12}\}&\{0\}\\[1mm]
V^{0,3}&\{\frac34,\frac{11}{12}\}&\{0,\frac16\}.
\end{array}
\]
Since $A$ and $B$ differ only in their last column, $T_1-I=A^{-1}(B-A)$ has rank one. Hence $T_1$ fixes a five-dimensional subspace pointwise. As $\det T_1=\det B/\det A=1$, its remaining eigenvalue is also $1$. Thus $T_1$ is unipotent, and its parabolic weights are zero.

We compute the ordinary degrees directly at these parameters. Let $j\colon X\hookrightarrow\mathbb P^1$, and put $\delta^p=\deg\left(\Xi_h(V^{p,3-p})^{\geqslant0}\right)$.
The dimension of $\ker(T_x-I)$ is zero at $0$, since $T_0$ has no eigenvalue $1$; five at $1$, since $T_1-I$ has rank one; and one at infinity, since the eigenvalue $1$ occurs in a single Jordan block.
The Euler--Poincar\'e formula \cite[proof of Theorem~1.1.2]{Katz96}, with $(j_*\mathbb{V})_x\cong\ker(T_x-I)$ \cite[\S1]{DettweilerWewers06}, gives, since $X$ is the sphere with three points removed and $\mathbb{V}$ has rank six,
\[
\chi(\mathbb P^1,j_*\mathbb{V})=(2-3)\operatorname{rank}\mathbb{V}
 +\sum_{x\in\{0,1,\infty\}}\dim\ker(T_x-I)=6(2-3)+0+5+1=0.
\]
Since $f$ and $g$ have no common root, the monodromy generated by $A$ and $B$ is irreducible by \cite[Proposition~3.3]{BeukersHeckman89}.
Both $\mathbb{V}$ and its dual $\mathbb{V}^\vee$ are therefore irreducible of rank six. A nonzero global flat section of either would span a monodromy-invariant line, contradicting irreducibility. Hence $H^0(X,\mathbb{V})=H^0(X,\mathbb{V}^\vee)=0$. Using the definition of $j_*$ for $H^0$ and Poincar\'e duality for $H^2$ (\cite[Remark~2.3.5]{DettweilerSabbah13}), we obtain
\[
H^0(\mathbb P^1,j_*\mathbb{V})=H^0(X,\mathbb{V})=0,\qquad H^2(\mathbb P^1,j_*\mathbb{V})\cong H^0(X,\mathbb{V}^\vee)^\vee=0.
\]
Together with $\chi(\mathbb P^1,j_*\mathbb{V})=0$, this gives $H^1(\mathbb P^1,j_*\mathbb{V})=0$.
Keeping the indexing convention of Remark~\ref{rem:C32-hodge-indices}, we write $\nu^q_{x,a,\ell}$ to indicate the puncture $x$.
For $a=0,1/6,5/6$, the table at infinity gives block sizes $4,1,1$, and the corresponding lines $E/N_\infty E$ have their only nonzero graded quotients at $p=3,3,0$, respectively.
In the preceding identity $\dim\operatorname{Gr}_F^p(E/N_\infty E)=\nu^{p+1}_{a,\dim E-1}$, these give $\ell=3,0,0$ and $q=p+1=4,4,1$. Hence
\[
\nu^4_{\infty,0,3}=\nu^4_{\infty,1/6,0}=\nu^1_{\infty,5/6,0}=1.
\]
And all other $\nu^q_{\infty,a,\ell}$ vanish. At $1$, the only Jordan block of size greater than one has size two. Put $E=E_{1,1}$ and $N=N_1$, so $N$ has rank one and $E/\ker N$ is a line. Choose a nonzero real vector $u$ with $\operatorname{im}N=\mathbb C u$. For the symplectic polarization $Q$, the identity $Q(Nv,w)+Q(v,Nw)=0$ gives $\ker N=u^{\perp_Q}$.
The nilpotent-orbit theorem \cite[Theorem~4.9]{Schmid73} gives $NF^pE\subset F^{p-1}E$ and makes $\mathrm e^{iyN}F^\bullet E$ a real weight-three Hodge filtration for sufficiently large real $y$.
Its polarization gives $Q(\mathrm e^{iyN}F^2E,\mathrm e^{iyN}F^2E)=0$. Since $\mathrm e^{iyN}$ preserves $Q$, this implies $Q(F^2E,F^2E)=0$. Moreover, $\dim F^2E=h^{3,0}+h^{2,1}=3=\frac12\dim E$, so $F^2E$ is Lagrangian.
If $NF^3E\neq0$, then $u\in F^2E$. Isotropy of $F^2E$ would give $F^3E\subset u^{\perp_Q}=\ker N$, a contradiction. Thus $NF^3E=0$.
If $NF^2E=0$, then $F^2E\subset u^{\perp_Q}$, and taking orthogonal complements gives $u\in F^2E$. But then $\mathrm e^{iyN}F^2E=F^2E$ contains the nonzero real vector $u$, contradicting
\[
\mathrm e^{iyN}F^2E\cap\overline{\mathrm e^{iyN}F^2E}=0.
\]
Therefore $NF^2E\neq0$. By definition, $F^p(E/\ker N)=(F^pE+\ker N)/\ker N$.
Since $NF^2E\neq0$, the subspace $F^2E$ is not contained in $\ker N$, so its image in $E/\ker N$ is nonzero. This quotient is one-dimensional, so that image is the whole quotient. On the other hand, $NF^3E=0$ means $F^3E\subset\ker N$, so its image in the quotient is zero.
Consequently,
\[
F^2(E/\ker N)=E/\ker N,\qquad F^3(E/\ker N)=0.
\]
The only nonzero graded quotient is at $p=2$. As above, $q=p+1=3$, and the block size two gives $\ell=1$. Hence
\[
\nu^3_{1,0,1}=1,
\]
and all other $\nu^q_{1,0,\ell}$ with $\ell>0$ vanish.
Apply \cite[Proposition~2.3.3]{DettweilerSabbah13} in the affine coordinate $z^{-1}$, whose point at infinity is $z=0$.
The correction term at the point at infinity in that formula is $-\sum_{\ell\geqslant0}\nu^p_{0,0,\ell}$ in our notation.
It vanishes because $T_0$ has no eigenvalue $1$. Expressing the remaining local terms using \cite[(7), (18)]{Fedorov18} gives
\[
0=\delta^{p-1}-\delta^p-h^{p,3-p}
 +\nu^p_{\infty,1/6,0}+\nu^p_{\infty,5/6,0}
 +\sum_{k=0}^{2}\nu^{p+1+k}_{\infty,0,3}
 +\nu^{p+1}_{1,0,1}.
\]
The local terms sum to $(0,2,2,1)$ for $p=0,1,2,3$. Starting with $\delta^{-1}=0$, the successive differences for $p=0,1,2,3$ are $(-2,1,1,-1)$. Consequently $(\delta^3,\delta^2,\delta^1,\delta^0)=(-1,0,-1,-2)$.
Adding the local weights yields
\begin{equation}
\operatorname{pardeg}V^{3,0}=\frac16,\quad
\operatorname{pardeg}V^{2,1}=\frac5{12},\quad
\operatorname{pardeg}V^{1,2}=-\frac5{12},\quad
\operatorname{pardeg}V^{0,3}=-\frac16,
\label{eq:C32-hodge-degrees}
\end{equation}
and therefore
\begin{equation}
\label{eq:C32-F2-degree}
\operatorname{pardeg}F^2=\frac7{12}.
\end{equation}

We pass to the sixfold cover $q(t)=t^6$. Put
\[
D_Y=\{t\in\mathbb C:t^6=1\}\cup\{\infty\},\qquad
Y=\mathbb P^1_t\setminus D_Y,\qquad Y^\circ=Y\setminus\{0\}.
\]
Then $q\colon Y^\circ\to X$ is an unramified cover of degree six.
Since $T_0^6=-I$, the pullback has monodromy $-I$ at $t=0$. Let $\mathbb L$ be the real rank-one local system on $Y^\circ$ with monodromy $-1$ at $0,\infty$ and $+1$ at each sixth root of unity, equipped with its unitary metric and Hodge type $(0,0)$. The twisted pullback $q^*\mathbb{V}\otimes\mathbb L$ has trivial monodromy at $0$, so it extends to a polarized real variation $\mathbb{V}_Y$ on $Y$ \cite[Theorem~4.9]{Schmid73}.
Write $V_Y^{p,3-p}$ and $F_Y^2$ for its Hodge bundles and filtration, and denote the component $V_Y^{2,1}\to V_Y^{1,2}\otimes\KK_Y$ of its Higgs field by $\Phi_{Y,2,1}$.

Let $\Gamma_Y$ and $\Gamma$ be the monodromy groups of $\mathbb{V}_Y$ and $\mathbb{V}$, respectively. The untwisted pullback image has index at most six in $\Gamma$ and contains $-I$. It equals $\Gamma_Y\cdot\{\pm I\}$, since the twist only changes signs. Thus $\Gamma_Y\subset\Gamma$ and $[\Gamma:\Gamma_Y]\leqslant12$.
In particular, $\Gamma_Y$ is integral and Zariski dense in $\mathrm{Sp}_6\mathbb R$, as is $\Gamma$ \cite{BDSS21,BDN25}, so $\mathbb{V}_Y$ is irreducible.

Tensoring by the unitary flat line preserves parabolic degrees, since that line has parabolic degree zero. Pullback by $q$ multiplies parabolic degrees by six. Indeed, the growth definition \eqref{eq:metric-extension} shows that a local ramification $q=t^e$ multiplies each parabolic weight by $e$, with the integer part incorporated into the underlying extension.
The smooth extension across $0$ contributes no parabolic weight. Consequently,
\[
\operatorname{pardeg}V_Y^{2,1}=\frac52,\qquad
\operatorname{pardeg}V_Y^{1,2}=-\frac52,\qquad
\operatorname{pardeg}F_Y^2=\frac72.
\]
The line bundles $V_Y^{2,1}$ and $V_Y^{1,2}$ have weight zero at each sixth root of unity because the monodromy there is unipotent. At infinity, $T_{\infty,s}^6=I$, so the semisimple part of the twisted monodromy $-T_\infty^6$ is $-I$. Since $-1=\mathrm e^{-2\pi i(1/2)}$, both line bundles have parabolic weight $1/2$ there. Their Hom line therefore has zero parabolic weights. Denote its metric extension to $\mathbb P^1_t$ by $\mathcal H$. The morphism $\Phi_{Y,2,1}$ defines a holomorphic section of $\mathcal H\otimes K_{\mathbb P^1_t}(D_Y)$, whose degree is
\[
\deg\left(\mathcal H\otimes K_{\mathbb P^1_t}(D_Y)\right)=-\frac52-\frac52+(-2+7)=0.
\]
The section $\Phi_{Y,2,1}$ is nonzero: otherwise $V_Y^{3,0}\oplus V_Y^{2,1}$ would be a Higgs subbundle of positive parabolic degree $7/2$, contradicting parabolic semistability. A nonzero holomorphic section of a degree-zero line bundle has no zeros, so $\Phi_{Y,2,1}$ is nowhere zero, including at the filled point $t=0$.

Because $V_Y^{2,1}$ and $V_Y^{1,2}$ have equal parabolic weights at each cusp, the nowhere-vanishing $\Phi_{Y,2,1}$ induces an isomorphism 
\[
\Phi_{Y,2,1}\colon(V_Y^{2,1})_\bullet\xrightarrow{\sim}
(V_Y^{1,2})_\bullet\otimes K_{\mathbb P^1_t}(D_Y).
\]
Equip $Y$ with its complete hyperbolic metric $g_Y$ of curvature $-4$. Theorem~\ref{thm:main-C} therefore applies to $\mathbb{V}_Y$. The real symplectic construction in the proof of Theorem~\ref{thm:main-C} gives a proper closed $\Gamma_Y$-invariant set $\Lambda$ of real Lagrangians transverse to $F_Y^2$ everywhere.

\begin{proposition}
\label{prop:C32-contact}
The C-32 monodromy group is thin. Let $\lambda_1(\mathbb{V})\geqslant\cdots\geqslant\lambda_6(\mathbb{V})$ be the Lyapunov exponents of $\mathbb{V}$ over the unit-speed geodesic flow of the complete curvature $-4$ metric $g_X$, with respect to normalized Liouville measure. Then
\begin{equation}\label{eq:C32-Lyapunov-equality}
\lambda_1(\mathbb{V})+\lambda_2(\mathbb{V})+\lambda_3(\mathbb{V})=\frac76.
\end{equation}
\end{proposition}

\begin{proof}
The proper closed $\Gamma_Y$-invariant set $\Lambda$ excludes the lattice property of $\Gamma_Y$ by Proposition~\ref{prop:lattice-minimality}. Since $[\Gamma:\Gamma_Y]<\infty$, the group $\Gamma$ is not a lattice either. Its integral symplectic structure places it in an arithmetic lattice, where it must have infinite index. Together with its Zariski density, this proves thinness.

For the equality, write $\lambda_1^{Y}(W_3)$ for the largest exponent of $\mathbb{V}_Y$ in the representation $W_3$, with the same convention on $(Y,g_Y)$. Since $Y$ is the sphere with seven points removed, $\deg K_{\mathbb P^1_t}(D_Y)=5$. Theorem~\ref{thm:main-C} gives
\[
\lambda_1^{Y}(W_3)=\frac{2\operatorname{pardeg}F_Y^2}{\deg K_{\mathbb P^1_t}(D_Y)}=\frac{2\cdot(7/2)}5=\frac75.
\]
It remains to compare the complete hyperbolic metrics on $Y$ and $Y^\circ$. The latter is $g_{Y^\circ}=q^*g_X$, since $q\colon Y^\circ\to X$ is unramified. The unitary twist preserves all transport singular values, and passage to this finite cover preserves Lyapunov exponents. Thus the exponents of $\mathbb{V}_Y|_{Y^\circ}$ for $g_{Y^\circ}$ are those of the original variation on $X$.

Write $g_Y=a\,g_{Y^\circ}$ on $Y^\circ$, and let $\omega_Y$ and $\omega_{Y^\circ}$ be their area forms.
In real dimension two,
\[
\Delta_{g_Y}=a^{-1}\Delta_{g_{Y^\circ}},\qquad \omega_Y=a\,\omega_{Y^\circ}.
\]
Start Brownian motion $Z_t$, with generator $\tfrac12\Delta_{g_{Y^\circ}}$, in normalized $g_{Y^\circ}$-area, and put $A_t=\int_0^t a(Z_s)\,\mathrm{d}s$.
The first identity shows that the time-changed path $\widehat Z_u=Z_{A^{-1}(u)}$ has generator $\tfrac12\Delta_{g_Y}$. By ergodicity \cite[\S1.1]{DanielDeroin19} and the second identity,
\[
\lim_{t\to\infty}\frac{A_t}{t}
=\frac{1}{\operatorname{Area}(Y^\circ)}
 \int_{Y^\circ}a\,\omega_{Y^\circ}
=\frac{\operatorname{Area}(Y)}{\operatorname{Area}(Y^\circ)}
=\frac{5\pi/2}{3\pi}=\frac56,
\]
where the areas follow from Gauss--Bonnet. In particular, $A_t\to\infty$, so $\widehat Z$ is defined for all time.
The two normalized area measures have the same null sets, so the almost-sure Lyapunov limits for $Y$ apply to $\widehat Z$. Parallel transport along $Z$ up to time $t$ equals that along $\widehat Z$ up to time $A_t$. Hence the Brownian exponents for $Y^\circ$ are $5/6$ times those for $Y$. Using the
Brownian--geodesic coincidence of Proposition~\ref{prop:lyap-coincidence}\romitem{1} on both surfaces, we obtain
\[
\lambda_1(\mathbb{V})+\lambda_2(\mathbb{V})+\lambda_3(\mathbb{V})
=\frac56\,\lambda_1^{Y}(W_3)
=\frac76.
\]
\end{proof}

\begin{remark}
\label{rem:C32-pingpong}
Riestenberg--Taha--Trettel--Wienhard \cite[Theorem~2]{RTTW26} prove C-32 thin by an exact ping--pong certificate and identify its projective monodromy group with $\mathbb Z*(\mathbb Z/6)$. The argument above obtains thinness from the invariant set of Lagrangians and establishes the Lyapunov equality \eqref{eq:C32-Lyapunov-equality}.
\end{remark}

\subsection{The fourteen hypergeometric Calabi--Yau systems}

We consider the fourteen hypergeometric Calabi--Yau systems in the Doran--Morgan list \cite[Corollary~2.3, \S2.2 \& Table~1]{DoranMorgan06}. Brav--Thomas proved thinness in seven cases \cite[Theorem~1.1 \& Table~1]{BravThomas14}. Arithmeticity in the remaining seven cases was proved by Singh--Venkataramana for three cases \cite[Corollaries~1.3 \&~1.4]{SinghVenkataramana14} and by Singh for the other four \cite[Theorem~1.1]{Singh15}.
Filip proved the Lyapunov equality for the seven thin cases \cite[Theorems~1 \&~6, Remark~2.2.10]{Filip21}. Combining this result with Theorem~\ref{thm:symplectic-contact} gives the following corollary.

\begin{corollary}\label{cor:fourteen-CY}
Let $\mathbb{V}$ be one of the fourteen classical rank-four hypergeometric Calabi--Yau local systems on $X=\mathbb P^1\setminus\{0,1,\infty\}$, with monodromy group $\Gamma$. Then
\begin{equation}\label{eq:CY-fourteen-eq}
\lambda_1(\mathbb{V})+\lambda_2(\mathbb{V})=2\operatorname{pardeg}F^2\quad\Longleftrightarrow\quad\Gamma\text{ is thin}.
\end{equation}
In particular, the seven arithmetic cases satisfy
\begin{equation}\label{eq:CY-fourteen-strict}
\lambda_1(\mathbb{V})+\lambda_2(\mathbb{V})>2\operatorname{pardeg}F^2.
\end{equation}
\end{corollary}

\begin{proof}
These systems are irreducible polarized real weight-three variations with Hodge numbers $(1,1,1,1)$ and Zariski-dense integral symplectic monodromy.
Filip's equality for the thin cases also holds after removing the additional points on the base \cite[p.~14]{Filip21}, so it applies to $X$. His curvature-$-1$ convention is converted to ours by $\lambda_i^{(-4)}=2\lambda_i^{(-1)}$, giving
\eqref{eq:CY-fourteen-eq} in these seven cases.

In each arithmetic case, $\Gamma$ is a lattice. Equality would contradict the Hodge-filtration assertion of Theorem~\ref{thm:symplectic-contact}, since there are four nonzero Hodge summands.
Since $\deg\KK_{\mathbb P^1}(0+1+\infty)=1$, the general lower bound \cite[Theorem~4.1]{EKMZ18} gives $\lambda_1(\mathbb{V})+\lambda_2(\mathbb{V})\geqslant2\operatorname{pardeg}F^2$, and hence strictness.
\end{proof}

This recovers the equality and strictness dichotomy observed numerically by Kontsevich and recorded in \cite[\S6 \& Table~1]{EKMZ18}.

\section{Comparison between Lyapunov and average Cartan spectra}
\label{sec:Lyapunov-Cartan}

We compare the Lyapunov spectrum of Section~\ref{sec:setup} with the average Cartan spectrum determined by the harmonic metric.

Let $\rho\colon\pi_1(X)\to G_{\mathbb R}$ be the monodromy of a $G_\mathbb{R}$-local system, which is semisimple with non-expanding cusp monodromy. As in Section~\ref{sec:hermitian-transversality}, $G_{\mathbb R}$ is a connected semisimple real Lie group with finite center and without compact factors.
Fix the compatible tame harmonic $K_{\mathbb R}$-reduction $h$ adapted to the Deligne extensions, and let $u\colon\widetilde X\to G_{\mathbb R}/K_{\mathbb R}$ be its $\rho$-equivariant harmonic map.
Write $D_h$ for its connection and $\mathfrak{g}_{\mathbb R}=\mathfrak{k}_{\mathbb{R}}\oplus\mathfrak{p}_{\mathbb{R}}$. Thus
\[
\nabla=D_h+\Phi+\Phi^{*_h}=D_h+\mathrm{d}u,
\]
with $\mathrm{d}u=\Phi+\Phi^{*_h}$ viewed as a $\mathfrak{p}_{\mathbb{R}}$-valued one-form. Proposition~\ref{prop:adapted-metric}\romitem{1}, applied to the adjoint harmonic bundle, implies that $\mathrm{d}u$ is bounded in hyperbolic norm.

For the group-level factorization below, let $\eta\colon G_{\mathbb R}\to G_{\mathbb C}$ be the universal complexification map of Section~\ref{sec:hermitian-transversality}, and set
\[
G_{\mathbb R}^{\mathrm{lin}}:=\eta(G_{\mathbb R}),\qquad
K_{\mathbb R}^{\mathrm{lin}}:=\eta(K_{\mathbb R}),\qquad
\rho^{\mathrm{lin}}:=\eta\circ\rho.
\]
Its finite central kernel lies in $K_{\mathbb R}$, so $\eta$ identifies $G_{\mathbb R}/K_{\mathbb R}$ with $G_{\mathbb R}^{\mathrm{lin}}/K_{\mathbb R}^{\mathrm{lin}}$ and preserves the Lie algebra, Cartan projection, and Lyapunov spectrum.
Every finite-dimensional real representation $R$ of $G_{\mathbb R}$ descends to $R^{\mathrm{lin}}$ on $G_{\mathbb R}^{\mathrm{lin}}$, by the universal property recalled in Section~\ref{sec:hermitian-transversality}.
This allows the global factorization to be stated without assuming that $G_{\mathbb R}$ is linear.

\subsection{Finite-time Cartan comparison}

We use the Cartan data, dominance order, Cartan projection $\kappa$, chamber-valued Lyapunov spectrum $\boldsymbol{\lambda}$, and triangle inequality \eqref{eq:Cartan-triangle} from Section~\ref{sec:setup}, with the maximal compact subgroup $K_{\mathbb R}$ fixed above.

Let $c\colon\mathbb R\to X$ be a unit-speed geodesic. In a $D_h$-parallel trivialization along $c$, the generator of flat transport is $-\mathrm{d}u(\dot c(t))$. For every $v\in T^1X$, define
\[
H(v):=\kappa\!\left(\exp(-\mathrm{d}u(v))\right)
\in\overline{\mathfrak a^+}.
\]
Here $\mathrm{d}u(v)$ is written in any local frame of the $K_{\mathbb R}$-reduction $h$. The definition is independent of this choice because $\kappa$ is invariant under $K_{\mathbb R}$-conjugation. Equivalently, $H(v)$ is the unique dominant element in the $K_{\mathbb R}$-orbit of $-\mathrm{d}u(v)$. The boundedness of $\mathrm{d}u$ gives a uniform bound for $H(v)$. Set 
\[\mathbf{C}:=\int_{T^1X}H(v)\,\mathrm{d}\mu_L(v)\in\overline{\mathfrak a^+},\] 
which we call the \emph{average Cartan spectrum}, and put $C_\alpha:=\omega_\alpha(\mathbf{C})$.

Let $U_{s,t}(v)\in G_{\mathbb R}$ be flat transport, written in the $D_h$-parallel trivialization, along the geodesic segment $g_rv$, $s\leqslant r\leqslant t$, and put $U_T(v)=U_{0,T}(v)$.

\begin{lemma}
\label{lem:finite-Cartan-inequality}
For every $v\in T^1X$ and every $s<t$,
\begin{equation}
	\kappa(U_{s,t}(v))\preceq\int_s^tH(g_rv)\,\mathrm{d}r.
\label{eq:finite-Cartan-inequality}
\end{equation}
\end{lemma}

\begin{proof}
Choose a partition $s=t_0<\cdots<t_N=t$. The solution of the transport ODE is the limit, as the mesh tends to zero, of the ordered products $\prod_{j=N-1}^{0}\exp(-(t_{j+1}-t_j)\,\mathrm{d}u(g_{t_j}v))$. Since $\kappa(\exp(-r\,\mathrm{d}u(g_{t_j}v)))=rH(g_{t_j}v)$ for $r\geqslant0$, iterating \eqref{eq:Cartan-triangle} bounds the Cartan projection of each product by the corresponding Riemann sum. Passing to the limit gives \eqref{eq:finite-Cartan-inequality}.
\end{proof}

Define the finite-time Cartan gap
\begin{equation}\label{eq:finite-Cartan-gap}
\mathfrak G_{s,t}(v):=\int_s^tH(g_rv)\,\mathrm{d}r-\kappa(U_{s,t}(v))\succeq0
\end{equation}
and its fundamental-weight coordinates $d_{\alpha;s,t}(v):=\omega_\alpha(\mathfrak G_{s,t}(v))\geqslant0$. The cocycle identity $U_{s,u}(v)=U_{t,u}(v)U_{s,t}(v)$ and \eqref{eq:Cartan-triangle} give
\begin{equation}\label{eq:Cartan-gap-superadditive}
\begin{aligned}
\mathfrak G_{s,u}(v)
&=\int_s^uH(g_rv)\,\mathrm{d}r-\kappa(U_{s,u}(v))\\
&\succeq\int_s^tH(g_rv)\,\mathrm{d}r+\int_t^uH(g_rv)\,\mathrm{d}r-\kappa(U_{s,t}(v))-\kappa(U_{t,u}(v))\\
&=\mathfrak G_{s,t}(v)+\mathfrak G_{t,u}(v).
\end{aligned}
\end{equation}
Thus every $d_{\alpha;s,t}$ is superadditive.

\begin{proposition}\label{prop:upper-LA}
	One has $\boldsymbol{\lambda}\preceq\mathbf{C}$, equivalently $\lambda_\alpha\leqslant C_\alpha$ for every $\alpha\in\Delta$.
\end{proposition}

\begin{proof}
For every $\alpha\in\Delta$, applying $\omega_\alpha$ to \eqref{eq:finite-Cartan-inequality} with $s=0$ and $t=T$ gives
\[\omega_\alpha(\kappa(U_T(v)))\leqslant\int_0^T\omega_\alpha(H(g_rv))\,\mathrm{d}r.\]
Hence, for almost every $v$,
\[\frac1T\omega_\alpha(\kappa(U_T(v)))\leqslant\frac1T\int_0^T\omega_\alpha(H(g_rv))\,\mathrm{d}r.\]
Letting $T\to\infty$, the definition of the Lyapunov spectrum and Birkhoff's theorem give
\[\lambda_\alpha\leqslant\int_{T^1X}\omega_\alpha(H(v))\,\mathrm{d}\mu_L(v)=\omega_\alpha(\mathbf{C})=C_\alpha.\]
Since this holds for every fundamental restricted weight, $\boldsymbol{\lambda}\preceq\mathbf{C}$.
\end{proof}

Equality in one asymptotic Cartan coordinate forces equality in that coordinate on every finite segment.

\begin{theorem}\label{thm:upper-asymptotic-finite}
	Fix $\alpha\in\Delta$. If $\lambda_\alpha=C_\alpha$, then $d_{\alpha;s,t}(v)=0$ for every $v\in T^1X$ and every finite interval $[s,t]$. Consequently, if $\boldsymbol{\lambda}=\mathbf{C}$, then
	\begin{equation}\label{eq:exact-Cartan-additivity}
		\kappa(U_{s,t}(v))=\int_s^tH(g_rv)\,\mathrm{d}r
	\end{equation}
	for every finite geodesic segment.
\end{theorem}

\begin{proof}
	It is enough to treat $[0,T]$. Fix $T>0$. By \eqref{eq:Cartan-gap-superadditive},
	\[d_{\alpha;0,nT}(v)\geqslant\sum_{j=0}^{n-1}d_{\alpha;0,T}(g_{jT}v).\]
	For fixed $T$, the function $d_{\alpha;0,T}$ is continuous, nonnegative, and bounded. Suppose that $\int_{T^1X}d_{\alpha;0,T}\,\mathrm{d}\mu_L>0$. Since the geodesic flow on a finite-area hyperbolic surface is mixing, the time-$T$ map $g_T$ is ergodic. Birkhoff's theorem therefore gives, for almost every $v$,
	\[\frac1n\sum_{j=0}^{n-1}d_{\alpha;0,T}(g_{jT}v)\longrightarrow\int_{T^1X}d_{\alpha;0,T}\,\mathrm{d}\mu_L>0.\]
	Together with the preceding superadditivity inequality, this yields
	\[\liminf_{n\to\infty}\frac{1}{nT}d_{\alpha;0,nT}(v)\geqslant\frac1T\int_{T^1X}d_{\alpha;0,T}\,\mathrm{d}\mu_L>0\]
	for almost every $v$. On the other hand,
	\[\frac{1}{nT}d_{\alpha;0,nT}(v)=\frac1{nT}\int_0^{nT}\omega_\alpha(H(g_rv))\,\mathrm{d}r-\frac1{nT}\omega_\alpha(\kappa(U_{nT}(v)))\longrightarrow C_\alpha-\lambda_\alpha=0\]
	for almost every $v$, a contradiction. Hence the integral of $d_{\alpha;0,T}$ is zero. Nonnegativity gives almost-everywhere vanishing, and continuity together with full support of $\mu_L$ gives vanishing everywhere. Shifting the orbit gives arbitrary $[s,t]$.
	
	If $\boldsymbol{\lambda}=\mathbf{C}$, this holds for every simple root. Since the fundamental weights form a basis of $\mathfrak a^*$, the Cartan gap $\mathfrak G_{s,t}(v)$ itself vanishes, giving \eqref{eq:exact-Cartan-additivity}.
\end{proof}

\subsection{Rigidity of the equality}

We analyze equality in \eqref{eq:finite-Cartan-inequality} using the Iwasawa cocycle. Let $G_{\mathbb R}=K_{\mathbb R}AN$ be the Iwasawa decomposition determined by the positive restricted roots fixed in Section~\ref{sec:setup}, with $A=\exp(\mathfrak a)$. Let $M:=Z_{K_{\mathbb R}}(\mathfrak a)$ and $P:=MAN$.
Write $W_{\mathrm{res}}:=N_{K_{\mathbb R}}(\mathfrak a)/M$ for the restricted Weyl group. Then $P$ is the corresponding minimal parabolic subgroup and the full real flag variety is
\[\mathcal F:=G_{\mathbb R}/P\cong K_{\mathbb R}/M.\]
Indeed, $G_{\mathbb R}=K_{\mathbb R}P$ and $K_{\mathbb R}\cap P=M$. Thus every flag may be written $\xi=kP$, or equivalently $kM$, with $k\in K_{\mathbb R}$ unique modulo $M$.

For $g\in G_{\mathbb R}$ and $\xi=kP$, write the Iwasawa decomposition
\[gk=k'\exp(\sigma(g,\xi))n,\qquad k'\in K_{\mathbb R},\quad \sigma(g,\xi)\in\mathfrak a,\quad n\in N.\]
The $\mathfrak a$-component is unchanged if $k$ is replaced by another representative of $\xi$, so $\sigma(g,\xi)$ is well defined. It is the Iwasawa cocycle and satisfies
\[\sigma(g_2g_1,\xi)=\sigma(g_2,g_1\xi)+\sigma(g_1,\xi).\]
Its infinitesimal form is the map
\[q(X,\xi):=\left.\frac{\mathrm{d}}{\mathrm{d}t}\right|_{t=0}\sigma(\exp(tX),\xi),\qquad X\in\mathfrak{p}_{\mathbb{R}},\]
where $\mathfrak{g}_{\mathbb{R}}=\mathfrak{k}_{\mathbb{R}}\oplus \mathfrak{p}_{\mathbb{R}}$ is the Cartan decomposition.
If $\mathfrak{g}_{\mathbb R}=\mathfrak{k}_{\mathbb R}\oplus\mathfrak a\oplus\mathfrak n$ is the Lie-algebra Iwasawa decomposition and $\operatorname{pr}^{\mathrm{Iw}}_{\mathfrak a}$ denotes projection onto $\mathfrak a$, then
\begin{equation}\label{eq:infinitesimal-Iwasawa-cocycle}
q(X,\xi)=\operatorname{pr}^{\mathrm{Iw}}_{\mathfrak a}\left(\operatorname{Ad}(k^{-1})X\right).
\end{equation}
The same flag determines the closed Weyl chamber
\[
\overline{\mathfrak a^+_\xi}:=\operatorname{Ad}(k)\overline{\mathfrak a^+}\subset\mathfrak{p}_{\mathbb{R}}.
\]
Both \eqref{eq:infinitesimal-Iwasawa-cocycle} and the chamber are independent of the representative $k$, since $M$ centralizes $\mathfrak a$ and preserves the positive restricted-root spaces.

\begin{lemma}
\label{lem:upper-CI-positive}
Let $X\in\mathfrak{p}_{\mathbb{R}}$, and let $H_X\in\overline{\mathfrak a^+}$ be the dominant element in its $K_{\mathbb R}$-orbit. Then
\[q(X,\xi)\in\operatorname{conv}(W_{\mathrm{res}}\!\cdot\!H_X),
\qquad
q(X,\xi)\preceq H_X.\]
Moreover,
\begin{equation}\label{eq:CI-equality-chamber}
q(X,\xi)=H_X\quad\Longleftrightarrow\quad X\in\overline{\mathfrak a^+_\xi}.
\end{equation}
\end{lemma}

\begin{proof}
Choose $k\in K_{\mathbb R}$ representing the flag $\xi$, so that $\xi=kP$ under $\mathcal F=G_{\mathbb R}/P\cong K_{\mathbb R}/M$. Then
\[q(X,\xi)=\operatorname{pr}^{\mathrm{Iw}}_{\mathfrak a}\left(\operatorname{Ad}(k^{-1})X\right),
\qquad
\overline{\mathfrak a^+_\xi}
=\operatorname{Ad}(k)\overline{\mathfrak a^+}.
\]
Since $\operatorname{Ad}(k^{-1})X\in K_{\mathbb R}\!\cdot H_X$, Kostant's linear convexity theorem gives
\[
q(X,\xi)\in\operatorname{conv}(W_{\mathrm{res}}\!\cdot\!H_X),
\]
see \cite[Ch.~\upRoman{4}, Thm.~10.2]{Helgason84}. Hence $q(X,\xi)\preceq H_X$.

For the equality case, the Iwasawa $\mathfrak a$-component on $\mathfrak{p}_{\mathbb{R}}$ is the orthogonal projection onto $\mathfrak a$ for the Cartan inner product. Therefore
\[
\|q(X,\xi)\|
\leqslant
\|\operatorname{Ad}(k^{-1})X\|
=
\|H_X\|.
\]
If $q(X,\xi)=H_X$, equality of norms forces $\operatorname{Ad}(k^{-1})X\in\mathfrak a$. Its $\mathfrak a$-projection is $H_X$, so $\operatorname{Ad}(k^{-1})X=H_X\in\overline{\mathfrak a^+}$. Thus $X\in\operatorname{Ad}(k)\overline{\mathfrak a^+}=\overline{\mathfrak a^+_\xi}$. The converse follows from the same two formulas above.
\end{proof}

If $g=k_1\exp(\kappa(g))k_2$ is a Cartan decomposition, then the right flag $\xi_g:=k_2^{-1}P\in\mathcal F$ satisfies 
\begin{equation}\label{eq:Iwasawa-maximizing-flag}
\sigma(g,\xi_g)=\kappa(g).
\end{equation}

\begin{proposition}
\label{prop:finite-segment-Cartan-equality}
Fix $v\in T^1X$ and a finite interval $[s,t]\subset\mathbb R$. Suppose that
\begin{equation}\label{eq:finite-segment-Cartan-equality}
\kappa(U_{s,t}(v))
=
\int_s^tH(g_rv)\,\mathrm{d}r.
\end{equation}
Then there is a full flag $\xi\in\mathcal F$, independent of $r$, such that
\[-\mathrm{d}u(g_rv)\in\overline{\mathfrak a^+_\xi}
\qquad\text{for every }r\in[s,t].\]
Conversely, if there is a full flag $\xi\in\mathcal F$ for which this inclusion holds for every $r\in[s,t]$, then \eqref{eq:finite-segment-Cartan-equality} holds.
\end{proposition}

\begin{proof}
After shifting the parameter, it is enough to treat the interval $[0,T]$. Choose a Cartan decomposition of $U_T(v)$ and let $\xi_0\in\mathcal F$ be its right flag as in \eqref{eq:Iwasawa-maximizing-flag}. Write $U_r:=U_{0,r}(v)$, put $\xi_r:=U_r\xi_0$, and set $X_r:=-\mathrm{d}u(g_rv)$. In the $D_h$-parallel trivialization, the transport equation gives
\[\dot U_rU_r^{-1}=X_r.\]
For small $h$, the cocycle identity gives
\[\sigma(U_{r+h},\xi_0)-\sigma(U_r,\xi_0)=\sigma(U_{r+h}U_r^{-1},U_r\xi_0)=\sigma(U_{r+h}U_r^{-1},\xi_r).\]
Moreover, $U_{r+h}U_r^{-1}=\exp(hX_r+o(h))$. By the definition of the infinitesimal Iwasawa cocycle $q$, it follows that
\[
\frac{\mathrm{d}}{\mathrm{d}r}\sigma(U_r,\xi_0)=q(X_r,\xi_r)=q(-\mathrm{d}u(g_rv),\xi_r).
\]
Since $U_0=e$ and $\sigma(e,\xi_0)=0$, integration from $0$ to $T$ gives
\[\sigma(U_T(v),\xi_0)=\int_0^Tq(-\mathrm{d}u(g_rv),\xi_r)\,\mathrm{d}r.
\]
The left-hand side is $\kappa(U_T(v))$. Under the hypothesis, therefore,
\[\int_0^T\left(H(g_rv)-q(-\mathrm{d}u(g_rv),\xi_r)\right)\,\mathrm{d}r=0.\]
By Lemma~\ref{lem:upper-CI-positive}, every fundamental-weight coordinate of the integrand is nonnegative. Hence all of them vanish for every $r$, by continuity, and $q(-\mathrm{d}u(g_rv),\xi_r)=H(g_rv)$. The equality case \eqref{eq:CI-equality-chamber} gives $-\mathrm{d}u(g_rv)\in\overline{\mathfrak a^+_{\xi_r}}$. Such an element fixes the flag $\xi_r$ under the infinitesimal $G_{\mathbb R}$-action on the full flag variety. Since $\xi_r'=-\mathrm{d}u(g_rv)\cdot\xi_r$, we obtain $\xi_r'=0$. Thus $\xi_r\equiv\xi_0$, and all $-\mathrm{d}u(g_rv)$ lie in the same closed chamber.

Conversely, suppose that $-\mathrm{d}u(g_rv)\in \overline{\mathfrak a^+_\xi}$ for every $r\in[0,T]$. These elements commute in the parallel trivialization. Hence $U_{0,T}(v)=\exp(\int_0^T-\mathrm{d}u(g_rv)\,\mathrm{d}r)$, and the integral lies in the same closed chamber. Its Cartan projection is therefore $\int_0^TH(g_rv)\,\mathrm{d}r$.
\end{proof}

\begin{corollary}\label{cor:orbitwise-common-flat}
Assume $\boldsymbol{\lambda}=\mathbf{C}$. Along every complete geodesic there is a single closed Weyl chamber containing all $-\mathrm{d}u(g_tv)$. In particular, $[\mathrm{d}u(g_sv),\mathrm{d}u(g_tv)]=0$ for all $s,t\in\mathbb R$.
\end{corollary}

\begin{proof}
Theorem~\ref{thm:upper-asymptotic-finite} and Proposition~\ref{prop:finite-segment-Cartan-equality} give a full flag with the required chamber property on every interval $[-N,N]$. The full flag variety is compact.
Choose a convergent subsequence of these flags. For each fixed $T$, the condition $-\mathrm{d}u(g_tv)\in\overline{\mathfrak a^+_\xi}$ for $|t|\leqslant T$ is closed in $\xi$, so the limiting flag has the same chamber property on every compact interval,
and hence on the whole orbit.
\end{proof}

\begin{theorem}\label{thm:complete-Cartan-equality}
	For $\rho$ and its adapted tame harmonic reduction as fixed above, the following are equivalent:
\begin{enumerate}[label=\textup{(\roman*)}]
\item $\boldsymbol{\lambda}=\mathbf{C}$;

\item for every $v\in T^1X$ and every finite interval $[s,t]$, $\kappa(U_{s,t}(v))=\int_s^tH(g_rv)\,\mathrm{d}r$;

\item for every $v\in T^1X$, there is one fixed closed Weyl chamber containing $-\mathrm{d}u(g_tv)$ for every $t\in\mathbb R$.
\end{enumerate}
\end{theorem}

\begin{proof}
The implication \romitem{1}$\Rightarrow$\romitem{2} is Theorem~\ref{thm:upper-asymptotic-finite}. The implication \romitem{2}$\Rightarrow$\romitem{3} follows from Proposition~\ref{prop:finite-segment-Cartan-equality} and the compactness argument of Corollary~\ref{cor:orbitwise-common-flat}. Conversely, under \romitem{3} all $-\mathrm{d}u(g_tv)$ lie in one fixed closed chamber and commute, so the transport is the exponential of their integral and \romitem{2} follows.
Dividing \romitem{2} by $t-s$ and applying Birkhoff together with the definition of $\boldsymbol{\lambda}$ gives \romitem{1}.
\end{proof}

\subsection{Nilpotent-Higgs rigidity}

Write the harmonic-bundle decomposition as
\[\nabla=D_h+\Phi+\Phi^{*_h},
\qquad
\mathrm{d}u=\Phi+\Phi^{*_h}.\]
We use the same symbol $D_h$ for the connection induced on the Higgs tensor by the Chern connection on the adjoint bundle and the Levi--Civita connection $\nabla^{\mathrm{LC}}$ of the hyperbolic metric. Thus
\[(D_h\Phi)(X,Y)=D_{h,X}(\Phi(Y))-\Phi(\nabla^{\mathrm{LC}}_X Y).\]
Assume in this subsection that $\Phi$ is pointwise nilpotent, meaning that $\Phi(\xi)\in\mathfrak{g}_{\mathbb C}$ is nilpotent for every $\xi\in T^{1,0}X$. For a variation of Hodge structure, nilpotency is automatic.

\begin{theorem}\label{thm:upper-parallel-Higgs}
If $\boldsymbol{\lambda}=\mathbf{C}$, then $D_h\Phi=0$.
\end{theorem}

\begin{proof}
By Corollary~\ref{cor:orbitwise-common-flat}, complete upper equality implies that along every geodesic the elements $-\mathrm{d}u(g_tv)$ commute in a $D_h$-parallel trivialization.
Differentiating at a point, with the tangent vector extended Levi--Civita parallel there, gives
\[
[\mathrm{d}u(v),D_{h,v}(\mathrm{d}u(v))]=0.
\]

We first record the orbit-tangency consequence of nilpotency. On a coordinate disk $U$, choose a holomorphic coordinate $z$ and a holomorphic trivialization of the complexified adjoint bundle, and put $e=\Phi(\partial_z)$. Then $e\colon U\to\mathfrak{g}_{\mathbb C}$ is holomorphic and takes values in the nilpotent cone. A complex semisimple Lie algebra has only finitely many nilpotent adjoint orbits. For each such orbit $\mathcal O$, the set $e^{-1}(\mathcal O)$ is locally closed analytic in $U$. Since a proper analytic subset of a disk is discrete and finitely many orbit preimages cover $U$, there is a dense open subset $U^\circ\subset U$ and one nilpotent orbit $\mathcal O$ such that $e(U^\circ)\subset\mathcal O$.

For $z\in U^\circ$, $\mathrm{d} e(T_zU^\circ)\subset T_{e(z)}\mathcal O$,
and hence
\[\partial_z e\in T_e\mathcal O=T_e(G_{\mathbb C}\cdot e)=[\mathfrak{g}_{\mathbb C},e].\]
If $A_z$ is the $(1,0)$-connection matrix of $D_h$ and $\nabla^{\mathrm{LC}}_{\partial_z}\partial_z=\Gamma^z_{zz}\partial_z$, then
\[(D_h\Phi)(\partial_z,\partial_z)=\partial_z e+[A_z,e]-\Gamma^z_{zz}e.\]
If $e=0$, then $e\in[\mathfrak{g}_{\mathbb C},e]$ trivially; if $e\neq0$, Jacobson--Morozov gives an element $h$ with $[h,e]=2e$.
Thus in either case $e\in[\mathfrak{g}_{\mathbb C},e]$. Therefore, on $U^\circ$,
\begin{equation}\label{eq:nilpotent-orbit-tangency}
(D_h\Phi)(\partial_z,\partial_z)\in[\mathfrak{g}_{\mathbb C},e]=\operatorname{im}(\operatorname{ad}e).
\end{equation}
Under a holomorphic change of coordinate $w$, both $(D_h\Phi)(\partial_w,\partial_w)$ and $\Phi(\partial_w)$ are obtained from the corresponding $z$-expressions by multiplication by nonzero scalars. In particular, $\operatorname{im}(\operatorname{ad}\Phi(\partial_w))=\operatorname{im}(\operatorname{ad}\Phi(\partial_z))$.
Hence \eqref{eq:nilpotent-orbit-tangency} is coordinate-independent, and it holds on a dense open subset $X^\circ\subset X$.

Now fix $x\in X^\circ$. Choose a holomorphic normal coordinate $z$ for the hyperbolic metric at $x$, so that
\[\nabla^{\mathrm{LC}}_{\partial_z}\partial_z\big|_x=0.\]
Set
\[e=\Phi(\partial_z)|_x,\qquad f=e^{*_h},\qquad\eta=(D_h\Phi)(\partial_z,\partial_z)|_x.\]
By the choice of coordinate, $\eta=D_{h,\partial_z}(\Phi(\partial_z))|_x$. The holomorphicity of Higgs field give, at $x$,
\[
D_{h,\partial_{\bar z}}e=0,
\qquad
D_{h,\partial_z}f=0,
\qquad
D_{h,\partial_{\bar z}}f=\eta^{*_h}.
\]
For $|\zeta|=1$, put $v_\zeta=\zeta\partial_z+\zeta^{-1}\partial_{\bar z}$, extended Levi--Civita parallel at $x$. Up to the common positive normalization of this real direction,
\[\mathrm{d}u(v_\zeta)=\zeta e+\zeta^{-1}f,\qquad D_{h,v_\zeta}(\mathrm{d}u(v_\zeta))=\zeta^2\eta+\zeta^{-2}\eta^{*_h}.\]
Substitution into the common-flat relation gives
\[
0=
\zeta^3[e,\eta]
+\zeta[f,\eta]
+\zeta^{-1}[e,\eta^{*_h}]
+\zeta^{-3}[f,\eta^{*_h}].
\]
Hence
\begin{equation}
[f,\eta]=0
\label{eq:nilpotent-Fourier}
\end{equation}
because $\zeta$ is arbitrary.

By \eqref{eq:nilpotent-orbit-tangency}, $\eta\in\operatorname{im}(\operatorname{ad}e)$. Since $f=e^{*_h}$, the adjoint relation between $\operatorname{ad}e$ and $\operatorname{ad}f$ gives
\[\ker(\operatorname{ad}f)=\left(\operatorname{im}(\operatorname{ad}e)\right)^\perp.\]
Together with \eqref{eq:nilpotent-Fourier}, this yields
\[\eta\in\operatorname{im}(\operatorname{ad}e)
\cap
\left(\operatorname{im}(\operatorname{ad}e)\right)^\perp=\{0\}.\]
Thus
\[(D_h\Phi)(\partial_z,\partial_z)|_x=0\] 
for every $x\in X^\circ$. Since $D_h\Phi$ is a tensor, if $w$ is any other holomorphic coordinate, then at $x$
\[
(D_h\Phi)(\partial_w,\partial_w)
=
\left(\frac{\partial z}{\partial w}\right)^2
(D_h\Phi)(\partial_z,\partial_z).
\]
Thus the vanishing obtained in the normal coordinate is intrinsic. By continuity, the $(2,0)$-component of $D_h\Phi$ vanishes on all of $X$. Higgs holomorphicity gives the $(0,1)$-covariant derivative of $\Phi$ equal to zero, and therefore $D_h\Phi=0$.
\end{proof}

To summarize, we have the following theorem.

\begin{theorem}
\label{thm:main-D}
In the standing setup of this section, with the adapted tame harmonic reduction $h$, the following conditions are equivalent:
\begin{enumerate}[label=\textup{(\roman*)}]
\item $\boldsymbol{\lambda}=\mathbf{C}$;

\item for every finite geodesic segment, $\kappa(U_{s,t}(v))=\int_s^tH(g_rv)\,\mathrm{d}r$;

\item along every complete geodesic all $-\mathrm{d}u(g_tv)$ lie in one fixed closed Weyl chamber.
\end{enumerate}
If, in addition, the Higgs field $\Phi$ is pointwise nilpotent, these conditions are also equivalent to the following:
\begin{enumerate}[label=\textup{(\roman*)},resume]
\item $D_h\Phi=0$;

\item the harmonic map
$u\colon\widetilde X\to G_{\mathbb R}/K_{\mathbb R}$ is totally geodesic;

\item either $u$ is constant, or there are a Lie-group homomorphism $\iota\colon\mathrm{SL}_2\mathbb R\to G_{\mathbb R}^{\mathrm{lin}}$ whose orbit map $\mathbb H\to G_{\mathbb R}/K_{\mathbb R}$ through some point is a totally geodesic embedding, and a representation $\rho_c$ into the compact pointwise stabilizer of that disk in $G_{\mathbb R}^{\mathrm{lin}}$, such that
\begin{equation}\label{eq:linear-disk-compact}
\rho^{\mathrm{lin}}(\gamma)=\iota(j(\gamma))\rho_c(\gamma),
\end{equation}
where $j\colon\pi_1(X)\to\mathrm{SL}_2\mathbb R$ is a lift of the uniformizing representation. The pointwise stabilizer centralizes $\iota(\mathrm{SL}_2\mathbb R)$.
\end{enumerate}
\end{theorem}

\begin{proof}
The equivalence of \romitem{1}--\romitem{3} is Theorem~\ref{thm:complete-Cartan-equality}. Assume now that $\Phi$ is pointwise nilpotent. Theorem~\ref{thm:upper-parallel-Higgs} gives \romitem{1}$\Rightarrow$\romitem{4}.

We next identify \romitem{4} and \romitem{5}. Since $\mathrm{d}u=\Phi+\Phi^{*_h}$, the covariant derivative of $\mathrm{d}u$, formed using $D_h$ on the target bundle and the Levi--Civita connection on the source, is the second fundamental form of $u$. Its type decomposition, together with Higgs holomorphicity, shows that it vanishes exactly when $D_h\Phi=0$. Thus \romitem{4} and \romitem{5} are equivalent. Moreover, if these conditions hold and $c(t)$ is a geodesic, then $\mathrm{d}u(\dot c(t))$ is constant in a $D_h$-parallel frame. Hence
\[U_T=\exp\left(-T\,\mathrm{d}u(\dot c(0))\right),\]
so equality holds in \eqref{eq:finite-Cartan-inequality} for every finite geodesic segment. This also gives \romitem{4}$\Rightarrow$\romitem{1}.

For \romitem{5}$\Rightarrow$\romitem{6}, the constant case is immediate. Otherwise the parallel tensor $u^*g_{G_{\mathbb R}/K_{\mathbb R}}$ is a positive multiple of the hyperbolic metric. The Lie-triple-system description of totally geodesic submanifolds \cite[Ch.~\upRoman{4}, \S7]{Helgason01} therefore gives an $\mathfrak{sl}_2\mathbb R$-embedding. Its complexification integrates to a homomorphism $\mathrm{SL}_2\mathbb C\to G_{\mathbb C}$; restricting to the connected real group gives $\iota\colon\mathrm{SL}_2\mathbb R\to G_{\mathbb R}^{\mathrm{lin}}$ with $u(gz)=\iota(g)u(z)$. Equivariance implies that $\rho_c(\gamma):=\iota(j(\gamma))^{-1}\rho^{\mathrm{lin}}(\gamma)$ fixes the disk pointwise. This stabilizer is compact, since it fixes a point. And it centralizes the image of $\iota$, since it fixes the tangent Lie triple system. Thus $\rho_c$ is a representation and \eqref{eq:linear-disk-compact} follows.

For \romitem{6}$\Rightarrow$\romitem{5}, the constant case is immediate. Otherwise, let $u_0$ be the totally geodesic orbit map supplied by \romitem{6}. Since $\rho_c$ fixes its image pointwise, $u_0$ is $\rho$-equivariant. Denote its harmonic reduction and Higgs field by $h_0$ and $\Phi_0$. The field $\Phi_0$ is induced from the nilpotent Higgs field of the uniformizing representation, while the compact factor contributes zero. Thus both $\Phi_0$ and $\Phi$ are pointwise nilpotent, so their harmonic reductions are tame and pure imaginary \cite[Definition~A.6]{Mochizuki07}. Apply \cite[Proposition~A.20]{Mochizuki07} to the induced adjoint harmonic bundles. Their harmonic decompositions coincide, and injectivity of the adjoint map on $\mathfrak g_{\mathbb R}$ gives
\[D_h=D_{h_0}, \qquad \Phi=\Phi_0.\]
Since $u_0$ is totally geodesic, $D_{h_0}\Phi_0=0$. Hence $D_h\Phi=0$, and the already established equivalence \romitem{4}$\Leftrightarrow$\romitem{5} shows that $u$ is totally geodesic.
\end{proof}

\begin{corollary}
Under the hypotheses of Theorem~\ref{thm:main-D}, assume that $\Phi$ is pointwise nilpotent and that $G_{\mathbb R}$ is noncompact, almost simple, and not locally isomorphic to $\mathrm{PSL}_2\mathbb R$. If the monodromy is Zariski dense in the sense of Section~\ref{sec:hermitian-transversality}, then $\boldsymbol{\lambda}\neq\mathbf{C}$.
Equivalently, $\lambda_\alpha<C_\alpha$ for at least one $\alpha\in\Delta$.
\end{corollary}

\begin{proof}
Complete equality would place $\rho^{\mathrm{lin}}(\pi_1(X))$ in a compact subgroup or in the product of $\iota(\mathrm{SL}_2\mathbb R)$ and its compact pointwise stabilizer. Under the stated assumptions, both have proper algebraic closures in $G_{\mathbb C}$, contradicting the Zariski density of $\eta(\rho(\pi_1(X)))$.
\end{proof}

\subsection{Reconstruction of VHS when \texorpdfstring{$\boldsymbol{\lambda}=\mathbf{C}$}{lambda=C}}
\label{subsec:hodge-reconstruction}

Assume that the tame harmonic bundle underlies a polarized variation of Hodge structure with nonconstant period map and $\boldsymbol{\lambda}=\mathbf{C}$. Use the factorization \eqref{eq:linear-disk-compact} from Theorem~\ref{thm:main-D}, and write $\mathcal{P}\colon\widetilde X\to\mathcal D$ for the period map, whose projection to $G_{\mathbb R}/K_{\mathbb R}$ is $u$. The $G_\mathbb{R}$-action on the period domain factors through $G_{\mathbb R}^{\mathrm{lin}}$. We use this action for $\iota$ and $\rho_c$ below.

\begin{proposition}
\label{prop:Hodge-compatible-disk}
After identifying $\widetilde X$ with $\mathbb H$,
\[
\mathcal{P}(gz)=\iota(g)\mathcal{P}(z)
\qquad(g\in\mathrm{SL}_2\mathbb R,\ z\in\mathbb H).
\]
Moreover, $\rho_c$ fixes this period-domain disk pointwise.
\end{proposition}

\begin{proof}
Fix $o=\iu\in\mathbb H$ and put $F_0=\mathcal{P}(o)$. Choose $\iota$ as constructed from $u$ in the proof of Theorem~\ref{thm:main-D}, so that $u(gz)=\iota(g)u(z)$.
Writing $E=\Phi_o(\partial_z)$, where $z$ is the standard coordinate on $\mathbb H$, this construction gives
\[
\mathrm{d}\iota(\mathfrak{sl}_2\mathbb C)
=
\mathbb C E
\oplus \mathbb C[E,E^{*_h}]
\oplus \mathbb C E^{*_h}.
\]
Relative to $F_0$, these summands have Hodge types $(-1,1)$, $(0,0)$, and $(1,-1)$, respectively. Moreover,
\[\mathrm{d}\iota(\mathfrak{so}_2\mathbb{R})\otimes_{\mathbb R}\mathbb C=\mathbb C[E,E^{*_h}].\]
Since endomorphisms of type $(0,0)$ preserve the Hodge filtration and $\mathrm{SO}_2\mathbb{R}$ is connected, $\iota(\mathrm{SO}_2\mathbb{R})$ fixes $F_0$. Thus
\[\begin{aligned}
\mathcal{P}_0\colon\mathbb H&\longrightarrow\mathcal{D}\\
go&\longmapsto\iota(g)F_0
\end{aligned}\]
is well-defined. Its differential at $o$ is complex linear, with holomorphic image $\mathbb C E$. Therefore, $\mathcal{P}_0$ is holomorphic and Griffiths-horizontal. Let $\pi\colon\mathcal D\to G_{\mathbb R}/K_{\mathbb R}$ be the natural projection. By equivariance,
\[\pi(\mathcal{P}_0(go))=\iota(g)\pi(F_0)=\iota(g)u(o)=u(go).\]
Hence $\mathcal{P}$ and $\mathcal{P}_0$ lift the same map $u$ and agree at $o$. Now uniqueness of horizontal path lifting gives $\mathcal{P}=\mathcal{P}_0$.
The equivariance of $\mathcal{P}_0$ proves the first assertion.

Finally, monodromy equivariance and \eqref{eq:linear-disk-compact} give
\[\rho_c(\gamma)\mathcal{P}(z)=\iota(j(\gamma))^{-1}\rho^{\mathrm{lin}}(\gamma)\mathcal{P}(z)=\iota(j(\gamma))^{-1}\mathcal{P}(\gamma z)=\mathcal{P}(z).\]
Therefore $\rho_c$ fixes the period-domain disk pointwise.
\end{proof}

For a finite-dimensional real representation $R\colon G_{\mathbb R}\to\mathrm{GL}(W)$, let $\mathcal E_{W,\bullet}$ be the metric extension of the associated harmonic bundle, with holomorphic structure $D_h^{0,1}$. Write $\boldsymbol{\mu}_{\mathrm{HN}}(\mathcal E_{W,\bullet})$ for its parabolic HN slopes, repeated according to rank, and put
\[\mathbf{C}_W(v):=\operatorname{Spec}\left(-\mathrm{d}R(\mathrm{d}u(v))\right),\]
with both lists ordered decreasingly.

\begin{theorem}\label{thm:main-E}
Under these hypotheses, for every $v\in T^1X$,
\[\frac{2}{\deg\KK_{\overline X}(D)}\boldsymbol{\mu}_{\mathrm{HN}}(\mathcal E_{W,\bullet})
=\boldsymbol{\lambda}(W)
=\mathbf{C}_W(v)
=\operatorname{Wt}_\iota(W),\]
where $\operatorname{Wt}_\iota(W)$ is the ordered list of $\mathrm{SL}_2$-weights of $R^{\mathrm{lin}}\circ\iota$, with multiplicities, normalized so that the standard representation has weights $1,-1$.
\end{theorem}

\begin{proof}
Let $U$ be the complexified standard uniformizing local system given by $j\colon\pi_1(X)\to\mathrm{SL}_2\mathbb R$. The standard $\mathrm{SL}_2$-representation theory, e.g. \cite{FultonHarris91}, implies that
\begin{equation}\label{eq:SL2-isotypical}
W_{\mathbb C}=\bigoplus_{m\geqslant0}\operatorname{Sym}^m U\otimes M_m,
\end{equation}
where $\rho_c$ gives compact monodromy on each $M_m$. In our curvature-$(-4)$ normalization, $U$ has Lyapunov and pointwise Cartan spectra $(1,-1)$. Hence the $m$-th summand contributes $(m,m-2,\ldots,-m)$, each repeated $\dim M_m$ times, to both spectra.

Moreover,
\[\mathcal E_{W,\bullet}=\bigoplus_{m\geqslant0}\ \bigoplus_{k=0}^{m}\mathcal L_\bullet^{\otimes(m-2k)}\otimes\mathcal M_{m,\bullet},\]
where $\mathcal L_\bullet^{\otimes2}\simeq\KK_{\overline X}(D)_\bullet$ arising from the lift of $j$, which is equivalent to a spin structure, and $\mathcal M_{m,\bullet}$ is the metric extension of the unitary bundle $M_m$. Each summand is polystable of slope $(m-2k)\deg\KK_{\overline X}(D)/2$. Ordering these slopes proves that the normalized HN spectrum is the same list of weights.
\end{proof}

\subsection{Maximal Hodge bundles over closed surfaces}\label{subsec:maximal-hodge-comparison}

Assume that $X=\overline{X}$ is compact of genus $g\geqslant2$, and let $\rho\colon\pi_1(X)\to G_{\mathbb R}$ be the monodromy of a polarized variation of Hodge structure with nonconstant period map. As in Section~\ref{subsec:hodge-reconstruction}, $G_{\mathbb R}$ is a connected semisimple real Lie group with finite center and without compact factors, and we use the harmonic Hodge metric. We apply the group-level rigidity statements below in $G_{\mathbb R}^{\mathrm{lin}}$.

Recall that a Hodge bundle is \emph{maximal} when
\[-\tau=r_T(2g-2),\]
where $\tau$ and $r_T$ are the Toledo invariant and the Toledo rank of its Hodge grading, with the sign and normalization of \cite{BCGT23}. The inequality defining maximality is \cite[Corollary~5.5]{BCGT23}. Here we use the corresponding rigidity theorem to determine the spectra.

By \cite[Theorems~6.5 \&~6.17]{BCGT23}, the monodromy of a maximal Hodge bundle has the factorization \eqref{eq:linear-disk-compact}, with $\iota$ horizontal and $\rho_c$ fixing the period-domain disk pointwise. The period map is a totally geodesic embedding whose holomorphic sectional curvature is maximal among horizontal directions for the normalized Hodge metric of \cite[\S3.3 \& Corollary~6.19]{BCGT23}.

\begin{corollary}\label{cor:maximal-tannakian}
	For a maximal Hodge bundle on $X$ as above, $\boldsymbol{\lambda}=\mathbf{C}$. Moreover, for every finite-dimensional real representation $R\colon G_{\mathbb R}\to\mathrm{GL}(W)$ and every $v\in T^1X$,
	\[\frac1{g-1}\boldsymbol{\mu}_{\mathrm{HN}}(\mathcal{E}_W)=\boldsymbol{\lambda}(W)=\mathbf{C}_W(v)=\operatorname{Wt}_\iota(W).\]
	Here $\mathcal{E}_W$ is the holomorphic bundle underlying the associated harmonic bundle, and the weights of the horizontal $\mathrm{SL}_2$-factor are normalized so that its standard representation has weights $1,-1$. Thus the normalized HN, Lyapunov, and pointwise Cartan polygons coincide.
\end{corollary}

\begin{proof}
	The factorization above gives a totally geodesic harmonic map to $G_{\mathbb R}/K_{\mathbb R}$. Since the Higgs field of a Hodge bundle is nilpotent, Theorem~\ref{thm:main-D} gives $\boldsymbol{\lambda}=\mathbf{C}$. Theorem~\ref{thm:main-E} then gives the spectral identity, using $D=\varnothing$ and $\deg\KK_X=2g-2$.
\end{proof}

Explicitly, each $\operatorname{Sym}^m U\otimes M_m$ in \eqref{eq:SL2-isotypical} contributes the weights $(m,m-2,\ldots,-m)$, each repeated $\dim M_m$ times. Its HN slopes are these weights multiplied by $g-1$.

Among the Hodge bundles satisfying complete upper equality, the maximal ones are precisely those whose horizontal disk has maximal holomorphic sectional curvature, by \cite[Corollary~6.19]{BCGT23} and Proposition~\ref{prop:Hodge-compatible-disk}.

\bibliographystyle{alpha}
\bibliography{LHN}

\end{document}